\documentclass{amsart}
\usepackage{amsmath}
\usepackage{amssymb}
\usepackage{extarrows}
\usepackage{mathabx}
\usepackage[all]{xy}
\usepackage{amsbsy}
\usepackage{graphicx}
\usepackage{appendix}
\usepackage{float}
\usepackage{subfigure}
\usepackage[numbers,sort&compress]{natbib}
\usepackage{amsthm}
\usepackage{geometry}
\usepackage{fancyhdr}
\usepackage{color} 
\usepackage{overpic}
\usepackage{enumerate}
\usepackage{mathrsfs}
\usepackage{colonequals}
\usepackage{microtype}
\usepackage{hyperref}
\usepackage{diagbox}

\newtheorem{thm}{Theorem}[section]
\newtheorem{cor}[thm]{Corollary}
\newtheorem{prop}[thm]{Proposition}
\newtheorem{lem}[thm]{Lemma}

\theoremstyle{definition}
\newtheorem{defn}[thm]{Definition}
\newtheorem{cons}[thm]{Construction}
\newtheorem{exmp}[thm]{Example}
\newtheorem{conj}[thm]{Conjecture}
\newtheorem*{fact}{Fact}

\newtheorem*{conv}{Convention}

\newtheorem*{org}{Organization}
\newtheorem*{ack}{Acknowledgement}

\theoremstyle{remark}
\newtheorem{rem}[thm]{Remark}

\numberwithin{equation}{section}

\newcommand{\beq}{\begin{equation*}\begin{aligned}}
\newcommand{\eeq}{\end{aligned}\end{equation*}}
\newcommand{\bpf}{\begin{proof}}
\newcommand{\epf}{\end{proof}}
\newcommand{\bthm}{\begin{thm}}
\newcommand{\ethm}{\end{thm}}
\newcommand{\bprop}{\begin{prop}}
\newcommand{\eprop}{\end{prop}}
\newcommand{\bcor}{\begin{cor}}
\newcommand{\ecor}{\end{cor}}
\newcommand{\blem}{\begin{lem}}
\newcommand{\elem}{\end{lem}}
\newcommand{\bdefn}{\begin{defn}}
\newcommand{\edefn}{\end{defn}}
\newcommand{\bcons}{\begin{cons}}
\newcommand{\econs}{\end{cons}}
\newcommand{\bexmp}{\begin{exmp}}
\newcommand{\eexmp}{\end{exmp}}
\newcommand{\brem}{\begin{rem}}
\newcommand{\erem}{\end{rem}}
\newcommand{\bfa}{\begin{fact}}
\newcommand{\efa}{\end{fact}}
\newcommand{\benu}{\begin{enumerate}[(1)]}
\newcommand{\eenu}{\end{enumerate}}

\newcommand{\bdia}{\begin{displaymath}\xymatrix}
\newcommand{\edia}{\end{displaymath}}

\newcommand{\shi}{\underline{\rm SHI}}

\newcommand{\khii}{\underline{\rm KHI}}

\newcommand{\deq}{\colonequals}

\newcommand{\cstar}{\mathbb{C}\backslash\{0\}}

\newcommand{\dimc}{\dim_\mathbb{C}}

\newcommand{\al}{\alpha}

\newcommand{\ga}{\gamma}
\newcommand{\Ga}{\Gamma}

\newcommand{\op}{\oplus}

\newcommand{\p}{\prime}
\newcommand{\pp}{{\prime\prime}}

\newcommand{\aand}{~{\rm and}~}

\newcommand{\intg}{\mathbb{Z}}
\newcommand{\ft}{{\mathbb{F}_2}}

\newcommand{\ra}{\rightarrow}

\newcommand{\xra}{\xrightarrow}

\DeclareMathOperator{\cok}{coker}
\DeclareMathOperator{\im}{Im}
\DeclareMathOperator{\ke}{ker}

\DeclareMathOperator{\cone}{Cone}

\newcommand{\sut}[1]{\mathbf{\Gamma}_{#1}}
\newcommand{\sutg}[2]{(\mathbf{\Gamma}_{#1},#2)}
\newcommand{\psp}[2]{\psi^{#1}_{+,#2}}
\newcommand{\psm}[2]{\psi^{#1}_{-,#2}}
\newcommand{\Psp}[2]{\Psi^{#1}_{+,#2}}
\newcommand{\Psm}[2]{\Psi^{#1}_{-,#2}}
\newcommand{\dehny}[1]{\mathbf{Y}_{#1}}

\begin{document}

\title{Knot surgery formulae for instanton Floer homology I: the main theorem}


\author{Zhenkun Li}
\address{Department of Mathematics and Statistics, University of South Florida}
\curraddr{}
\email{zhenkun@usf.edu}
\thanks{}

\author{Fan Ye}
\address{Department of Mathematics, Harvard University}
\curraddr{}
\email{fanye@math.harvard.edu}
\thanks{}

\date{}
\dedicatory{}

\begin{abstract}
We prove an integral surgery formula for framed instanton homology $I^\sharp(Y_m(K))$ for any knot $K$ in a $3$-manifold $Y$ with $[K]=0\in H_1(Y;\mathbb{Q})$ and $m\neq 0$. Although the statement is similar to Ozsv\'ath-Szab\'o's integral surgery formula for Heegaard Floer homology, the proof is new and based on sutured instanton homology $SHI$ and the octahedral lemma in the derived category. As byproducts, we obtain a formula computing instanton knot homology of the dual knot analogous to Eftekhary's and Hedden-Levine's work, and also an exact triangle between $I^\sharp(Y_m(K))$, $I^\sharp(Y_{m+k}(K))$ and $k$ copies of $I^\sharp(Y)$ for any $m\neq 0$ and large $k$. In the proof of the formula, we discover many new exact triangles for sutured instanton homology and relate some surgery cobordism map to the sum of bypass maps, which are of independent interest. In a companion paper, we derive many applications and computations based on the integral surgery formula.
\end{abstract}
\maketitle
\tableofcontents

\section{Introduction}\label{sec: introduction}


The framed instanton homology $I^\sharp(Y)$ for a closed 3-manifold $Y$ was introduced by Kronheimer and Mrowka in \cite{kronheimer2011knot} and has been conjectured to be isomorphic to the hat version of Heegaard Floer homology $\widehat{HF}(Y)$. This conjecture is still widely open and, due to the computational difficulty of instanton Floer homology, not many examples has been known. In recent years, many people have done computations of the framed instanton homology of special families of $3$-manifolds, see for example \cite{lidman2020framed,baldwin2019lspace,baldwin2020concordance}. Yet, most results have focused on computing the dimension of framed instanton Floer homology and many techniques only work for $S^3$ or rational homology spheres, however, a general structural theorem that relates the framed instanton homology of Dehn surgeries to the information from the knot complement still remains elusive. 

In \cite{LY2021large}, the authors of the current paper proved a large surgery formula for framed instanton homology which led to a series of applications in computing the framed instanton homology and studying the representations of the fundamental groups of Dehn surgeries of some families of knots. However, in that work, the Dehn surgery slope must be large (at least $2g+1$ where $g$ is the Seifert genus of the knot), and thus still not much is known about the framed instanton homology of small Dehn surgery slopes. In this paper, we further prove an integral surgery formula for rationally null-homologous knots, inspired by Ozsv\'ath-Szab\'o's surgery formula for Heegaard Floer homology \cite{Ozsvath2008integral,Ozsvath2011rational}. For simplicity, in the introduction, we present only the discussions and results for (integral) null-homologous knots (\textit{e.g.} knots in $S^3$) and leave the general setup to Section \ref{subsec: integral surgery formula for I-sharp}.

First, let us recall the results from \cite{LY2021large}. Suppose $K\subset Y$ is a null-homologous knot. Let $Y\backslash N(K)$ be the knot complement, and let $\Ga_{\mu}$ be the union of two oppositely oriented meridians of the knot on $\partial (Y\backslash N(K))$. Let $SHI(-Y\backslash N(K),-\Ga_{\mu})$ be the corresponding sutured instanton homology introduced by Kronheimer-Mrowka \cite{kronheimer2010knots}, where the minus sign denotes orientation reversal for technical needs (note that $SHI(-M,-\ga)\cong SHI(M,\ga)$ and in particular $I^\sharp (-Y_{-m}(K))\cong I^\sharp(Y_{-m}(K))$). A Seifert surface of $K$ induces a $\intg$-grading on $SHI(-Y\backslash N(K),-\Ga_{\mu})$. In \cite{LY2021large}, we constructed a set of differentials on $SHI(-Y\backslash N(K),-\Ga_{\mu})$
$$d^i_j:SHI(-Y\backslash N(K),-\Ga_{\mu},i)\ra SHI(-Y\backslash N(K),-\Ga_{\mu},j)$$
for any gradings $i\neq j\in\intg$. We then constructed bent complexes
$$A_s=\bigg(SHI(-Y\backslash N(K),-\Ga_{\mu}),\sum_{s\leq i<j}d^i_j+\sum_{s\geq i>j}d^i_j\bigg),$$
$$B^+=\bigg(SHI(-Y\backslash N(K),-\Ga_{\mu}),\sum_{i<j}d^i_j\bigg),~{\rm and}~ B^-=\bigg(SHI(-Y\backslash N(K),-\Ga_{\mu}),\sum_{i>j}d^i_j\bigg).$$

From \cite{LY2021large}, the homologies of these complexes are related to the Dehn surgeries of $K$ as follows:
\begin{equation}\label{eq: intro, homology of B_s}
	H(B^+)\cong H(B^-)\cong I^{\sharp}(-Y),
\end{equation}
\begin{equation}\label{eq: intro, homology of A_s}
	I^{\sharp}(Y_{-m}(K))\cong\bigoplus_{s=\lfloor\frac{1-m}{2}\rfloor}^{\lfloor\frac{m-1}{2}\rfloor}H(A_s) \text{ for any integer }m\ge 2g(K)+1.
\end{equation}

To state the integral surgery formula, we introduce more notations. For $s\in\intg$, let $B^{\pm}_s$ be identical copies of $B^{\pm}$. Define chain maps
$$\pi^{\pm,s}:A_s\ra B^\pm_{s}$$
as follows. For $x\in (SHI(-Y\backslash N(K),-\Ga_{\mu}),i)$,
\[\pi^{+,s}(x)=\begin{cases}
x&i\ge s,\\
0&i< s,\end{cases}\aand \pi^{-,s}(x)=\begin{cases}
x&i\le s,\\
0&i>s.\end{cases}\]

Let $\pi^\pm$ denote the direct sum of all $\pi^{\pm,s}$. While this slightly abuses the notation, we also use them to denote the induced maps on the homologies. The main result of the paper is the following.

\bthm[Integral surgery formula]\label{thm: intro, 1}
Suppose $K\subset Y$ is a null-homologous knot. Let $A_s$, $B_s^{\pm}$, and $\pi^{\pm}$ be defined as above. Then for any $m\in \mathbb{Z}\backslash\{ 0\}$, there exists an isomorphism $$\Xi_{m}:\bigoplus_{s\in\mathbb{Z}}H(B^+_s)\xra{\cong}\bigoplus_{s\in\mathbb{Z}}H(B^-_{s+m})$$as the direct sum of isomorphisms$$\Xi_{m,s}:H(B^+_s)\xra{\cong}H(B^-_{s+m})$$ such that
$$I^{\sharp}(-Y_{-m}(K))\cong H\bigg({\rm Cone}(\pi^{-}+\Xi_m\circ \pi^{+}:\bigoplus_{s\in\intg}H(A_s)\to\bigoplus_{s\in\intg}H(B^-_s))\bigg).$$
\ethm
\brem
As an analog of the surgery formula in Heegaard Floer homology, the map $\pi^-$ is related to the vertical projection map $v$, and the map $\Xi_m\circ\pi^+$ is related to the map $h$, which is defined by the composition of the horizontal projection and a chain homotopy equivalence between $C\{j\ge 0\}$ to $C\{i\ge 0\}$. Here we bend the horizontal part of the hook complex to become vertical, so the differentials go upwards. The homotopy equivalence in Heegaard Floer homology depends on many auxiliary choices (\textit{c.f.} Construction before \cite[Lemma 2.16]{hedden21surgery}). The same situation applies to $\Xi_m$. Hence we only state the existence of the isomorphism.
\erem
\brem\label{rem: scalar ambiguity}
The hypothesis of Theorem \ref{thm: intro, 1} excludes the case where $m=0$. This is due to the sign ambiguity in the definition of sutured instanton homology. The original version of sutured instanton homology defined by Kronheimer-Mrowka \cite{kronheimer2010knots} was only well-defined up to isomorphisms, and then Baldwin-Sivek \cite{baldwin2015naturality} proved that they are well-defined up to a scalar in $\mathbb{C}$. As a result, all related maps are only well-defined up to scalars. When $m\neq 0$, the maps $\pi^{+,s}$ and $\Xi_m\circ \pi^{-,s}$ have distinct target spaces, namely $B_s$ and $B_{s+m}$. As a result, the scalar ambiguity for individual maps does not influence the dimension of the homology of the mapping cone. However, when $m=0$, different scalars would indeed make differences. For an example of this subtlety, see the end of Section \ref{sec: (-1)-surgery and bypass}.
\erem
\brem
We also obtain a formula computing instanton knot homology $KHI(-Y_{-m}(K),K_{-m})$ of the dual knot $K_{-m}$ inside the resulting manifold in Theorem \ref{thm: mapping cone for KHI}, which is analogous to the results by Eftekhary for knot Floer homology $\widehat{HFK}$ \cite[Proposition 1.5]{eft18bordered} and by Hedden-Levine \cite{hedden21surgery}. In this formula, we may assume $m=0$ because the scalar issue in Remark \ref{rem: scalar ambiguity} does not appear.
\erem
With the isomorphisms in (\ref{eq: intro, homology of A_s}) and (\ref{eq: intro, homology of B_s}), we can truncate the above formula for $I^{\sharp}(Y_{-m}(K))$ to obtain the following exact triangle.
\bcor[Generalized surgery exact triangle]\label{cor: intro, 1}
Suppose $K\subset Y$ is a null-homologous knot, and $m$ is a fixed non-zero integer. Then for any sufficiently large integer $k$, there exists an exact triangle
\begin{equation}\label{eq: intro, 2}
	\xymatrix{
	I^{\sharp}(-Y_{-m-k}(K))\ar[rr]&&\mathop{\bigoplus}_{i=1}^{k}I^{\sharp}(-Y)\ar[dl]\\
	&I^{\sharp}(-Y_{-m}(K))\ar[ul]&
	}
\end{equation}
\ecor

\brem
The analogous result of the exact triangle (\ref{eq: intro, 2}) in Heegaard Floer theory was proved by Ozsv\'ath-Szab\'o \cite{Ozsvath2008integral} using twisted coefficients, which is a crucial step towards proving the integral surgery formula in their setup. The proof cannot be applied to instanton theory directly. Thus, in this paper, we adopt a reversed approach: we use sutured instanton theory to prove Theorem \ref{thm: intro, 1} and derive Corollary \ref{cor: intro, 1} as a direct application. The strategy to prove Theorem \ref{thm: intro, 1} can be found in Section \ref{subsec: A strategy to prove the formula} and Section \ref{subsec: a strategy}.
\erem

The analogs of $\pi^{\pm,s}$ in Heegaard Floer theory can be interpreted as cobordism maps associated to some particular spin${}^c$ structures. In instanton theory, there is a decomposition of cobordism maps along basic classes. However, currently such a decomposition is only known to exist for cobordisms whose first Betti number is zero. So for the moment let us assume the ambient $3$-manifold $Y$ is a rational homology sphere. For any integer $m$, there is a natural cobordism $W_m$ from $-Y^3_{-m}(K)$ to $-Y^3$. From \cite[Section 1.2]{baldwin2019lspace}, there exists a decomposition of the cobordism map $I^{\sharp}(W_m)$ along basic classes
	$$I^{\sharp}(W_m)=\sum_{s\in\intg}I^{\sharp}(W_m,[s]),$$
	where $[s]\in H^2(W)$ denotes the class that satisfies the equality
	$$[s]([\widebar{S}])=2s-m.$$
We make the following conjecture.

\begin{conj}\label{conj: intro, 1}
Suppose $K\subset Y$ is a null-homologous knot. Suppose $b_1(Y)=0$ and $m\in\intg$ with $m\geq 2g(K)+1$. Let $A_s$, $B^-_s$, $\pi^{\pm,s}$, $W_m$, $I^{\sharp}(W_m,[s])$ be defined as above. Then for any $s\in[\lfloor\frac{m-1}{2}\rfloor,\lfloor\frac{1-m}{2}\rfloor]\cap \intg$, there are commutative diagrams
\begin{equation*}
	\xymatrix@R=6ex{
	H(A_s)\ar@{^{(}->}[rr]\ar[d]^{\pi^{-,s}}&&I^{\sharp}(-Y_{-m}(K))\ar[d]^{I^{\sharp}(W_m,[s])}\\
	H(B_s)\ar[rr]^{\cong}&&I^{\sharp}(-Y)
	}
	\xymatrix@R=6ex{
	H(A_s)\ar@{^{(}->}[rr]\ar[d]^{\pi^{+,s}}&&I^{\sharp}(-Y_{-m}(K))\ar[d]^{I^{\sharp}(W_m,[s+m])}\\
	H(B_s)\ar[rr]^{\cong}&&I^{\sharp}(-Y)
	}
\end{equation*}
\end{conj}

\begin{rem}
	In Heegaard Floer theory, the large surgery formula in \cite[Theorem 4.1]{ozsvath2004holomorphicknot} states that the homology of the bent complex $A_s$ is isomorphic to the Heegaard Floer homology of $Y_{-m}(K)$ together with a spin${}^c$ structure specified by $s$. In instanton theory, we do not have the spin${}^c$ structures in the construction of instanton Floer homology but a similar decomposition was introduced in \cite{LY2020,LY2021enhanced}. However, involving the spin${}^c$-type decomposition in the statement of Conjecture \ref{conj: intro, 1} would make the statement more complicated. So here we only write the top horizontal map in each commutative diagram as an inclusion. 
\end{rem}

The obstacle to obtaining a decomposition of the instanton cobordism map, in general, is one of the difficulties in exporting the original proof of the integral surgery formula in Heegaard Floer theory to the instanton setup. To overcome this problem, we need to work with a suitable setup for which some kind of decompositions do exist. A good candidate is sutured instanton theory. In sutured instanton theory, properly embedded surfaces induce $\intg$-gradings on the homology, and bypass maps relating different sutures are homogeneous with respect to such gradings. We have already used this setup to construct spin${}^c$-like decompositions for the framed instanton homology of Dehn surgeries of knots, construct bent complexes in instanton theory, and establish a large surgery formula in our previous work \cite{LY2020,LY2021enhanced,LY2021large}.

In this paper, to prove the integral surgery formula, we further study the relations between different sutures on the knot complement and establish some new exact triangles and commutative diagrams that may be of independent interest. Then these new and old algebraic structures relating to different sutures enable us to apply the octahedral lemma to prove the desired integral surgery formula. It is worth mentioning that ultimately the whole proof in the current paper depends only on some most fundamental properties of Floer theory: the surgery exact triangle, the functoriality of the cobordism maps, and the adjunction inequality. This implies that the existence of the surgery formula is an inherent property of Floer theory.


The surgery formula developed in the current paper is a powerful tool to study the Dehn surgeries along knots. It enables us to do explicit computations in many cases, even when the ambient $3$-manifold has a positive first Betti number. In a companion paper \cite{LY2022integral2}, we will use the surgery formula and the techniques developed in this paper to derive many new applications and computations. We sketch the results as follows.
\begin{enumerate}
    \item We study the behavior of the integral surgery formula under the connected sum with a core knot in a lens space (whose complement is a solid torus) and then derive a rational surgery formula for framed instanton homology.
    \item We study the $0$-surgery on a knot $K$ inside $S^3$. We derive a formula computing the non-zero grading part of $I^\sharp(S_0(K))$ with respect to the grading induced by the Seifert surface.
    \item We study non-zero integral surgeries on Boromean knots inside $Y=\#^{2g}S^1\times S^2$, which gives  nontrivial circle bundles over surfaces. In this case the bent complexes $A_s$ and $B^{\pm}_s$ can be computed directly and the maps $\pi^{\pm}$ between them can be fixed with the help of the $\Lambda^* H_1(Y;\mathbb{C})$-action. Moreover, we show $\dim_\mathbb{C} I^\sharp(Y)=\dim_{\ft}\widehat{HF}(Y)$ for most Seifert fibered manifolds $Y$ with non-zero orbifold degrees.
   
    \item We study a family of alternating knots. Using an inductive argument by the oriented skein relation, we can describe their bent complexes explicitly and then the surgery formula applies routinely. 
    
     \item Using the same technique as above, we also study the non-zero integral surgery of twisted Whitehead doubles. The results for Whitehead doubles can also tell us the framed instanton Floer homology of the splicing of two knot complements in $S^3$, where one knot is a twist knot, \textit{i.e.} the Whitehead double of the unknot.

     \item We study almost L-space knots, \textit{i.e.}, a non-L-space knot $K$ such that there exists $n\in\mathbb{N}_+$ with $\dim I^{\sharp}(S^3_n(K))=n+2$ (see \cite{baldwin22knot52} for the results in Heegaard Floer theory). We prove a genus one almost $L$-space knot is either the figure-eight or the mirror of the knot $5_2$. We also show that almost $L$-space knots of genus at least $2$ are fibered and strongly quasi-positive.
\end{enumerate}
\begin{org}
The paper is organized as follows. In Section \ref{sec: basic setups}, we introduce basic setup, the notations in sutured instanton homology, and deal with the scalar ambiguity mentioned in Remark \ref{rem: ambiguity}. We also present some algebraic lemmas including the octahedral lemma in the derived category that are used in latter sections. In Section \ref{sec: integral surgery formulae}, we present the strategy to prove the integral surgery formula. We first restate the integral surgery formula using sutured instanton homology, and explain how to apply the octahedral lemma to prove it. Then we explain how to translate the integral surgery formula from the language of sutured instanton theory to the language of bent complexes, which coincides with the discussions in the introduction. All the rest of the sections are devoted to prove the three exact triangles and three commutative diagrams that are involved in the octahedral lemma, {\it i.e.}, Equation (\ref{eq: intro, exact, 1}) to Equation (\ref{eq: intro, commute, 3}). In Section \ref{sec: (-1)-surgery and bypass}, we study the relation between the $(-1)$-Dehn surgery map associated to a curve intersecting the suture twice and the two natural bypass maps associated to that curve. This helps us to prove Equation (\ref{eq: intro, exact, 1}) and Equation (\ref{eq: intro, commute, 1}). In Section \ref{sec: Surgery exact triangles and commutative diagrams}, Equation (\ref{eq: intro, exact, 2}), Equation (\ref{eq: intro, commute, 2}) and part of Equation (\ref{eq: intro, exact, 3}) are proved. The last two sections are devoted to prove Equation (\ref{eq: intro, exact, 3}) and Equation (\ref{eq: intro, commute, 3}), which is the most technical part of the paper. In Section \ref{sec: Some technical lemmas} we prove some technical lemmas that are finally used in Section \ref{sec: third triangle} to finish the proof.
\end{org}

\begin{ack}
The authors thank John A. Baldwin and Steven Sivek for the discussion on the proof of Proposition \ref{prop: -1 surgery and bypass}, and thank Zekun Chen and Linsheng Wang for the discussion on homological algebra. The authors would like to thank Ciprian Manolescu and Jacob Rasmussen, and the anonymous referee for helpful comments. The authors also thank Sudipta Ghosh, Jianfeng Lin, Yi Xie and Ian Zemke for valuable discussions. The second author is also grateful to Yi Liu for inviting him to BICMR, Peking University when he was writing the early version of this paper.
\end{ack}

\section{Basic setup}\label{sec: basic setups}
\subsection{Conventions}

If it is not mentioned, all manifolds are smooth, oriented, and connected. Homology groups and cohomology groups are with $\mathbb{Z}$ coefficients. We write $\mathbb{Z}_n$ for $\mathbb{Z}/n\mathbb{Z}$ and $\mathbb{F}_2$ for the field with two elements.

A knot $K\subset Y$ is called \textbf{null-homologous} if it represents the trivial homology class in $H_1(Y;\mathbb{Z})$, while it is called \textbf{rationally null-homologous} if it represents the trivial homology class in $H_1(Y;\mathbb{Q})$.

For any oriented 3-manifold $M$, we write $-M$ for the manifold obtained from $M$ by reversing the orientation. For any surface $S$ in $M$ and any suture $\ga\subset \partial M$, we write $S$ and $\ga$ for the same surface and suture in $-M$, without reversing their orientations. For a knot $K$ in a 3-manifold $Y$, we write $(-Y,K)$ for the induced knot in $-Y$ with induced orientation, called the \textbf{mirror knot} of $K$. The corresponding balanced sutured manifold is $(-Y\backslash N(K),-\ga_K)$.

\subsection{Sutured instanton homology}\label{sec: preliminaries}

For any \textbf{balanced sutured manifold} $(M,\ga)$ \cite[Definition 2.2]{juhasz2006holomorphic}, Kronheimer-Mrowka \cite[Section 7]{kronheimer2010knots} constructed an isomorphism class of $\mathbb{C}$-vector spaces $SHI(M,\ga)$. Later, Baldwin-Sivek \cite[Section 9]{baldwin2015naturality} dealt with the naturality issue and constructed (untwisted and twisted vesions of) projectively transitive systems related to $SHI(M,\ga)$. We will use the twisted version, which we write as $\shi(M,\ga)$ and call \textbf{sutured instanton homology}. 

In this paper, when considering maps between sutured instanton homology, we can regard them as linear maps between actual vector spaces, at the cost that equations (or commutative diagrams) between maps only hold up to a non-zero scalar due to the projectivity. A more detailed discussion on the projectivity can be found in the next subsection. 

Moreover, there is a relative $\mathbb{Z}_2$-grading on $\shi(M,\ga)$ obtained from the construction of sutured instanton homology, which we consider as a \textbf{homological grading} and use to take Euler characteristic.

\bdefn
Suppose $K$ is a knot in a closed 3-manifold $Y$. Let $Y(1)\deq Y\backslash B^3$ and let $\delta$ be a simple closed curve on $\partial Y(1)\cong S^2$. Let $Y\backslash N(K)$ be the knot complement and let $\Ga_{\mu}$ be two oppositely oriented meridians of $K$ on $\partial (Y\backslash N(K))\cong T^2$. Define\[I^\sharp(Y)\deq \shi(Y(1),\delta)\aand \khii(Y,K)\deq \shi(Y\backslash N(K),\Ga_{\mu}).\]
\edefn
\brem\label{rem: basept}
By the naturality results, we should specify the places of the removing ball, the neighborhood of the knot, and the sutures to define $I^\sharp(Y)$ and $\khii(Y,K)$. These data can be fixed by choosing a basepoint in $Y$ or $K$. For simplicity, we omit those choices in the notations.
\erem
From now on, we will suppose $K\subset Y$ is a rationally null-homologous knot and fix some notations. Let $\mu$ be the meridian of $K$ and pick a longitude $\lambda$ (such that $\lambda\cdot \mu=1$) to fix a framing of $K$.  We will always assume $Y\backslash N(K)$ is irreducible, but many results still hold due to the following connected sum formula of sutured instanton homology \cite[Remark 1.6]{li2018contact}:
$$\shi(Y^\p\# Y\backslash N(K),\ga)\cong I^\sharp(Y^\p)\otimes \shi(Y\backslash N(K),\ga).$$

Given coprime integers $r$ and $s$, let $\Gamma_{r/s}$ be the suture on $\partial (Y\backslash N(K))$ consists of two oppositely oriented, simple closed curves of slope $-r/s$, with respect to the chosen framing (\textit{i.e.} the homology of the curves are $\pm(-r\mu+s\lambda)\in H_1(\partial N(K))$). Pick $S$ to be a minimal genus Seifert surface of $K$, with genus $g=g(S)$. Note that $\partial S$ may have multiple components. 

\begin{conv}
For a fixed pair $(\lambda,\mu)$ as above, we write $p=\partial S\cdot \mu$ and $q=\partial S\cdot \lambda$. Note that when an orientation of the knot $K$ is chosen, the orientation of $S$ is induced by the knot. The orientation of $\mu$ is chosen such that $p>0$ and the orientation of $\lambda$ is chosen such that $\lambda\cdot \mu=1$. Note that $p$ is the order of $[K]\in H_1(Y)$, \textit{i.e.}, $p$ is the minimal positive integer satisfying $p [K]=0\in H_1(Y)$. When $K$ is null-homologous, we always choose the Seifert framing $\lambda=\partial S$. In such case, we have $(p,q)=(1,0)$. 
\end{conv}
\brem\label{rem: notations}
The meanings of $p$ and $q$ above are different from our previous papers \cite{LY2020,LY2021large}. Before, we used $\hat{\mu}$ and $\hat{\lambda}$ to denote the meridian of the knot $K$ and the preferred framing. In particular, the framing is fixed by \cite[Definition 4.2]{LY2020}. Note that in that case, we assume that $\partial S$ is connected, and hence it is the same as the homological longitude (with the notation $\lambda$ in previous papers, while we use $\mu$ to denote the homological meridian). Also, the numbers $p$ and $q$ in this paper should be $q$ and $q_0$ in previous papers.
\erem

For simplicity, we use the bold symbol of the suture to represent the sutured instanton homology of the balanced sutured manifold with the reversed orientation:
$$\mathbf{\Gamma}_{r/s}\deq \shi(-(Y\backslash N(K)),-\Gamma_{r/s}).$$
When $(r,s)=(\pm 1,0)$, we write $\Ga_{r/s}=\Ga_{\mu}$. When $s=\pm 1$, we write $\Ga_{n}=\Ga_{n/1}=\Ga_{(-n)/(-1)}$. We also write $\mathbf{\Gamma}_{\mu}$ and $\mathbf{\Gamma}_{n}$ for the corresponding sutured instanton homologies.
\brem\label{rem: ambiguity}
Strictly speaking, the sutures corresponding to $(r,s)=(1,0)$ and $(-1,0)$ are not identical because the orientations are opposite. Since both sutures are on $\partial (Y\backslash N(K))$ of the same slope, they are isotopic. Moreover, we can choose a canonical isotopy by rotating the suture along the direction specified by the orientation of the knot. Due to discussion in Heegaard Floer theory \cite{sarkar15moving,Zemke2019} and the conjectural relation between Heegaard Floer theory and instanton theory \cite{kronheimer2010knots}, it is expected that rotating the suture back to the original position induces a nontrivial isomorphism of the sutured instanton homology. So we pick the canonical isotopy to be the minimal rotation of the suture. Hence we can abuse notations and write $\Ga_\mu$ for both sutures. The same discussion also applies to the relation between $\Ga_{n/1}$ and $\Ga_{(-n)/(-1)}$.
\erem

We always assume that $\partial S$ has minimal intersections with $\Ga_{r/s}$, i.e. $|\partial S\cap \ga|=2|rp-sq|$. When the intersection number $rp-sq$ is odd, then $S$ induces a $\intg$-grading on $\sut{r/s}$. When $rp-sq$ is even, we need to perform either a positive stabilization or negative stabilization on $S$ to induce a $\intg$-grading, and the two gradings are related by an overall grading shift of $1$. To get rid of stabilizations, we can equivalently regard that, in this case, the surface $S$ induces a $(\intg+\frac{1}{2})$-grading. We write the graded part of $\sut{r/s}$ as
$$\sutg{r/s}{i}\deq \shi(-(Y\backslash N(K)),-\Gamma_{r/s},S,i)$$
with $i\in\intg$ or $i\in\intg+\frac{1}{2}$, depending on the parity of the intersection number. From the construction of the grading in \cite{li2019direct}, we have the following vanishing theorem due to the adjunction inequality.
\blem\label{lem: vanishing grading}
We have $\sutg{r/s}{i}=0$ when $$|i|>\frac{|rp-sq|-\chi(S)}{2}.$$
\elem
\bpf
This follows from \cite[Theorem 2.21 (1)]{LY2020} (which is ultimately based on \cite[Proposition 7.5]{kronheimer2010knots}) and a direct computation in the new notations.
\epf

The bypass exact triangle for sutured instanton homology was introduced by Baldwin-Sivek in \cite[Section 4]{baldwin2018khovanov}. In \cite[Section 4.2]{LY2020}, we applied the triangle to sutures on knot complements and computed the grading shifts. We restate the results in the notation introduced before.
\blem\label{lem: bypass n,n+1,mu}
For any $n\in\intg$, there are two graded bypass exact triangles
\begin{equation*}
\xymatrix@R6ex{
\sutg{n}{i+\frac{p}{2}}\ar[rr]^{\psp{n}{n+1}}&&\sutg{n+1}{i}\ar[dl]^{\psp{n+1}{\mu}}\\
&\sutg{\mu}{i-\frac{np-q}{2}}\ar[ul]^{\psp{\mu}{n}}&
}	
\end{equation*}
\begin{equation*}
\xymatrix@R6ex{
\sutg{n}{i-\frac{p}{2}}\ar[rr]^{\psm{n}{n+1}}&&\sutg{n+1}{i}\ar[dl]^{\psm{n+1}{\mu}}\\
&\sutg{\mu}{i+\frac{np-q}{2}}\ar[ul]^{\psm{\mu}{n}}&
}	
\end{equation*}
where the maps are homogeneous with respect to the homological $\mathbb{Z}_2$-gradings.  
\elem
\bpf
This is \cite[Proposition 4.14]{LY2020} in the new notations. The idea of the proof can be found in \cite[Lemma 3.18]{LY2020} (see also \cite[Remark 3.19]{LY2020}). Roughly, we perturb the surface $S$ by stabilizations so that its boundary is disjoint from the bypass arc. Then the grading shifts are obtained by counting the number of positive or negative stabilizations. 

Unlike the setup in \cite[Section 4]{LY2020}, here $K$ is not necessarily a dual knot of the Dehn surgery on a null-homologous knot, so we adopt the remarks in the beginning of \cite[Section 5]{LY2020}. For example, when $n$ is large enough so that $np-q\ge 0$ and \cite[Proposition 4.14 (1)]{LY2020} applies, we have \[\hat{i}_{\max}^n=\frac{np-q-\chi(S)}{2},~\hat{i}_{\min}^n=-\frac{np-q-\chi(S)}{2},~\hat{i}_{\max}^\mu=\frac{p-\chi(S)}{2},~\hat{i}_{\min}^\mu=-\frac{p-\chi(S)}{2},\]where we omit $\lceil\cdot\rceil$ since we think about $\mathbb{Z}+\frac{1}{2}$ if necessary. Then \[\hat{i}_{\min}^{n+1}-\hat{i}_{\min}^n=-\frac{p}{2}\aand \hat{i}_{\max}^{n+1}-\hat{i}_{\max}^{\mu}=\frac{np-q}{2}\]

An easy way to understand the grading shift was described in \cite[Remark 4.15]{LY2020}. Note that the grading shift of a map between two spaces equals half of the intersection number between $\partial S$ and the curve corresponding to the third space up to the sign, while the sign depends on the choice of the sign in the bypass map. For example, we have $\partial S\cdot \mu=p$, so the grading shifts of $\psi_{\pm,n+1}^n$ are $\mp p/2$.

\epf
\brem
The reason to use balanced sutured manifolds with reversed orientation is because of the above bypass exact triangles.
\erem
\brem
 If we do not mention gradings, the above triangles and the results in the rest of this subsection (except Corollary \ref{cor: psi^n_+,n+1 is an iso} and Lemma \ref{lem: F_n and G_n are iso when n large} since the statements involve gradings) also hold without the assumption that $K$ is rationally null homologous since the proofs only involve the neighborhood of $\partial(-Y\backslash N(K))$.
\erem

\bcor\label{cor: psi^n_+,n+1 is an iso}
For any sufficiently large integer $n$, we have the following properties for restrictions of maps.
\begin{enumerate}
	\item The map $\psp{n}{n+1}|_{\sutg{n}{i}}$ is an isomorphism when 
	\[i< \frac{1}{2}\bigg(np-q+\chi(S)\bigg).\]
	\item The map $\psm{n}{n+1}|_{\sutg{n}{i}}$ is an isomorphism when 
	\[i> -\frac{1}{2}\bigg(np-q+\chi(S)\bigg)\]
	\item For any grading $i$ such that 
	\[-\frac{1}{2}\bigg(np-q+\chi(S)\bigg) < i< \frac{1}{2}\bigg((n-2)p-q+\chi(S)\bigg),\]
	there is an isomorphism
	\[(\psp{n}{n+1})^{-1}\circ\psm{n}{n+1}:\sutg{n}{i}\xra{\cong}\sutg{n}{i+p}.\]
	\item The map $\psm{-n}{1-n}|_{\sutg{-n}{i}}$ is an isomorphism when 
	\[i<\frac{1}{2}\bigg((n-2)p+q+\chi(S)\bigg)\]
	\item The map $\psp{-n}{1-n}|_{\sutg{-n}{i}}$ is an isomorphism when 
	\[i>-\frac{1}{2}\bigg((n-2)p+q+\chi(S)\bigg)\]
	\item For any grading $i$ such that
	\[-\frac{1}{2}\bigg((n-2)p+q+\chi(S)\bigg)< i<\frac{1}{2}\bigg((n-4)p+q+\chi(S)\bigg)\]
	there is an isomorphism
	\[(\psp{-n}{1-n})^{-1}\circ\psm{-n}{1-n}:\sutg{-n}{i}\xra{\cong}\sutg{-n}{i+p}.\]

\end{enumerate}
\ecor
\bpf
It is a combination of Lemma \ref{lem: vanishing grading} and Lemma \ref{lem: bypass n,n+1,mu}.
\epf
\bdefn
The maps in Lemma \ref{lem: bypass n,n+1,mu} are called \textbf{bypass maps}. The ones with subscripts $+$ and $-$ are called \textbf{positive} and \textbf{negative bypass maps}, respectively. We will use $\pm$ to denote one of the bypass maps. For any integer $n$ and any positive integer $k$, define $$\Psi_{\pm,n+k}^n\deq \psi_{\pm,n+k}^{n+k-1}\circ\cdots\circ \psi_{\pm,n+1}^n:\mathbf{\Gamma}_{n}\rightarrow \mathbf{\Gamma}_{n+k}.$$
\edefn
In \cite[Section 4.4]{LY2020}, we proved many commutative diagrams for bypass maps, which we restate as follows by notations introduced before.
\begin{lem}\label{lem: com diag for n,n+1,n+2}
For any $n\in\intg$, we have the following commutative diagrams up to scalars.
\begin{equation*}
	\xymatrix{
	\sut{n}\ar[rr]^{\psp{n}{n+1}}\ar[dd]_{\psm{n}{n+1}}&&\sut{n+1}\ar[dd]^{\psm{n+1}{n+2}}\\
	&&\\
	\sut{n+1}\ar[rr]^{\psp{n+1}{n+2}}&&\sut{n+2}
	}
	\xymatrix{
	\sut{n+2}\ar[rr]^{\psp{n+2}{\mu}}\ar[dd]_{\psm{n+2}{\mu}}&&\sut{\mu}\ar[dd]^{\psp{\mu}{n}}\\
	&&\\
	\sut{\mu}\ar[rr]^{\psm{\mu}{n}}&&\sut{n}
	}
\end{equation*}
\end{lem}
\bpf
The first diagram follows from \cite[Lemma 4.33]{LY2020}. Note that the proof only used the functionality of the contact gluing map and did not depend on the assumption that $K$ is rationally null-homologous. The second diagram is obtained from the first diagram by changing the choice of the framed knot. Explicitly, let $K^\p$ be the dual knot corresponding to $\sut{n+1}$. Let $\mu^\p=-(n+1)\mu+\lambda$ denote its meridian. Then $\lambda^\p=-\mu$ is a framing of $K^\p$. Applying the first diagram to $K^\p$, we will obtain the second diagram for the original $K$. Note that the sign of the bypass map may switch when we regard it as the bypass map for the original knot. That is the reason for the signs in the second diagram. This can be double-checked by keeping track of the grading shifts.
\epf

\begin{lem}\label{lem: comm diag for n,n+1,mu}
For any $n\in\intg$, we have the following commutative diagrams up to scalars
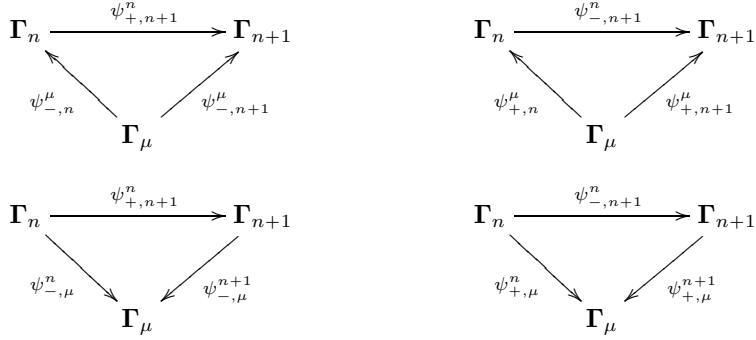
\begin{center}
\begin{minipage}{0.4\textwidth}
	\begin{equation*}
\xymatrix{
\sut{n}\ar[rr]^{\psp{n}{n+1}}&&\sut{n+1}\\
&\sut{\mu}\ar[lu]^{\psm{\mu}{n}}\ar[ru]_{\psm{\mu}{n+1}}&
}	
\end{equation*}
\end{minipage}
\begin{minipage}{0.4\textwidth}
	\begin{equation*}
\xymatrix{
\sut{n}\ar[rr]^{\psm{n}{n+1}}&&\sut{n+1}\\
&\sut{\mu}\ar[lu]^{\psp{\mu}{n}}\ar[ru]_{\psp{\mu}{n+1}}&
}	
\end{equation*}
\end{minipage}
\end{center}

\begin{center}
\begin{minipage}{0.4\textwidth}
	\begin{equation*}
\xymatrix{
\sut{n}\ar[rr]^{\psp{n}{n+1}}\ar[dr]_{\psm{n}{\mu}}&&\sut{n+1}\ar[dl]^{\psm{n+1}{\mu}}\\
&\sut{\mu}&
}	
\end{equation*}
\end{minipage}
\begin{minipage}{0.4\textwidth}
		\begin{equation*}
\xymatrix{
\sut{n}\ar[rr]^{\psm{n}{n+1}}\ar[dr]_{\psp{n}{\mu}}&&\sut{n+1}\ar[dl]^{\psp{n+1}{\mu}}\\
&\sut{\mu}&
}	
\end{equation*}\end{minipage}
\end{center}
\end{lem}

There are more bypass triangles involving more complicated sutures, which are obtained from changing the choice of the framed knot as in the proof of Lemma \ref{lem: com diag for n,n+1,n+2}.
\blem\label{lem: bypass n+1,n,2n+1/2}
For a knot $K\subset Y$ and $n\in\intg$, there are two graded bypass exact triangles
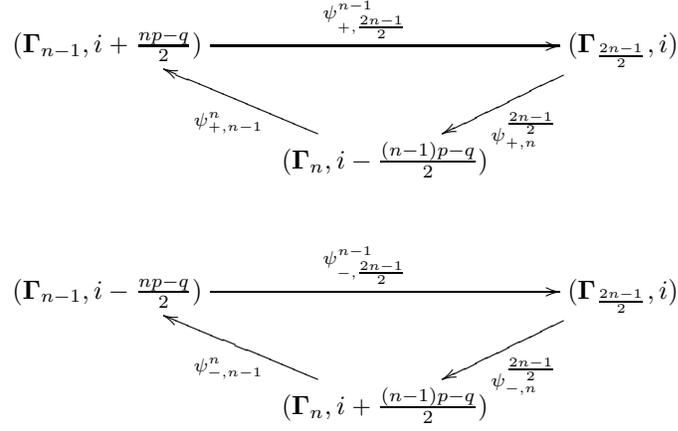
\begin{equation*}
\xymatrix{
\sutg{n-1}{i+\frac{np-q}{2}}\ar[rr]^{\psp{n-1}{\frac{2n-1}{2}}}&&\sutg{\frac{2n-1}{2}}{i}\ar[dl]^{\psp{\frac{2n-1}{2}}{n}}\\
&\sutg{n}{i-\frac{(n-1)p-q}{2}}\ar[ul]^{\psp{n}{n-1}}&
}	
\end{equation*}
\begin{equation*}
\xymatrix{
\sutg{n-1}{i-\frac{np-q}{2}}\ar[rr]^{\psm{n-1}{\frac{2n-1}{2}}}&&\sutg{\frac{2n-1}{2}}{i}\ar[dl]^{\psm{\frac{2n-1}{2}}{n}}\\
&\sutg{n}{i+\frac{(n-1)p-q}{2}}\ar[ul]^{\psm{n}{n-1}}&
}	
\end{equation*}
\elem



\blem\label{lem: comm diag for n+1 to n}
For a knot $K\subset Y$ and $n\in\intg$, there are commutative diagrams up to scalars
\begin{center}
\begin{minipage}{0.4\textwidth}
	\begin{equation*}
\xymatrix{
	\sut{\mu}\ar[rr]^{\psp{\mu}{n-1}}\ar[dd]_{\psm{\mu}{n-1}}&&\sut{n-1}\ar[dd]^{\psp{n-1}{\frac{2n-1}{2}}}\\
	&&\\
	\sut{n-1}\ar[rr]^{\psm{n-1}{\frac{2n-1}{2}}}&&\sut{\frac{2n-1}{2}}
	}
	\end{equation*}
\end{minipage}
\begin{minipage}{0.4\textwidth}
	\begin{equation*}
\xymatrix{
	\sut{\frac{2n-1}{2}}\ar[rr]^{\psp{\frac{2n-1}{2}}{n}}\ar[dd]_{\psm{\frac{2n-1}{2}}{n}}&&\sut{n}\ar[dd]^{\psm{n}{\mu}}\\
	&&\\
	\sut{n}\ar[rr]^{\psp{n}{\mu}}&&\sut{\mu}
	}
	\end{equation*}
\end{minipage}

\begin{minipage}{0.4\textwidth}
	\begin{equation*}
	\xymatrix{
\sut{\mu}\ar[rr]^{\psp{\mu}{n-1}}&&\sut{n-1}\\
&\sut{n}\ar[lu]^{\psm{n}{\mu}}\ar[ru]_{\psp{n}{n-1}}&
}
\end{equation*}
\end{minipage}
\begin{minipage}{0.4\textwidth}
	\begin{equation*}
\xymatrix{
\sut{\mu}\ar[rr]^{\psm{\mu}{n-1}}&&\sut{n-1}\\
&\sut{n}\ar[lu]^{\psp{n}{\mu}}\ar[ru]_{\psm{n}{n-1}}&
}
\end{equation*}
\end{minipage}
\end{center}

\begin{center}
\begin{minipage}{0.4\textwidth}
	\begin{equation*}
\xymatrix{
\sut{n-1}\ar[rr]^{\psp{\frac{2n-1}{2}}{n-1}}\ar[dr]_{\psp{n-1}{n}}&&\sut{\frac{2n-1}{2}}\ar[dl]^{\psm{\frac{2n-1}{2}}{n}}\\
&\sut{n}&
}	
\end{equation*}
\end{minipage}
\begin{minipage}{0.4\textwidth}
		\begin{equation*}
\xymatrix{
\sut{n-1}\ar[rr]^{\psm{\frac{2n-1}{2}}{n-1}}\ar[dr]_{\psm{n-1}{n}}&&\sut{\frac{2n-1}{2}}\ar[dl]^{\psp{\frac{2n-1}{2}}{n}}\\
&\sut{n}&
}		
\end{equation*}\end{minipage}
\end{center}
\elem
\brem\label{rem: changing the framed knot}
The choices of positive and negative bypass maps in Lemma \ref{lem: comm diag for n+1 to n} seem to be different from Lemma \ref{lem: com diag for n,n+1,n+2} and Lemma \ref{lem: comm diag for n,n+1,mu}. But indeed they are the same up to changing the framed knot. In particular, the grading shifts match. Note that similar to the second diagram in Lemma \ref{lem: com diag for n,n+1,n+2}, we always use the notations of the bypass maps for the original knot, while the signs may change if the maps are regarded as the bypass maps of the dual knot.
\erem
Suppose $\alpha$ is a connected non-separating simple closed curve on $\partial (Y\backslash N(K))$. We can push $\alpha$ into the interior of $Y\backslash N(K)$. For any fixed suture on $\partial(Y\backslash N(K))$ and a closure of the sutured manifold, the push-off of $\alpha$ is inside the closure, which is a closed $3$-manifold. We can then take the framing on $\alpha$ induced by the surface $\partial (Y\backslash N(K))$ and there is an exact triangle associated to the instanton Floer homology of the $(-1)$- $0$- and $\infty$-surgeries along the push-off of $\alpha$. Since the push-off of $\alpha$ is disjoint from $\partial(Y\backslash N(K))$, the exact triangle descends to one between corresponding sutured instanton Floer homologies.

According to \cite[Section 4]{baldwin2016contact}, when $\alpha$ intersects the suture at two points, the $0$-surgery along the push-off of $\alpha$ (with framing induced by $\partial (Y\backslash N(K))$) corresponds to a $2$-handle attachment along $\alpha$. Note that attaching a $2$-handle along $\alpha\subset \partial (Y\backslash N(K))$ will change the $3$-manifold from $Y\backslash N(K)$ to $Y_{\alpha}(K)\backslash B^3$, where $Y_{\alpha}(K)$ is the Dehn surgery along $K$ with slope specified by $\alpha$. We write
$$\mathbf{Y_{r/s}}\deq I^{\sharp}(-Y_{-r/s}(K)),$$
and in particular
$$\mathbf{Y}_n\deq I^{\sharp}(-Y_{-n}(K))\aand \dehny{}\deq I^{\sharp}(-Y).$$
\begin{lem}[{\cite[Lemma 3.21]{LY2020}}]\label{lem: surgery triangles}
	For any $n\in\intg$, we have the following exact triangles.
\begin{center}
\begin{minipage}{0.4\textwidth}
	\begin{equation*}
\xymatrix{
\sut{n}\ar[rr]^{H_n}&&\sut{n+1}\ar[dl]^{F_{n+1}}\\
&\dehny{}\ar[ul]^{G_n}&
}	
\end{equation*}
\end{minipage}
\begin{minipage}{0.4\textwidth}
		\begin{equation*}
\xymatrix{
\sut{\mu}\ar[rr]^{A_{n-1}}&&\sut{n-1}\ar[dl]^{B_{n-1}}\\
&\dehny{n}\ar[ul]^{C_n}&
}	
\end{equation*}\end{minipage}
\end{center}
\end{lem}
\bpf
To obtain the first exact triangle, we can take the sutured manifold $(-(Y\backslash N(K)),-\Gamma_n)$, and take a meridian $\alpha\subset \partial(Y\backslash N(K))$. As explained before the lemma, there is a surgery exact triangle associated to the sutured instanton Floer homology of the three sutured manifolds obtained by taking $(-1)$-, $0$-, and $\infty$-surgeries \cite[Theorem 2.1]{scaduto2015instanton}; see also \cite{floer1990knot} for the original construction and \cite[Proof of Theorem 1.21, especially (16)-(19)]{baldwin2018khovanov} for the resolution of the subtlety of the bundle data. 

The $\infty$-surgery will keep the manifold $(-(Y\backslash N(K)),-\Gamma_n)$. The $(-1)$-surgery changes the framing and hence we obtain $(-(Y\backslash N(K)),-\Gamma_{n+1})$. The $0$-surgery, as discussed above, gives rise to the manifold $Y_{\alpha}(K)\backslash B^3$ which is $Y\backslash B^3$ since $\alpha$ is the meridian. Hence we obtain the expected triangle. The second exact triangle in the statement of the lemma is obtained similarly by taking $\alpha$ to be a curve on $\partial (Y\backslash N(K))$ having slope $-n$ instead of a meridian.
\epf
\brem\label{rem: unknot}
From \cite[Section 4]{baldwin2016contact}, we know the $0$-surgery corresponds to a $2$-handle attachment and a $1$-handle attachment. Hence $\dehny{}$ in the above lemma is indeed $\khii(-Y,U)$, where $U$ is the unknot in $-Y$ bounding an embedded disk. By \cite[Section 4]{baldwin2016contact}, a $1$-handle attachment does not change the closure of the balanced sutured manifold, and then there is a canonical identification between $\khii(-Y,U)$ and $I^\sharp(-Y)$. Hence we can abuse the notations. The same discussion also applies to $\mathbf{Y}_{n}$.
\erem
Furthermore, we proved the following properties in \cite{LY2020}. Note that the assumption that $K$ is the dual knot of a null-homologous knot in that paper is inessential by remarks in the beginning of \cite[Section 5]{LY2020}. The inequalities of the gradings are from Corollary \ref{cor: psi^n_+,n+1 is an iso}.
\blem[{\cite[Lemma 3.21 and Lemma 4.9]{LY2020}}]\label{lem: comm diag for n,n+1,dehn}
For any $n\in\intg$, we have the following commutative diagrams up to scalars
\begin{center}
\begin{minipage}{0.4\textwidth}
	\begin{equation*}
\xymatrix{
\sut{n}\ar[rr]^{\psp{n}{n+1}}&&\sut{n+1}\\
&\dehny{}\ar[lu]^{G_n}\ar[ru]_{G_{n+1}}&
}	
\end{equation*}
\end{minipage}
\begin{minipage}{0.4\textwidth}
		\begin{equation*}
\xymatrix{
\sut{n}\ar[rr]^{\psm{n}{n+1}}&&\sut{n+1}\\
&\dehny{}\ar[lu]^{G_n}\ar[ru]_{G_{n+1}}&
}	
\end{equation*}
\end{minipage}
\end{center}

\begin{center}
\begin{minipage}{0.4\textwidth}
	\begin{equation*}
\xymatrix{
\sut{n}\ar[rr]^{\psp{n}{n+1}}\ar[dr]_{F_{n}}&&\sut{n+1}\ar[dl]^{F_{n+1}}\\
&\dehny{}&
}	
\end{equation*}
\end{minipage}
\begin{minipage}{0.4\textwidth}
		\begin{equation*}
\xymatrix{
\sut{n}\ar[rr]^{\psm{n}{n+1}}\ar[dr]_{F_{n}}&&\sut{n+1}\ar[dl]^{F_{n+1}}\\
&\dehny{}&
}	
\end{equation*}\end{minipage}
\end{center}
\elem

\blem[{\cite[Lemma 4.17, Proposition 4.26, Lemma 4.29 and Proposition 4.32]{LY2020}}]\label{lem: F_n and G_n are iso when n large}
Let $F_n$ and $G_n$ be defined as in Lemma \ref{lem: surgery triangles}. Then for any sufficiently large integer $n$, we have the following properties.
\begin{enumerate}
	\item The map $G_{n-1}$ is zero and $F_n$ is surjective. Moreover, for any grading $$-(np-q+\chi(S))/2< i_0<(np-q+\chi(S))/2-p+1,$$the restriction of the map
	$$F_n:\bigoplus_{i=0}^{p-1}\sutg{n}{i_0+i}\to \dehny{}$$
	is an isomorphism.
	\item The map $F_{-n+1}$ is zero and $G_{-n}$ is injective. Moreover, for any grading $$-((n-2)p+q+\chi(S))/2< i_0< ((n-2)p+q+\chi(S))/2-p+1,$$the map
	$${\rm Proj}\circ G_{-n}:\dehny{}\to \bigoplus_{i=0}^{p-1}\sutg{-n}{i_0+i},$$
	is an isomorphism, where 	$${\rm Proj}:\sut{-n}\ra\bigoplus_{i=0}^{p-1}\sutg{-n}{i_0+i}$$
	is the projection.
\end{enumerate}
\elem
The following lemma is a special case of Proposition \ref{prop: -1 surgery and bypass}, which we will prove later.
\blem\label{lem: -1 surgery and bypass}
For any $n\in\intg$, let the maps $H_{n}$ and $\psi^{n}_{\pm,n+1}$ be defined as in Lemma \ref{lem: surgery triangles} and Lemma \ref{lem: bypass n,n+1,mu} respectively. Then there exist $c_1,c_2\in\mathbb{C}\backslash\{0\}$ such that
$$H_{n}=c_1\psp{n}{n+1}+c_2\psm{n}{n+1}$$
\elem

\subsection{Fixing the scalars}\label{subsec Fixing the scalars}
By construction, sutured instanton homology forms a projectively transitive system, which means all the spaces and maps between spaces are well-defined only up to non-zero scalars. When the balanced sutured manifold is obtained from a framed knot as in the last subsection, we can make some canonical choices to reduce the projective ambiguities.

Suppose $K\subset Y$ is a framed knot with the meridian $\mu$ and the framing $\lambda$. Fix a knot complement $Y\backslash N(K)$ and the suture $\Ga_\mu$. We fix a special choice of a (marked odd) closure of $(Y\backslash N(K), \Ga_\mu)$ following the construction in \cite[Formula (18)]{kronheimer2010knots}. 

Let $F$ be a closed, oriented, connected surface of genus at least $2$. Suppose $c_0$ is a non-separating curve in $F$. Let $c=\operatorname{pt}\times c_0\subset S^1\times F$ and let $(\mu_c,\lambda_c)$ be the meridian and the longitude of $c$ (the latter comes from the surface framing). Let\begin{equation}\label{eq: special closure}
    (\widebar{Y}_\mu,R)=(S^1\times F\backslash N(c)\cup_{(\mu_c,\lambda_c)\sim (\lambda,\mu)}Y\backslash N(K),\operatorname{pt}^\p\times F),
\end{equation}where $\operatorname{pt}^\p$ is a point different from $\operatorname{pt}$. We can pick $\alpha=S^1\times \operatorname{pt}^\pp$ and $\eta\subset R$ be a curve intersecting $\operatorname{pt}^\p \times c_0$ once. Since $\alpha\cdot R=1$ and $\eta\cdot R=0$, the pair $(\widebar{Y}_\mu, \alpha\cup \eta)$ defines an instanton Floer homology in the setting of \cite[Section 7.1]{kronheimer2010knots}. Moreover, $\mathcal{D}_\mu=(Y,R,\eta,\al)$ forms a marked odd closure in \cite[Definition 9.2]{baldwin2015naturality}, which was used in the naturality result \cite[Theorem 9.17]{baldwin2015naturality}. The reason for $g(F)\ge 2$ is to apply the naturality result (\textit{c.f.} \cite[Remark 9.4]{baldwin2015naturality}).

Similarly, for $(Y\backslash N(K), \Ga_n)$ and $(Y\backslash N(K), \Ga_{\frac{2n-1}{2}})$, we fix closures $\mathcal{D}_{n}$ and $\mathcal{D}_{\frac{2n-1}{2}}$ as in (\ref{eq: special closure}), except replacing the gluing map $(\mu_c,\lambda_c)\sim (\lambda,\mu)$ by $(\mu_c,\lambda_c)\sim (-\mu,-n\mu+\lambda)$ and $(\mu_c,\lambda_c)\sim (-n\mu+\lambda,(1-2n)\mu+2\lambda)$, respectively.

For the sutured manifold $(Y\backslash B^3,\delta)$, we regard it as $(Y\backslash N(U),\Ga_{\mu,U})$ by Remark \ref{rem: unknot}, where $U$ is the unknot and $\Ga_{\mu,U}$ is meridian suture on the unknot complement. Then we apply the above construction to obtain a special closure of $(Y\backslash B^3,\delta)$. We reverse the orientations of the chosen closures when the orientations of the sutured manifolds are reversed. Note that we do not choose canonical closures for $(Y_n(K)\backslash B^3,\delta)$ since we only care about the dimension of its framed instanton homology.

After fixing the choices of closures, we can view $\sut{n}$ and $\dehny{}$ as actual vector spaces, and then the elements in them are well-defined. Strictly speaking, we also need to choose extra data such as the metric and the perturbation on the closure to define the instanton Floer homology of the closure, but different choices of metrics and perturbations now induce a transitive system of vector spaces, from which we can construct an actual vector space. So, we omit the discussion on those extra data. 

The construction of bypass maps and surgery maps may not be realized as cobordism maps between the chosen closures, but the construction of the projectively transitive system (\textit{c.f.} \cite[Definition 9.18]{baldwin2015naturality}) guarantees the existence of such maps up to scalars. Now we make (non-canonical) choices of the maps to get rid of the scalar ambiguities in the commutative diagrams mentioned in the last subsection. 

We first assume that $I^{\sharp}(Y)\neq 0$. When $I^{\sharp}(Y) = 0$, the first exact triangle in Lemma \ref{lem: surgery triangles} is trivial for any $n$ and hence the maps $F_n$ and $G_n$ that play an important role in later sections are both zero. This makes fixing the scalars a somewhat straightforward job: we just need to fix the scalars for the bypass maps. Also, it is worth mentioning that, the Euler characteristic result in \cite[Corollary 1.4]{scaduto2015instanton} implies that $I^{\sharp}(Y)\neq 0$ for any rational homology sphere $Y$ and, up to author's knowledge, there is no known closed oriented $3$-manifold $Y$ with $I^{\sharp}(Y) = 0$.

To help us fixing the scalars, suppose the maps $F_n$ and $G_n$ are defined as in the proof of Lemma \ref{lem: surgery triangles}, and we define
\begin{equation}\label{eq: scalars, n_G and n_F}
	n_{G} = \min \{n\in\intg~|~ G_n = 0\} \text{ and } n_F = \max \{ n\in\intg~|~F_n = 0 \}.
\end{equation}
We have the following basic properties for these indices.
\blem\label{lem: n_F geq n_G}
Assuming $I^{\sharp}(Y)\neq 0$. Suppose $n_G$ and $n_F$ are defined as in Equation (\ref{eq: scalars, n_G and n_F}). Then we have
\[
-\infty<n_F\leq n_G<\infty.
\]

Moreover, we have $G_n = 0$ if and only if $n\geq n_G$ and $F_n = 0$ if and only if $n\leq n_F$.
\elem
\bpf
The fact that $-\infty < n_F < \infty$ and $-\infty < n_G < \infty$ follow from Lemma \ref{lem: F_n and G_n are iso when n large} and the fact that they fits into an exact triangle as in Lemma \ref{lem: surgery triangles}. Next, the commutative diagrams in Lemma \ref{lem: comm diag for n,n+1,dehn} implies that $G_{n+1}=0$ whenever $G_n=0$ and thus we know $G_n = 0$ if and only if $n\geq n_G$. The argument for $F_n$ is similar. Finally, by definition we know $G_{n_G}=0$ and hence from the exact triangle we know that 
\[\im F_{n+1}=\ker G_{n_G}=I^{\sharp}(Y)\neq 0.\]
Hence we conclude that $n_F\leq n_G$.
\epf

By Lemma \ref{lem: F_n and G_n are iso when n large}, we can pick a sufficiently large integer $n_0$ such that $-n_0< n_F\leq n_G$. Pick arbitrary representatives of the maps 
$$G_{-n_0},\psp{-n_0}{\mu},\psm{-n_0}{\mu},\psp{\mu}{-n_0},\psm{\mu}{-n_0}$$
and we also pick arbitrary representatives of the maps
$$\psp{n}{n+1},\psp{n-1}{\frac{2n-1}{2}}$$ for all $n\in\intg$. 

Now we explain how to fix the scalars for maps with $n\ge n_0$. Note that we have already chosen a representative for $\psi^{-n_0}_{+,-n_0+1}$ and $G_{-n_0}$. From Lemma \ref{lem: comm diag for n,n+1,dehn}, we have
\[
G_{-n_0+1}\doteq\psp{-n_0}{-n_0+1}\circ G_{-n_0}
\]
where $\doteq$ means commutative up to a non-zero scalar. We can choose a representative of $G_{-n_0+1}$ to obtain an equality
\[
G_{-n_0+1}=\psp{-n_0}{-n_0+1}\circ G_{-n_0}
\]
We then choose a representative of $\psm{-n_0}{-n_0+1}$ so that
\[
G_{-n_0+1}=\psm{-n_0}{-n_0+1}\circ G_{-n_0}.
\]

Next, we pick representatives of the maps $\psm{n}{n+1}$ inductively, with the base case $n=-n_0$ constructed above, so that 
\[
\psm{n}{n+1}\circ\psp{n-1}{n}=\psp{n}{n+1}\circ\psm{n-1}{n}
\]
hold for all $n\geq -n_0+1$. If the compositions happen to be zero, we could pick an arbitrary representative since the diagram will be trivially satisfied. We will discuss the ambiguity arising from the possibility that $\psp{n}{n+1}\circ\psm{n-1}{n}=0$ more carefully later. 

Similarly, we pick the maps $\psp{n}{\mu}$ inductively to satisfy the commutative diagram
\[
\psp{n}{\mu}\circ\psm{n-1}{n}=\psp{n-1}{\mu}.
\]
We can choose representatives of $\psm{n}{\mu}$, $\psp{\mu}{n}$, and $\psm{\mu}{n}$ in a similar manner.

Furthermore, the representatives of 
$$\psm{\frac{2n-1}{2}}{n},\psp{n}{n-1},\psm{n}{n-1},\psm{n-1}{\frac{2n-1}{2}},\psp{\frac{2n-1}{2}}{n}$$
can be chosen according to Lemma \ref{lem: comm diag for n+1 to n}. As mentioned in Remark \ref{rem: changing the framed knot}, we always use the notations of the bypass maps for the original knot even though we consider some dual knots in the proofs. Hence here we first fix the knot $K\subset Y$ and then fix the representatives, while we do not fix the representatives by any commutative diagrams for the dual knot in the proofs.

Next, we deal with maps $G_n$ and $F_n$ in Lemma \ref{lem: surgery triangles}. We choose representatives of the maps $G_n$ inductively so that
\[
G_{n+1}=\psp{n}{n+1}\circ G_{n}
\]
is satisfied for all $n\geq -n_0$. We pick arbitrary representatives of the map $F_{n_F+1}$ and then pick $F_n$ inductively so that
\[
F_{n+1}\circ \psp{n}{n+1}= F_{n}
\]
is satisfied for all $n\geq n_{F}+1$. We can then use induction to prove the following two equalities.
\begin{itemize}
	\item $G_{n+1}=\psm{n}{n+1}\circ G_{n}$ for all $n\geq -n_0$, and
	\item $F_{n+1}\circ \psm{n}{n+1}=c\cdot F_{n}$ for a non-zero scalar $c$ that is independent of $n$ with $n\geq n_F+1$.
\end{itemize}

We verify the equality for $G$ first. The base case $n=n_0$ is by construction. Assuming we have already established the equality for $n$, from Lemma \ref{lem: com diag for n,n+1,n+2} and Lemma \ref{lem: comm diag for n,n+1,dehn}, we have
\[
\begin{aligned}
	G_{n+2}&=\psp{n+1}{n+2}\circ G_{n+1}\\
	(\text{Inductive hypothesis})&=\psp{n+1}{n+2}\circ\psm{n}{n+1}\circ G_{n}\\
	(\text{Lemma }\ref{lem: com diag for n,n+1,n+2})&=\psm{n+1}{n+2}\circ\psp{n}{n+1}\circ G_{n}\\
	&=\psm{n+1}{n+2}\circ G_{n+1}.
\end{aligned}
\]
The argument for $F_n$ is similar, once we take $c\neq 0$ to be the complex number such that
\begin{equation}\label{eq: scalars, c}
	F_{n_{F}+2}\circ \psm{n_F+1}{n_F+2}=c\cdot F_{n_F+1}.
\end{equation}

The remaining issues are summarized as follows:
\begin{itemize}
	\item [(i)] When choosing representatives of $\psm{n}{n+1}$, we use the commutative diagram
	\[\psm{n}{n+1}\circ\psp{n-1}{n}\doteq\psp{n}{n+1}\circ\psm{n-1}{n}\]
	However, when $\psp{n}{n+1}\circ\psm{n-1}{n}=0$, there is no unique choice of $\psm{n}{n+1}$ and one might worry that different choices of $\psm{n}{n+1}$ may affect the commutative diagrams in Lemma \ref{lem: comm diag for n,n+1,dehn}.
	\item [(ii)] By Proposition \ref{prop: -1 surgery and bypass}, we can assume that the map $H_n$ in the exact triangle in Lemma \ref{lem: surgery triangles} to have the form 
	\[H_n=\psp{n}{n+1}-c_n\cdot \psm{n}{n+1}.\]
	We want to pin down the values of $c_n$.
	\item [(iii)] We want to get rid of the scalar $c$ in Equation \ref{eq: scalars, c}.
\end{itemize}
We treat these issues in several different cases.

{\bf Case 1}. $n_F+3\leq n_G$. In this case, we know that 
\[\psp{n}{n+1}\circ\psm{n-1}{n}\neq 0\]
for any integer $n$. Indeed, if $n+1<n_G$, we know from Lemma \ref{lem: comm diag for n,n+1,dehn} that
\[\psp{n}{n+1}\circ\psm{n-1}{n}\circ G_{n-1}=G_{n+1}\neq 0.\]
If $n+1\geq n_G\geq n_F+3$, instead of the above equation involving $G$, we have
\[F_{n+1}\circ \psp{n}{n+1}\circ\psm{n-1}{n}=c\cdot F_{n-1}\neq 0.\]
This implies that inductively we can fix a unique representative of $\psm{n}{n+1}$ for all $n\geq -n_0+1$.

For any integer $n$ with $-n_0\le n<n_G-1$, we have $G_{n+1}\neq 0$. Then we take an element $\alpha\in\dehny{}$ so that $G_{n+1}(\alpha)\neq 0$. Then we can solve the scalar $c_n$ as follows.
\[
\begin{aligned}
	0&=H_n\circ G_n(\alpha)\\
&=(\psp{n}{n+1}-c_{n}\cdot \psm{n}{n+1})\circ G_n(\alpha)\\
&=\psp{n}{n+1}\circ G_n(\alpha)-c_{n}\cdot \psm{n}{n+1}\circ G_n(\alpha)\\
&=(1-c_n)\cdot G_{n+1}(\alpha)\\
\end{aligned}
\]
Hence we conclude that $c_n=1$. In particular, we can take $n=n_F+1<n_G-1$. Note $F_{n}\neq 0$ so we can take $x\in\sut{n}$ so that $F_{n}(x)\neq 0$. Then we have
\[
\begin{aligned}
	0&=F_{n+1}\circ H_{n}(x)\\
	&=F_{n+1}\circ(\psp{n}{n+1}-\psm{n}{n+1})(x)\\
	&=F_{n+1}\circ\psp{n}{n+1}(x)-F_{n+1}\circ \psm{n}{n+1}(x)\\
	&=(1-c)\cdot F_{n}(x)\\
\end{aligned}
\]
This implies that $c=1$ as well. Now for any $n\geq n_G-1\geq n_F+2$, we can take $x\in\sut{n}$ such that $F_{n}(x)\neq 0$. we have
\[
\begin{aligned}
	0&=F_{n+1}\circ H_n(x)\\
	&=F_{n+1}\circ(\psp{n}{n+1}-c_n\cdot \psm{n}{n+1})(x)\\
	&=(1-c_n)\cdot F_n(x)\\
\end{aligned}
\]
Hence we conclude that $c_n=1$. In summary, in Case 1, we have the following.
\begin{itemize}
	\item We can fix a unique representative of $\psm{n}{n+1}$ for any $n\geq -n_0$.
	\item We have $c_n=1$ for any $n\geq -n_0$.
	\item We have $c=1$.
\end{itemize}

{\bf Case 2}. $n_G=n_F+2$. Note that some arguments in Case 1 still apply. We summarize as follows.
\begin{itemize}
	\item For $n<n_G-1$, we have $G_{n+1}\neq 0$ so there is a unique choice of $\psm{n}{n+1}$.
	\item For $n\geq n_G= n_F+2$, we have $F_{n-1}\neq 0$ so again there is a unique choice of $\psm{n}{n+1}$. 
	\item For $n< n_G-1$, we have $G_{n+1}\neq 0$ so $c_n=1$ in the expression of $H_n$.
	\item For $n\geq n_G-1=n_F+1$, we have $F_{n}\neq 0$ so $c_n=c^{-1}$.
\end{itemize}
To resolve the issue (i), the only nonfixed index is $n=n_F+1=n_G-1$. In case
\[\psp{n}{n+1}\circ\psm{n-1}{n}\neq 0,\]
there is a unique choice of $\psm{n}{n+1}$ so that
\[\psm{n}{n+1}\circ\psp{n-1}{n}=\psp{n}{n+1}\circ\psm{n-1}{n}.\]
Otherwise, we just fix any representative of $\psm{n}{n+1}$.

To resolve issues (ii) and (iii), we rescale the maps $F_n$ according to the grading on $\sut{n}$. To do this, for an integer $n$ and a grading $i$, define
\[j(n,i)=\lfloor\frac{i}{p}-\frac{n}{2}\rfloor.\]
It is straightforward that
\[j(n+1,i+\frac{p}{2})=j(n,i)\text{ and }j(n+1,i-\frac{p}{2})=j(n,i)-1.\]
Hence we define
\begin{equation}\label{eq: scalars, F-tilde}
	\widetilde{F}_n=\sum_{i}c^{j(n,i)}\cdot F_{n}\circ \text{Proj}^i_n,
\end{equation}
where the map
\[
\text{Proj}^i_n: \sut{n}\to\sutg{n}{i}
\]
is the projection. Equivalently, if we have an element $x\in \sutg{n}{i}$, we take
\[
\widetilde{F}_n(x)=c^{j(n,i)}\cdot F_n(x).
\]
We then need to verify the following two equalities for $\widetilde{F}$:
\begin{itemize}
	\item $\widetilde{F}_{n+1}\circ \psp{n}{n+1}=\widetilde{F}_{n}$, and
	\item $\widetilde{F}_{n+1}\circ \psm{n}{n+1}=\widetilde{F}_{n}$.
\end{itemize}
For the first one, note that $\psp{n}{n+1}$ increases the grading by $p/2$ from Lemma \ref{lem: bypass n,n+1,mu}. So assume $x\in \sutg{n}{i}$, we have
\[
\begin{aligned}
	\widetilde{F}_{n+1}\circ \psp{n}{n+1}(x)&=c^{j(n+1,i+\frac{p}{2})}\cdot F_{n+1}(x)\circ \psp{n}{n+1}\\
	&=c^{j(n,i)}\cdot F_{n}(x)\\
	&=\widetilde{F}_{n}(x)
\end{aligned}
\]
On the other hand, the map $\psm{n}{n+1}$ decreases the grading by $p/2$ from Lemma \ref{lem: bypass n,n+1,mu}, and so for $x\in\sutg{n}{i}$,
\[
\begin{aligned}
	\widetilde{F}_{n+1}\circ \psm{n}{n+1}(x)&=c^{j(n+1,i-\frac{p}{2})}\cdot F_{n+1}(x)\circ \psm{n}{n+1}\\
	&=c^{j(n,i)-1}\cdot c\cdot F_{n}(x)\\
	&=\widetilde{F}_{n}(x)
\end{aligned}
\]
Hence we can use $\widetilde{F}_n$ instead of $F_n$ to get rid of the scalar $c$. However, changing $F_n$ to $\widetilde{F}_n$, we might break the exact triangle in Lemma \ref{lem: surgery triangles}. Hence we need to alter the definition of $H_n$ as well. We define 
\begin{equation}\label{eq: scalar, H-tilde}
	\widetilde{H}_n=\psp{n}{n+1}-\psm{n}{n+1}
\end{equation}
and it remains to establish the following exact triangle:
\begin{equation*}
	\xymatrix{
	\sut{n}\ar[rr]^{\widetilde{H}_n}&&\sut{n+1}\ar[dl]^{\widetilde{F}_{n+1}}\\
	&\dehny{}\ar[ul]^{G_n}&
	}
\end{equation*}
When $n<n_G-2=n_F$, the triangle automatically holds. Indeed, for such $n$, we have $c_{n}=1$ so $\widetilde{H}_{n}=H_{n}$; the fact $F_{n+1}=0$ implies that $\widetilde{F}_{n+1}=0$, so the new exact triangle is exactly the original one in Lemma \ref{lem: surgery triangles}. 

When $n=n_F$, we still have $c_n=1$ and hence $\widetilde{H}_n=H_n$. So the exactness at $\sut{n}$ is from Lemma \ref{lem: surgery triangles}. By construction we have $\im \widetilde{F}_{n+1}=\im F_{n+1}$ for any $n$ so we conclude the exactness at $\dehny{}$. This also implies that
\[\dim \ker \widetilde{F}_{n+1}=\dim \ker F_{n+1}=\dim \im {H}_{n}.\]
Hence to show the exactness at $\sut{n+1}$, it remains to show that
\[
\im H_{n}\subset \ker\widetilde{F}_{n+1}.
\]
For any $x\in \sut{n}$, we have
\[
\begin{aligned}
	\widetilde{F}_{n+1}\circ H_n(x)&=\\
&=\widetilde{F}_{n+1}\circ(\psp{n}{n+1}-\psm{n}{n+1})(x)\\
&=\widetilde{F}_{n+1}\circ\psp{n}{n+1}(x)-\widetilde{F}_{n+1}\circ\psm{n}{n+1}(x)\\
&=\widetilde{F}_n(x)-\widetilde{F}_n(x)\\
&=0.
\end{aligned}
\]
Hence we are done.

Finally, we verify the exact triangle for $n\geq n_F+1$. Define a homomorphism
\[
\iota_n: \sut{n}\to\sut{n}
\]
as
\[\iota_n(x)=\sum_{i}c^{-j(n,i)}\text{Proj}^i_n(x).\]
It is clear that $\iota$ is an isomorphism as its inverse is
\[
\iota_n^{-1}(x)=\sum_{i}c^{j(n,i)}\text{Proj}^i_n(x),
\]
and from the construction of $\widetilde{F}$ in Equation (\ref{eq: scalars, F-tilde}), we have
\[
\iota_{n+1}(\ker \widetilde{F}_{n+1})= \ker F_{n+1}.
\]
From the fact that $H_n= \psp{n}{n+1}-{c}^{-1}\cdot \psm{n}{n+1}$ (we have $c_n=c^{-1}$), $\widetilde{H}_{n}=\psp{n}{n+1}-\psm{n}{n+1}$, and that $\psi^n_{\pm, n+1}$ shift the grading by $\mp p/2$, we conclude that
\[
\iota_{n+1}(\im \widetilde{H}_{n}) = \im H_n.
\]
As a result, we conclude the exactness at $\Gamma_{n+1}$. The exactness at $\dehny{}$ holds as above, since we still have
\[\im \widetilde{F}_{n+1}=\im F_{n+1}=\ker G_n.\]
Dimension counting similar to the above argument then implies that
\[\dim \ker \widetilde{H}_n=\dim \im G_n.\]
We then verify that
\[
\im G_n\subset\ker\widetilde{H}_n.
\]
For any $\alpha\in\dehny{}$, we have
\[
\begin{aligned}
	\widetilde{H}_{n}\circ G_{n}(\alpha)&=(\psp{n}{n+1}-\psm{n}{n+1})\circ G_n(\alpha)\\
	&=\psp{n}{n+1}\circ G_n(\alpha)-\psm{n}{n+1}\circ G_n(\alpha)\\
	&=G_{n+1}(\alpha)-G_{n+1}(\alpha)\\
	&=0.
\end{aligned}
\]

In summary, in Case 2, we did the following (extra things):
\begin{itemize}
	\item We choose representatives of $\psm{n}{n+1}$ for all $n\geq -n_0+1$ so that
	\[
	\psm{n}{n+1}\circ\psp{n-1}{n}=\psp{n}{n+1}\circ\psm{n-1}{n}.
	\]
	\item We define new maps $\widetilde{F}_n$ for all $n$ so that
	\begin{equation}\label{eq: scalar, equality for F-tilde}
		\widetilde{F}_{n}=\widetilde{F}_{n+1}\circ\psi^n_{\pm,n+1}.
	\end{equation}
	\item We define new maps $\widetilde{H}_n$ for all $n$ so that we have the following exact triangle for all $n$.
	\begin{equation}\label{eq: scalar, triangle}
	\xymatrix{
	\sut{n}\ar[rr]^{\widetilde{H}_n}&&\sut{n+1}\ar[dl]^{\widetilde{F}_{n+1}}\\
	&\dehny{}\ar[ul]^{G_n}&
	}
\end{equation}
\end{itemize}

{\bf Case 3}. $n_{G}=n_{F}+1$ or $n_G=n_F$. The situation and argument are similar to those in Case 2. We summarize the differences here:
\begin{itemize}
	\item In Case 3, the composition $\psp{n}{n+1}\circ\psm{n-1}{n}$ could only be zero when $n_G-1\le n \le n_F+1$. In case the composition is indeed $0$, we choose an arbitrary representative of the map $\psm{n}{n+1}$.
	\item We still define the maps $\widetilde{F}_n$ as in Equation (\ref{eq: scalars, F-tilde}) and the maps $\widetilde{H}_n$ as in Equation (\ref{eq: scalar, H-tilde}), and can verify the Equation (\ref{eq: scalar, equality for F-tilde}) and the exact triangle (\ref{eq: scalar, triangle}) as in Case 2.
\end{itemize}

Note from Lemma \ref{lem: n_F geq n_G}, we must have $n_G\geq n_F$, so the above three cases cover all situations.

Finally, we could extend the choice of representatives for all relevant maps for the indices $n<-n_0$. Note that when $n<-n_0$, we have that $G_n$ is injective and $F_n=0$. Hence we do not need to worry about the issues (i), (ii), and (iii).

\begin{conv}
 From now on, we write the maps $\widetilde{H}_n$ and $\widetilde{F}_n$ simply as $H_n$ and $F_n$, respectively. From the above discussion, when $K\subset Y$ is a fixed rationally null-homologous knot, we can assume the first commutative diagram in Lemma \ref{lem: com diag for n,n+1,n+2} and all commutative diagrams in Lemma \ref{lem: comm diag for n,n+1,mu}, Lemma \ref{lem: comm diag for n+1 to n} and Lemma \ref{lem: comm diag for n,n+1,dehn} hold without introducing scalars.
\end{conv}

\subsection{Algebraic lemmas}\label{subsec: Algebraic lemmas}

In this subsection, we introduce some lemmas in homological algebra. All graded vector spaces in this subsection are finite dimensional and over $\mathbb{C}$ and all maps are complex linear maps. For convenience, we will switch freely between long exact sequences and exact triangles.

From Section \ref{sec: preliminaries}, we know the sutured instanton homology is usually $\mathbb{Z}\oplus\mathbb{Z}_2$-graded, where we regard the $\mathbb{Z}_2$-grading as a homological grading. Many results in this subsection come from properties of the derived category of vector spaces over $\mathbb{C}$, for which people usually consider cochain complexes. However, for a $\mathbb{Z}_2$-graded space there is no difference between the chain complex and the cochain complex. Hence by saying a \textbf{complex} we mean a $\mathbb{Z}_2$-graded (co)chain complex, though all results apply to $\mathbb{Z}$-graded cochain complexes verbatim.

For a complex $C$ and an integer $n$, we write $C^n$ for its grading $n$ part (under the natural map $\mathbb{Z}\to \mathbb{Z}_2$). With this notation, we suppose the differential $d_C$ on $C$ sends $C^{n}$ to $C^{n+1}$. For any integer $k$, we write $C\{k\}$ for the complex obtained from $C$ by the grading shift $C\{k\}^n=C^{n+k}$. We write $H(C,d_C)$ or $H(C)$ for the homology of a complex $C$ with differential $d_C$. A $\mathbb{Z}_2$-graded vector space is regarded as a complex with the trivial differential. 

For a chain map $f:C\to D$, we write $\cone(f)$ for the \textbf{mapping cone} of $f$, \textit{i.e.}, the complex consisting of the space $D\oplus C\{1\}$ and the differential $$d_{\cone(f)}\deq \begin{bmatrix}
    d_D & -f \\
   0& -d_C
\end{bmatrix}.$$Then there is a long exact sequence$$\cdots\to H(C)\xra{f}H(D)\xra{i} H(\cone(f))\xra{p} H(C)\{1\}\to \cdots$$where $i$ sends $x\in D$ to $(x,0)$ and $p$ sends $(x,y)\in D\oplus C\{1\}$ to $-y$. If the differentials of $C$ and $D$ are trivial, then we know \begin{equation}\label{eq: ker cok}H(\cone(f))\cong \ker(f)\oplus \cok(f).\end{equation}
\brem\label{rem: different defn of mapping cone}
Our definitions about mapping cones follow from \cite{homologicalalgebra94}, which are different from those in \cite{Ozsvath2008integral,Ozsvath2011rational}.
\erem

Note that the derived category is a triangulated category, so it satisfies the octahedral lemma (for example, see \cite[Proposition 10.2.4]{homologicalalgebra94}). 
\blem[octahedral lemma]\label{lem: octahedral lemma}
Suppose $X,Y,Z,X^\p,Y^\p,Z^\p$ are $\mathbb{Z}_2$-graded vector spaces satisfying the following long exact sequences$$\begin{aligned}
        \cdots&\to X\xra{f} Y\xra{h} Z^\p\to X\{1\}\to\cdots\\
        \cdots&\to X\xra{g\circ f} Z\xra{h^\p} Y^\p\xra{l^\p} X\{1\}\to\cdots\\
        \cdots&\to Y\xra{g} Z\to X^\p\xra{l} Y\{1\}\to\cdots\\
    \end{aligned}$$Then we have the fourth long exact sequence$$\cdots\to Z^\p\xra{\psi} Y^\p\xra{\phi} X^\p \xra{h\{1\}\circ l} Z^\p\{1\}\to\cdots$$such that the following diagram commutes
    \begin{equation}\label{eq: octahedral 0}
\xymatrix{
	&Y^\p\ar@{..>}[dr]^{\phi}\ar[ddl]&\\
	Z^\p \ar@{..>}[ur]^{\psi}\ar[d]&&X^\p\ar@{..>}[ll]^{h\{1\}\circ l}\ar[ddl]\\
X\ar[dr]_{f}\ar[rr]^{g\circ f}&&Z\ar[u]\ar[uul]\\
	&Y\ar[uul]\ar[ur]_{g}&
	}
    \end{equation}
    where we omit the grading shifts and the notation for the maps $h,l,j$. We can also write (\ref{eq: octahedral 0}) in another form so that there is enough room to write the maps
    \begin{equation}\label{eq: octahedral}
        \xymatrix@R=6ex{
    &&&Z^\p\ar@{..>}[dr]^{\psi}\ar[rr]&&X\{1\}\ar[ddr]^{f\{1\}}&&
    \\
    &&Y\ar[dr]^{g}\ar[ur]^h&&Y^\p\ar@{..>}[ddrr]^{\phi}\ar[ur]^{l^\p}&&&
    \\
 &&&Z\ar[drrr]\ar[ur]^{h^\p}&&&Y\{1\}\ar[dr]^{h\{1\}}&
    \\
    X\ar[uurr]^{f}\ar[urrr]^{g\circ f}&&&&&&X^\p\ar@{..>}[r]^{ h\{1\}\circ l}\ar[u]^{l}&Z^\p\{1\}
    }
\end{equation}
\elem

The maps $\psi$ and $\phi$ in (\ref{eq: octahedral}) can be written explicitly as follows. By the long exact sequences in the assumption of Lemma \ref{lem: octahedral lemma}, we know that $Z^\p, X^\p,Y^\p$ are chain homotopic to the mapping cones $\cone(f),\cone(g),\cone(g\circ f)$, respectively. Under such homotopies, we can write$$\begin{aligned}
    \psi:Y\oplus X\{1\}&\to Z\oplus X\{1\}\\
    \psi(y,x)&\mapsto (g(y),x)
\end{aligned}$$and$$\begin{aligned}
    \phi:Z\oplus X\{1\}&\to Z\oplus Y\{1\}\\
    \phi(z,x)&\mapsto (z,f\{1\}(x))
\end{aligned}$$

However, the chain homotopies are not canonical, and hence the maps $\psi$ and $\phi$ are also not canonical. Thus, usually we cannot identify them with other given maps. Fortunately, with an extra $\mathbb{Z}$-grading, we may identify $H(\cone(\phi))$ with $H(\cone(\phi^\p))$ for another map $\phi^\p:Y^\p\to X^\p$.

First, we introduce the following lemma to deal with the projectivity of maps (\textit{i.e.} maps well-defined only up to scalars). Note that the $\mathbb{Z}$-grading in the following lemma is not the homological grading used before.

\blem\label{lem: projectivity}
Suppose $X$ and $Y$ are $\mathbb{Z}$-graded vector spaces and suppose $f,g:X\to Y$ are homogeneous maps with different grading shifts $k_1$ and $k_2$. Then $\cone(f+g)$ is isomorphic to $\cone(c_1f+c_2g)$ for any $c_1,c_2\in \mathbb{C}\backslash\{0\}$. 
\elem
\bpf
For simplicity, we can suppose $k_1=0$ and $k_2=1$. The proof for the general case is similar. For $i,j\in\mathbb{Z}$, we write $(X,i)$ and $(Y,j)$ for grading summands of $X$ and $Y$, respectively. Suppose $T$ is an automorphism of $X\oplus Y$ that acts by $$\frac{c_2^i}{c_1^{i+1}}{\rm Id}\text{ on }(X,i)\text{ and }\frac{c_2^j}{c_1^j}{\rm Id}\text{ on }(Y,j).$$Then $T$ is an isomorphism between $\cone(f+g)$ and $\cone(c_1f+c_2g)$.
\epf

Then we state the lemma that relates the map $\phi$ in Lemma \ref{lem: octahedral lemma} to another map $\phi^\p$ constructed explicitly.
\blem\label{lem: replacing maps}
Suppose $X,Y, Z,X^\p,Y^\p$ are $\mathbb{Z}\oplus\mathbb{Z}_2$-graded vector spaces satisfying the following horizontal exact sequences. 
\begin{equation*}
\xymatrix@R=3ex{
Z\ar[rr]^{h^\p}\ar[dd]^{=}&&Y^\p\ar[rr]^{l^\p}\ar@<-.5ex>[dd]_{\phi} \ar@<.5ex>[dd]^{\phi^{\prime}=a^\p+b^\p}&&X\{1\}\ar[dd]^{f\{1\}=a+b}\\
&&\\
Z\ar[rr]^{\phi\circ h^\p=\phi^\p\circ h^\p}&&X^\p\ar[rr]^{l}&&Y\{1\}
}
\end{equation*}
where the shift $\{1\}$ is for the $\mathbb{Z}_2$-grading. Suppose $\phi:Y^\p\to X^\p$ satisfies the two commutative diagrams and suppose $\phi^\p:Y^\p\to X^\p$ satisfies the left commutative diagram. Suppose $l$ and $l^\p$ are homogeneous with respect to the $\mathbb{Z}$-grading. Suppose $f\{1\}=a+b$ and $\phi^\p=a^\p+b^\p$ are sums of homogeneous maps with different grading shifts with respect to the $\mathbb{Z}$-grading. Moreover, suppose the following diagrams hold up to scalars.

\begin{center}
\begin{minipage}{0.4\textwidth}
	\begin{equation*}
\xymatrix@R=3ex{
Y^\p\ar[rr]^{l^\p} \ar[dd]^{a^{\prime}}&&X\{1\}\ar[dd]^{a}\\
&&\\
X^\p\ar[rr]^{l}&&Y\{1\}
}
\end{equation*}
\end{minipage}
\begin{minipage}{0.4\textwidth}
		\begin{equation*}
\xymatrix@R=3ex{
Y^\p\ar[rr]^{l^\p} \ar[dd]^{b^{\prime}}&&X\{1\}\ar[dd]^{b}\\
&&\\
X^\p\ar[rr]^{l}&&Y\{1\}
}
\end{equation*}
\end{minipage}
\end{center}
Then there is an isomorphism between $H(\cone(\phi))$ and $H(\cone(\phi^\p))$.
\elem
\bpf
Since $\phi$ and $\phi^\p$ share the same domain and codomain, it suffices to show that they have the same rank. Fix inner products on $Y^\p$ and $X^\p$ such that we have orthogonal decompositions$$Y^\p=\im( h^\p)\oplus Y^\pp\aand X^\p=\im(\phi\circ h^\p)\oplus X^\pp.$$By commutativity, we know both $\phi$ and $\phi^\p$ send $\im( h^\p)$ onto $\im(\phi\circ  h^\p)$. Hence if we choose bases with respect to the decompositions so that linear maps are represented by matrices (we use row vectors), then we have $$\phi=\begin{bmatrix} A & B \\ 0 & C \end{bmatrix}\aand \phi^\p=\begin{bmatrix} A^\p & B^\p \\ 0 & C^\p \end{bmatrix},$$where $A=A^\p:\im( h^\p)\to \im(\phi\circ  h^\p)$ has full row rank. Then it suffices to show $C$ and $C^\p$ have the same row rank. 

By the exactness at $Y^\p$ and $X^\p$, we know the restriction of $l^\p$ on $Y^\pp$ is an isomorphism between $Y^\pp$ and $\im(l^\p)$ and the restriction of $l$ on $X^\pp$ is an isomorphism between $X^\pp$ and $\im(l)$. By commutativity, we know that both $a$ and $b$ send $\im(l^\p)$ to $\im(l)$ and $$\operatorname{rowrank}(C)=\operatorname{rank}(f\{1\}|{\im(l^\p)})\aand \operatorname{rowrank}(C^\p)=\operatorname{rank}((c_1a+c_2b)|{\im(l^\p)})$$for some $c_1,c_2\in\cstar$. Since $l$ and $l^\p$ are homogeneous, there exist induced $\mathbb{Z}$-gradings on $\im(l)$ and $\im(l^\p)$. The maps $a$ and $b$ are still homogeneous with different grading shifts with respect to these induced gradings. Then we can apply Lemma \ref{lem: projectivity} to obtain that the ranks of the restrictions of $f\{1\}=a+b$ and $c_1a+c_2b$ on $\im(l^\p)$ are the same.
\epf

\section{Integral surgery formulae}\label{sec: integral surgery formulae}


\subsection{A formula for framed instanton homology}\label{subsec: A strategy to prove the formula}

In this subsection, we propose an integral surgery formula based on sutured instanton homology and package it into the language of bent complexes in a later subsection. 

Suppose $K\subset Y$ is a framed rationally null-homologous knot, and we adopt the notations introduced in Section \ref{sec: preliminaries}. Define $$\pi^{\pm}_{m,k}\deq \Psi_{\pm,m-1+2k}^{m+k}\circ \psi^{\frac{2m+2k-1}{2}}_{\mp,m+k}:\mathbf{\Ga}_{\frac{2m+2k-1}{2}}\ra\mathbf{\Ga}_{m+2k-1}$$
and write $\pi^{\pm,i}_{m,k}$ as the restriction of $\pi^{\pm}_{m,k}$ on $\sutg{\frac{2m+2k-1}{2}}{i}$. From Lemma \ref{lem: bypass n+1,n,2n+1/2} and Lemma \ref{lem: bypass n,n+1,mu}, we can verify that $\pi^{\pm}_{m,k}$ shifts grading by $\pm(mp-q)/2$, and then the integral surgery formula can be stated as follow.

\bthm\label{thm: general mapping cone} Suppose $m$ is a fixed integer such that $mp-q\neq 0$. Then for any sufficiently large integer $k$, there exists an exact triangle
\begin{equation*}\label{eq: generalized mapping cone}
	\xymatrix{
	\mathbf{\Ga}_{\frac{2m+2k-1}{2}}\ar[rr]^{\pi_{m,k}^++\pi_{m,k}^-}&&\mathbf{\Ga}_{m+2k-1}\ar[dl]\\
	&\mathbf{Y}_m\ar[ul]&
	}
\end{equation*}
Hence we have an isomorphism$$\mathbf{Y}_m\cong H(\cone(\pi_{m,k}^++\pi_{m,k}^-)).$$
\ethm
\brem
Let $\mu$ and $\lambda$ represent the meridian and the longitude of the knot $K$, respectively. Then, $mp-q\neq 0$ is equivalent to the fact that $-m\mu+\lambda$ is not isotopic to a connected component of the boundary of the Seifert surface. Specifically, if $K$ is null-homologous, we must have $m\neq 0$.
\erem

In the rest of this subsection and in the next subsection, we state the strategy to prove Theorem \ref{thm: general mapping cone}, and defer the proofs of some propositions to the remaining sections. An important step is to apply the octahedral axiom mentioned in Section \ref{subsec: Algebraic lemmas} to the following diagram:

\begin{equation}\label{eq: intro, octahedron}
	\xymatrix{
	&\sut{\frac{2m+2k-1}{2}}\ar@{..>}[dr]\ar[ddl]&\\
	\dehny{m}\ar@{..>}[ur]\ar[d]&&\sut{m-1+2k}\ar@{..>}[ll]\ar[ddl]\\
	\sut{\mu}\ar[dr]\ar[rr]&&\sut{m-1+k}\oplus \sut{m-1+k}\ar[u]\ar[uul]\\
	&\sut{m-1}\ar[uul]\ar[ur]&
	}
\end{equation}

To obtain the dotted exact triangle, we need to establish the following three exact triangles:
\begin{equation}\label{eq: intro, exact, 1}
	\xymatrix{
	I^{\sharp}(-S^3_{-m}(K))\ar[d]&\\
	\sut{\mu}\ar[dr]&\\
	&\sut{m-1}\ar[uul]
	}
\end{equation}

\begin{equation}\label{eq: intro, exact, 2}
	\xymatrix{
	&\sut{\frac{2m+2k-1}{2}}\ar[ddl]&\\
	&&\\
		\sut{\mu}\ar[rr]&&\sut{m-1+k}\oplus \sut{m-1+k}\ar[uul]
	}
\end{equation}

\begin{equation}\label{eq: intro, exact, 3}
\xymatrix{
	&\sut{m-1+2k}\ar[ddl]\\
	&\sut{m-1+k}\oplus \sut{m-1+k}\ar[u]\\
	\sut{m-1}\ar[ur]&
	}
\end{equation}
and establish the following commutative diagram:
\begin{equation}\label{eq: intro, commute, 1}
	\xymatrix{
	\sut{\mu}\ar[dr]\ar[rr]&&\sut{m-1+k}\oplus \sut{m-1+k}\\
	&\sut{m-1}\ar[ur]&
	}
\end{equation}
The octahedral lemma then implies the existence of the dotted triangle and ensure that all diagrams in (\ref{eq: intro, octahedron}) other than exact triangles commute.

We will then use Lemma \ref{lem: replacing maps} to identify the map coming from the octahedral lemma with $\pi_{m,k}^++\pi_{m,k}^-$. We also require the following two extra diagrams to commute, where the maps other than  $\pi_{m,k}^++\pi_{m,k}^-$ come from (\ref{eq: intro, octahedron}).
\begin{equation}\label{eq: intro, commute, 2}
	\xymatrix{
	\sut{\frac{2m+2k-1}{2}}\ar[dr]^{\pi^+_{m,k}+\pi^-_{m,k}}&\\
	&\sut{m-1+2k}\\
	&\sut{m-1+k}\oplus \sut{m-1+k}\ar[u]\ar[uul]
	}
\end{equation}

\begin{equation}\label{eq: intro, commute, 3}
	\xymatrix{
	&\sut{\frac{2m+2k-1}{2}}\ar[ddl]\ar[dr]^{\pi^+_{m,k}+\pi^-_{m,k}}&\\
	&&\sut{m-1+2k}\ar[ddl]\\
	\quad\quad\sut{\mu}\ar[dr]\quad\quad&&\\
	&\sut{m-1}&
	}
\end{equation}

Indeed, by applying Lemma \ref{lem: replacing maps}, it suffices to prove some weaker commutative diagrams involving $\pi^\pm_{m,k}$ separately.

\subsection{A strategy of the proof}\label{subsec: a strategy}

In this subsection, we provide more details of the strategy mentioned in Section \ref{subsec: A strategy to prove the formula}. For simplicity, we fix the scalar ambiguities of commutative diagrams as in Section \ref{subsec Fixing the scalars}. To write down the maps, we redraw the octahedral diagram (\ref{eq: intro, octahedron}) as follows:

\begin{equation}\label{eq: octahedral 3}
    \xymatrix@C=4ex{
    &&&\dehny{m}\ar@{..>}[dr]^{\psi}\ar[rr]^{}&&\sut{\mu}\ar[ddr]^{\begin{aligned}
        \psi_{+,m-1}^\mu\\+\psi_{-,m-1}^\mu
    \end{aligned}}&
    \\
    &&\sut{m-1}\ar[dr]^{\begin{aligned}
       (\Psi_{+,m-1+k}^{m-1},\\\Psi_{-,m-1+k}^{m-1}) 
    \end{aligned}}\ar[ur]^{}&&\sut{\frac{2m+2k-1}{2}}\ar@{..>}[ddrr]^{\phi}\ar[ur]^{l^\p}&&
    \\
 &&&\sut{m-1+k}\oplus\sut{m-1+k}\ar[drrr]_{\begin{aligned}
     \Psi_{-,m-1+2k}^{m-1+k}\\-\Psi_{+,m-1+2k}^{m-1+k}
 \end{aligned}}\ar[ur]_{h^\p}&&&\sut{m-1}
    \\
    \sut{\mu}\ar[uurr]^{\begin{aligned}
        \psi_{+,m-1}^\mu\\+\psi_{-,m-1}^\mu
    \end{aligned}}\ar[urrr]_{\begin{aligned}
        (\psi_{-,m-1+k}^\mu,\\\psi_{+,m-1+k}^\mu)
    \end{aligned}}&&&&&&\sut{m-1+2k}\ar[u]^{l}
    }
     \end{equation}
where $$h^\p=\psi_{-,\frac{2m+2k-1}{2}}^{m+k-1}-\psi_{+,\frac{2m+2k-1}{2}}^{m+k-1}.$$
   
The reader can compare (\ref{eq: octahedral 3}) with (\ref{eq: octahedral 0}) and (\ref{eq: octahedral}). We omit the term corresponding to $Z^\p\{1\}$ because there is not enough room, and the maps involving it are not important in our proof.    

The first exact sequence of (\ref{eq: octahedral 3})\begin{equation}\label{eq: first sequence} \sut{\mu}\xra{\psi_{+,m-1}^\mu+\psi_{-,m-1}^\mu}\sut{m-1}\xra{}\dehny{m}\xra{}\sut{\mu}\end{equation} follows from the second exact triangle in Lemma \ref{lem: surgery triangles}. Though the map $A_{m-1}$ may not be the same as the sum $\psi_{+,m-1}^\mu+\psi_{-,m-1}^{n}$, we can use the following proposition and Lemma \ref{lem: projectivity} (another special case of Proposition \ref{prop: -1 surgery and bypass}) to identify $\cone(A_{m-1})$ with $\cone(\psi_{+,m-1}^\mu+\psi_{-,m-1}^n)$. Here we use the assumption that $mp-q\neq 0$.
\bprop[]\label{prop: -1 surgery}
Suppose $A_{n-1}$ is the map in Lemma \ref{lem: surgery triangles}. For any integer $n$, there exist scalars $c_1,c_2\in\mathbb{C}\backslash\{0\}$ such that $$A_{n-1}=c_1\psi_{+,n-1}^\mu+c_2\psi_{-,n-1}^\mu.$$
\eprop

The exactness at $$\sut{m-1+k}\oplus\sut{m-1+k}$$in the second and the third exact sequences are both special cases of the following proposition, which will be proved in Section \ref{subsec: the direct summand space} by diagram chasing. 


\bprop[]\label{prop: exactness at direct summand}
Fixing the scalars as in Section \ref{subsec Fixing the scalars}, and given $n\in\mathbb{Z}$ and $k_0\in\mathbb{N}_+$. Then, for any $c_1,c_2,c_3,c_4$ satisfying the equation\begin{equation*}
    c_1c_3=-c_2c_4,
\end{equation*}
the following sequence is exact
$$\sut{n}\xra{(c_1\Psi_{+,n+k_0}^n,c_2\Psi_{-,n+k_0}^n)}\sut{n+k_0}\oplus\sut{n+k_0}\xra{c_3\Psi_{-,n+2k_0}^{n+k_0}+c_4\Psi_{+,n+2k_0}^{n+k_0}}\sut{n+2k_0}$$
\eprop
\brem
The exactness at the direct summand for the second exact sequence (the one involving $\sut{\mu}$) might not be as clear from Proposition \ref{prop: exactness at direct summand}. Explicitly, we apply the proposition to the dual knot $K^\p$ corresponding to $\Ga_{m+k}$ with framing $\lambda^\p=-\mu$ and $n=0,k_0=1$.
\erem

The exactness at $\sut{\mu}$ and $\sut{\frac{2m+2k-1}{2}}$ in the second exact sequence of (\ref{eq: octahedral 3}) \begin{equation}\label{eq: second sequence} \sut{\mu}\xra{(\psi_{-,m-1+k}^{\mu},\psi_{+,m-1+k}^{\mu})}\sut{m-1+k}\oplus\sut{m-1+k}\xra{\psi_{-,\frac{2m+2k-1}{2}}^{m+k-1}-\psi_{+,\frac{2m+2k-1}{2}}^{m+k-1}}\sut{\frac{2m+2k-1}{2}}\xra{l^\p}\sut{\mu}\end{equation}will also be proved by diagram chasing. We can explicitly construct the map $l^\p$ by the composition of bypass maps
    $$l^\p\deq\psi_{-,\mu}^{m+k}\circ\psi_{+,m+k}^{\frac{2m+2k-1}{2}}=\psi_{+,\mu}^{m+k}\circ\psi_{-,m+k}^{\frac{2m+2k-1}{2}},$$where the last equation follows from Lemma \ref{lem: comm diag for n+1 to n} and the conventions in Section \ref{subsec Fixing the scalars}.
The following proposition will be proved in Section \ref{subsec: second triangle} by diagram chasing.

\bprop[]\label{prop: exactness at other two, 2}
Suppose $l^\p$ is constructed as above. For any $c_1,c_2,c_3,c_4\in\cstar$, the following sequence is exact
$$\begin{aligned}
    \sut{m-1+k}\oplus\sut{m-1+k}\xra{c_3\psi_{-,\frac{2m+2k-1}{2}}^{m+k-1}+c_4\psi_{+,\frac{2m+2k-1}{2}}^{m+k-1}}\sut{\frac{2m+2k-1}{2}}\xra{l^\p}\sut{\mu}\\\xra{(c_1\psi_{-,m-1+k}^{\mu},c_2\psi_{+,m-1+k}^{\mu})}\sut{m-1+k}\oplus\sut{m-1+k}
\end{aligned}$$
\eprop
\brem
In the proof of \cite[Theorem 3.23]{LY2021large}, we obtained a long exact sequence $$\sut{\mu}\xra{(\psi_{+,n-1}^\mu,\psi_{-,n-1}^\mu)}\sut{n-1}\oplus\sut{n-1}\xra{}\sut{\frac{2n-1}{2}}\xra{} \sut{\mu}$$by the octahedral lemma. However, we did not know the two maps involving $\sut{\frac{2n-1}{2}}$ explicitly. Thus, the second exact sequence here is stronger than the one from octahedral lemma.
\erem
\brem\label{rem: coefficients}
The reason that Proposition \ref{prop: exactness at other two, 2} holds for any choices of $c_1,c_2,c_3,c_4$ is because $$\ker((c_1\psi_{-,m-1+k}^{\mu},c_2\psi_{+,m-1+k}^{\mu}))=\ker(c_1\psi_{-,m-1+k}^{\mu})\cap\ker(c_2\psi_{+,m-1+k}^{\mu})$$and$$\im(c_3\psi_{-,\frac{2m+2k-1}{2}}^{m+k-1}+c_4\psi_{+,\frac{2m+2k-1}{2}}^{m+k-1})=\im(c_3\psi_{-,\frac{2m+2k-1}{2}}^{m+k-1})+\im(c_4\psi_{+,\frac{2m+2k-1}{2}}^{m+k-1}),$$where the right hand sides of the equations are independent of scalars.
\erem
The exactness at $\sut{m-1}$ and $\sut{m-1+2k}$ in the third exact sequence of (\ref{eq: octahedral 3})\begin{equation}\label{eq: third sequence}
     \sut{m-1}\xra{(\Psi_{+,m-1+k}^{m-1},\Psi_{-,m-1+k}^{m-1})}\sut{m-1+k}\oplus\sut{m-1+k}\xra{\Psi_{-,m-1+2k}^{m-1+k}-\Psi_{+,m-1+2k}^{m-1+k}}\sut{m-1+2k}\xra{l}\sut{m-1}
\end{equation} is harder to prove since the map $l$ cannot be constructed by  bypass maps. We expect that there are many equivalent constructions of $l$ and we will use the one for which the exactness is easiest to prove. Even so, we only prove the exactness with the assumption that $k$ is large. See Section \ref{subsec: constructing Phi^n+2k_n} for more details.
\bprop[]\label{prop: exactness at other two}
Suppose $c_1,c_2,c_3,c_4\in\cstar$ and suppose $k_0$ is a large integer. For any $n\in \intg$, there exists a map $l$ such that the following sequence is exact
$$\sut{n+k_0}\oplus\sut{n+k_0}\xra{c_3\Psi_{-,n+2k_0}^{n+k_0}+c_4\Psi_{+,n+2k_0}^{n+k_0}}\sut{n+2k_0}\xra{l}\sut{n}\xra{(c_1\Psi_{+,n+k_0}^n,c_2\Psi_{-,n+k_0}^n)}\sut{n+k_0}\oplus\sut{n+k_0}$$
\eprop
\brem
 In the first arXiv version of this paper, we only proved Proposition \ref{prop: exactness at other two} for knots in $S^3$ because we had to use the fact that $\dimc I^\sharp(-S^3)=1$ and $S^3$ has an orientation-reversing involution. The construction of $l$ for knots in general $3$-manifolds is inspired by the original proof for $S^3$ and the proof in Section \ref{sec: third triangle} is a generalization of the previous proof.
\erem
\brem
For the same reason as in Remark \ref{rem: coefficients}, the coefficients in Proposition \ref{prop: exactness at other two} are not important.
\erem

Then we consider the commutative diagrams mentioned in Section \ref{subsec: A strategy to prove the formula}. By Lemma \ref{lem: bypass n,n+1,mu} and Lemma \ref{lem: comm diag for n,n+1,mu}, we have  \begin{equation*}\label{eq: commu in octa}
    (\Psi_{+,m-1+k}^{m-1},\Psi_{-,m-1+k}^{m-1})\circ (\psi_{+,m-1}^\mu+\psi_{-,m-1}^\mu)=(\psi_{-,m-1+k}^\mu,\psi_{+,m-1+k}^\mu),
\end{equation*}
which verifies the commutative diagram in the assumption of the octahedral axiom.

Define $$\phi^\p\deq \pi_{m,k}^++\pi_{m,k}^-=\Psi_{+,m-1+2k}^{m+k}\circ \psi^{\frac{2m+2k-1}{2}}_{-,m+k}+\Psi_{-,m-1+2k}^{m+k}\circ \psi^{\frac{2m+2k-1}{2}}_{+,m+k}$$
By Lemma \ref{lem: bypass n+1,n,2n+1/2} and Lemma \ref{lem: comm diag for n+1 to n} with $n=m+k$, we have
$$\begin{aligned}\phi^\p\circ h^\p=&(\Psi_{+,m-1+2k}^{m+k}\circ \psi^{\frac{2m+2k-1}{2}}_{-,m+k}+\Psi_{-,m-1+2k}^{m+k}\circ \psi^{\frac{2m+2k-1}{2}}_{+,m+k})\circ(\psi_{-,\frac{2m+2k-1}{2}}^{m+k-1}-\psi_{+,\frac{2m+2k-1}{2}}^{m+k-1})\\=& \Psi_{-,m-1+2k}^{m+k}\circ \psi^{\frac{2m+2k-1}{2}}_{+,m+k}\circ \psi_{-,\frac{2m+2k-1}{2}}^{m+k-1}-\Psi_{+,m-1+2k}^{m+k}\circ \psi^{\frac{2m+2k-1}{2}}_{-,m+k}\circ \psi_{+,\frac{2m+2k-1}{2}}^{m+k-1}\\=&\Psi_{-,m-1+2k}^{m+k}\circ \psi^{m+k-1}_{-,m+k}-\Psi_{+,m-1+2k}^{m+k}\circ \psi^{m+k-1}_{+,m+k}\\=&\Psi_{-,m-1+2k}^{m-1+k}-\Psi_{+,m-1+2k}^{m-1+k}\end{aligned}$$
This verifies the second commutative diagram mentioned in Section \ref{subsec: A strategy to prove the formula}.

Finally, we state a weaker version of the third commutative diagram mentioned in Section \ref{subsec: A strategy to prove the formula}, which is enough to apply Lemma \ref{lem: replacing maps}. The following proposition will be proved in Section \ref{subsec: commutative diagram}. 

\bprop[]\label{prop: commute 2}
Suppose $l^\p$ and $l$ are the maps in Proposition \ref{prop: exactness at other two, 2} and Proposition \ref{prop: exactness at other two}. Then, there are two commutative diagrams up to scalars.
\begin{center}
\begin{minipage}{0.4\textwidth}
	\begin{equation*}
\xymatrix@R=3ex{
\sut{\frac{2m+2k-1}{2}}\ar[rr]^{l^\p}\ar[dd]^{\pi_{m,k}^+}&&\sut{\mu}\ar[dd]^{\psi_{+,m-1}^\mu}\\
&&\\
\sut{m-1+2k}\ar[rr]^{l}&&\sut{m-1}
}
\end{equation*}
\end{minipage}
\begin{minipage}{0.4\textwidth}
		\begin{equation*}
\xymatrix@R=3ex{
\sut{\frac{2m+2k-1}{2}}\ar[rr]^{l^\p}\ar[dd]^{\pi_{m,k}^-}&&\sut{\mu}\ar[dd]^{\psi_{-,m-1}^\mu}\\
&&\\
\sut{m-1+2k}\ar[rr]^{l}&&\sut{m-1}
}
\end{equation*}
\end{minipage}
\end{center}

\eprop

\bpf[Proof of Theorem \ref{thm: general mapping cone}]
We verified all assumptions of the octahedral lemma (Lemma \ref{lem: octahedral lemma}) for the diagram (\ref{eq: octahedral 3}). Hence, there exists a map $\phi$ such that 
$$\dehny{m}\cong H(\cone(\phi)).$$
We also verified all assumptions of Lemma \ref{lem: replacing maps} for $\phi^\p=\pi_{m,k}^++\pi_{m,k}^-$. Thus, we have 
$$H(\cone(\phi^\p))\cong H(\cone(\phi))\cong \dehny{m}.$$
Then the desired triangle in the theorem holds.
\epf

\subsection{Reformulation by bent complexes}\label{subsec: integral surgery formula for I-sharp}
\quad

In this subsection, we restate Theorem \ref{thm: general mapping cone} using the language of bent complexes introduced in \cite{LY2021large}. Suppose $K$ is a rationally null-homologous knot in a closed 3-manifold $Y$. We continue to adopt the notations and conventions from Section \ref{sec: preliminaries} and Section \ref{subsec Fixing the scalars}.

Putting bypass triangles in Lemma \ref{lem: bypass n,n+1,mu} for different $n$ together, we obtain the following diagram:
\begin{equation}\label{eq: useful triangles}
\xymatrix@R=6ex{
\cdots&\ar[l]\sut{n+1}\ar[dr]_{\psi_{+,\mu}^{n+1}}&&\sut{n}\ar[ll]_{\psi_{+,n+1}^{n}}\ar[dr]_{\psi_{+,\mu}^n}&&\sut{n-1}\ar[ll]_{\psi_{+,n}^{n-1}}\ar[dr]_{\psi_{+,\mu}^{n-1}}&&\sut{n-2}\ar[ll]_{\psi_{+,n-1}^{n-2}}&\ar[l]\cdots
\\
&\cdots&\sut{\mu}\ar[ur]^{\psi_{+,n}^\mu}\ar[dl]^{\psi_{-,n-2}^\mu}&&\sut{\mu}\ar[ur]^{\psi_{+,n-1}^\mu}\ar[dl]^{\psi_{-,n-1}^\mu}&&\sut{\mu}\ar[ur]^{\psi_{+,n-2}^\mu}\ar[dl]^{\psi_{-,n}^\mu}&\cdots&\\
\cdots\ar[r]&\sut{n-2}\ar[rr]_{\psi_{-,n-1}^{n-2}}&&\sut{n-1}\ar[rr]_{\psi_{-,n}^{n-1}}\ar[ul]_{\psi_{-,\mu}^{n-1}}&&\sut{n}\ar[rr]_{\psi_{-,n+1}^n}\ar[ul]_{\psi_{-,\mu}^n}&&\sut{n+1}\ar[ul]_{\psi_{-,\mu}^{n+1}}\ar[r]&\cdots
}
\end{equation}
where the $\mathbb{Z}$-grading shift of $\psi_{\pm,k}^\mu\circ\psi_{\pm,\mu}^k$ is $\pm p$ for any $k\in\mathbb{Z}$. From (\ref{eq: useful triangles}), we constructed in \cite[Section 3.4]{LY2021large} two spectral sequences $\{E_{r,+},d_{r,+}\}_{r\ge 1}$ and $\{E_{r,-},d_{r,-}\}_{r\ge 1}$ from $\sut{\mu}$ to $\dehny{}$, where $d_{r,\pm}$ is roughly \begin{equation}\label{eq: dr}
   \psi_{\pm,\mu}^{k}\circ(\Psi_{\pm,k+r}^k)^{-1}\circ\psi_{\pm,k+r}^\mu\text{ for any }k\in\mathbb{Z}.
\end{equation}The composition with the inverse map is well-defined on the $r$-th page, and the independence of $k$ (and hence $n$ in (\ref{eq: dr})) follows from Lemma \ref{lem: comm diag for n,n+1,mu}. The $\mathbb{Z}$-grading shift of $d_{r,\pm}$ is $\pm rp$. By fixing an inner product on $\sut{\mu}$, we then lifted those spectral sequences to two differentials $d_+$ and $d_-$ on $\sut{\mu}$ such that $$H(\sut{\mu},d_+)\cong H(\sut{\mu},d_-)\cong \dehny{}.$$In such way, the inverses of $\Psi_{\pm,k+r}^k$ are also well-defined, which we will use freely later.

Then we propose an integral surgery formula for $\dehny{m}$ using differentials $d_+$ and $d_-$ on $\sut{\mu}$. To state the formula, we introduce the following notations.
\bdefn[{\cite[Construction 3.27 and Definition 5.12]{LY2021large}}]\label{defn: complex B}
For any integer $s$, define the complexes\[B^\pm(s)\deq(\bigoplus_{k\in\mathbb{Z}}(\sut{\mu},s+kp),d_\pm),\quad B^+({\ge s})\deq (\bigoplus_{k\ge 0}(\sut{\mu},s+kp),d_+),\]\[\aand B^-({\le s})\deq (\bigoplus_{k\le 0}(\sut{\mu},s+kp),d_-).\]
Furthermore, define $$I^+(s):B^+(\ge s)\to B^+(s)\aand I^-(s):B^-(\le s)\to B^-(s)$$to be the inclusion maps. We also write the same notation for the induced map on homology.
\edefn
\brem\label{rem: total complex}
By Lemma \ref{lem: vanishing grading}, we know that the nontrivial gradings of $\sut{\mu}$ are finite. Then, for any sufficiently large integer $s_0$ satisfying $$s-s_0p \le -\frac{p-\chi(S)}{2}\aand s+s_0p \ge \frac{p-\chi(S)}{2},$$we have \begin{equation*}\label{eq: total complex}
    B^+(s)=B^+(\ge s-s_0p)\aand B^-(s)=B^-({\le s+s_0p}).
\end{equation*}In such case, $I^+(s-s_0p)$ and $I^-(s+s_0p)$ are identities.
\erem
By splitting the diagram (\ref{eq: useful triangles}) into $\mathbb{Z}$-gradings, we can calculate homologies of the complexes defined in Definition \ref{defn: complex B}.
\bprop\label{prop: subcomplex convergence}
Suppose $n\in\mathbb{N}_+$ and $i$ is a grading. Fix an inner product on $\sut{n}$. If $i>(p-\chi(S))/2-np$, then there exists a canonical isomorphism $$H(B^+({\ge i}))\cong (\sut{n},i+\frac{(n-1)p-q}{2}).$$If $i<-(p-\chi(S))/2+np$, then there exists a canonical isomorphism $$H(B^-({\le i}))\cong (\sut{n},i-\frac{(n-1)p-q}{2}).$$
\eprop
\bpf
The proof mirrors that of \cite[Lemma 5.13]{LY2021large}. Following the notation in \cite[(3.9) and (3.10)]{LY2021large}, if $$i>\hat{i}_{max}^\mu-nq=(p-\chi(S))/2-np,$$then $\sut{0}^{i,+}=0$ (the corresponding grading summand of $\sut{0}$) and the isomorphism follows from the convergence theorem of the unrolled spectral sequence \cite[Theorem 2.4]{LY2021large} (see also \cite[Theorem 6.1]{Boardman99}). Note that the unrolled spectral sequence induces a filtration on $\sut{n}$, and the homology is canonically isomorphic to the direct sum of all associated graded objects of the filtration. Then we use the inner product to identify the direct sum with the total space $\sut{n}$. The other statement holds for the same reason.
\epf

\bdefn[{\cite[Construction 3.27 and Definition 5.12]{LY2021large}}]\label{defn: complex A}
For any integer $s$, define the \textbf{bent complex}\[A(s)\deq (\bigoplus_{k\in\mathbb{Z}}(\sut{\mu},s+kp),d_s),\]where for any element $x\in (\sut{\mu},s+kp)$, 
\[
d_s(x)=\begin{cases}
d_+(x)&k>0,\\
d_+(x)+d_-(x)&k=0,\\
d_-(x)&k<0.
\end{cases}\]
Define\[\pi^+({s}):A(s)\to B^+(s)\aand \pi^-({s}):A(s)\to B^-(s)\]by
\[\pi^+(s)(x)=\begin{cases}
x&k\ge 0,\\
0&k< 0,\end{cases}\aand \pi^-(s)(x)=\begin{cases}
x&k\le 0,\\
0&k>0,\end{cases}\]where $x\in (\sut{\mu},s+kp)$. Define $$\pi^\pm:\bigoplus_{s\in\mathbb{Z}}A(s)\to \bigoplus_{s\in\mathbb{Z}}B^\pm(s)$$by putting $\pi^\pm(s)$ together for all $s$ . We also use the same notation for the induced map on homology.
\edefn
\brem\label{rem: iso when large}
Similar to Remark \ref{rem: total complex}, according to Lemma \ref{lem: vanishing grading}, there are only finitely many the nontrivial gradings of $\sut{\mu}$. Then, for any sufficiently large integer $s_0$ such that $s_0\ge (p-\chi(S))/2$, we have \begin{equation*}\label{eq: total complex A}
A(s_0)=B^-(s_0)\aand A(-s_0)=B^+(-s_0).
\end{equation*}
In such case, $\pi^-(s_0)$ and $\pi^+(-s_0)$ are identities.
\erem
Now, we state the integral surgery formula in the above setup.

\bthm\label{thm: integral bent complex}
Suppose $m$ is a fixed integer such that $mp-q\neq 0$. Then there exists an isomorphism $$\Xi_{m}:\bigoplus_{s\in\mathbb{Z}}H(B^+(s))\xra{\cong}\bigoplus_{s\in\mathbb{Z}}H(B^-(s+mp-q))$$as the direct sum of isomorphisms$$\Xi_{m,s}:H(B^+_s)\xra{\cong}H(B^-_{s+mp-q})$$ so that $$\dehny{m}\cong H\bigg(\cone(\pi^-+\Xi_m\circ \pi^+:\bigoplus_{s\in\mathbb{Z}}H(A(s))\to  \bigoplus_{s\in\mathbb{Z}}H(B^-(s)))\bigg).$$ 
\ethm

\bpf
According to Remark \ref{rem: iso when large}, we only need to consider the maps $\pi^\pm(s)$ for $|s|$ less than a fixed integer. For such values of $s$, we can apply the following proposition.

\bprop[{\cite[Proposition 3.28]{LY2021large}}]\label{prop: identify cone}
Fix $m,s\in\mathbb{Z}$ such that $|s|\le (p-\chi(S))/2$. For any large integer $k$, fix inner products on $\sut{\frac{2m+2k-1}{2}}$ and $\sut{m-1+2k}$. Then there exist $s_1,s_2^+,s_2^-,s_3^+,s_3^-\in\mathbb{Z}$ such that the following diagram commutes

\begin{equation*}
\xymatrix@R=3ex{
H(A(s))\ar[rr]^{\pi^\pm(s)}\ar[dd]^{\cong}&&H(B^\pm(s))\ar[dd]^{\cong}
\\\\
(\sut{\frac{2m+2k-1}{2}},s_1)\ar[rr]^{\pi^{\pm,s_1}_{m,k}}&&(\sut{m-1+2k},s_3^\pm)
}
\end{equation*}
where the maps $\pi^{\pm,s_1}_{m,k}$, defined in Section \ref{subsec: A strategy to prove the formula}, factor through $\sutg{m+k}{s_2^\pm}$.
\eprop
\brem
The maps $\pi^\pm(s)$ factor through $I^\pm(s)$ constructed in Definition \ref{defn: complex B}. We denote $$\pi^\pm(s)=I^+(s)\circ \pi^{\pm,\prime}(s).$$This corresponds to the factorization about $(\sut{m+k},s_2^\pm)$ in Proposition \ref{prop: identify cone} (we fix an inner product on $\sut{m+k}$ to apply Proposition \ref{prop: subcomplex convergence}), \textit{i.e.}, the following diagrams commute
\begin{equation*}
\xymatrix@R=3ex{
H(A(s))\ar[rr]^{\pi^{+,\prime}(s)}\ar[dd]^{\cong}&&H(B^+(\ge s))\ar[dd]^{\cong}\ar[rr]^{I^+(s)}&&H(B^+(s))\ar[dd]^{\cong}
\\\\
(\sut{\frac{2m+2k-1}{2}},s_1)\ar[rr]^{\psi_{-,m+k}^{\frac{2m+2k-1}{2}}}&&(\sut{m+k},s_2^+)\ar[rr]^{\Psi_{+,m-1+2k}^{m+k}}&&(\sut{m-1+2k},s_3^+)
}
\end{equation*}
\begin{equation*}
\xymatrix@R=3ex{
H(A(s))\ar[rr]^{\pi^{-,\prime}(s)}\ar[dd]^{\cong}&&H(B^-(\le s))\ar[dd]^{\cong}\ar[rr]^{I^-(s)}&&H(B^-(s))\ar[dd]^{\cong}
\\\\
(\sut{\frac{2m+2k-1}{2}},s_1)\ar[rr]^{\psi_{+,m+k}^{\frac{2m+2k-1}{2}}}&&(\sut{m+k},s_2^-)\ar[rr]^{\Psi_{-,m-1+2k}^{m+k}}&&(\sut{m-1+2k},s_3^-)
}
\end{equation*}
\erem
From the calculation in \cite[Remark 3.29]{LY2021large} (we replace $n$ and $l$ there by $m+k$ and $k-1$, and note that there is a typo about sign in the first arXiv version of \cite{LY2021large}), the difference of the grading shifts is $$s_3^+-s_3^-=(m+k-(k-1)-1)p-q=mp-q.$$
Note that the notations in this paper and \cite{LY2021large} are different (\textit{c.f.} Remark \ref{rem: notations}).

Then we can construct the isomorphism$$\Xi_{m,s}:H(B^+(s))\xra{\cong}H(B^-(s+mp-q))$$ for $|s|\le (p-\chi(S))/2$ by identifying both $H(B^+(s))$ and $H(B^-(s+mp-q))$ with $(\sut{m-1+2k},s_3^+)$ for a sufficiently large $k$. \textit{A priori}, this isomorphism depends on inner products on $$\sut{\mu},\sut{\frac{2m+2k-1}{2}},\sut{m-1+2k}\aand \sut{m+k}.$$ For other $s$, we can take any isomorphism $\Xi_{m,s}$ since the choice does not affect the computation of the mapping cone.

Consequently, we obtain 
$$H(\cone(\pi^-+\Xi_m\circ \pi^+))\cong H(\cone(\pi_{m,k}^-+\pi_{m,k}^+))\cong \dehny{m},$$
where the last isomorphism comes from Theorem \ref{eq: generalized mapping cone}.
\epf

\brem
Theorem \ref{thm: integral bent complex} is slightly weaker than Theorem \ref{eq: generalized mapping cone}. Indeed, when we use the integral surgery formula to calculate surgeries on the Boromean knot in the companion paper \cite{LY2022integral2}, we have to study the $H_1(Y)$ action on sutured instanton homology, where $Y=\#^{2n}S^1\times S^2$ is the ambient manifold of the knot. This action vanishes on $\sut{\mu}$ so vanishes on the bent complex. But it is nonvanishing on $\sut{m+k}$ and $\sut{m-1+2k}$ and we use this information to realize the computation. This issue for the bent complex might be resolved by introducing some $E_0$-pages for differentials $d_+$ and $d_-$ such that the action is nontrivial on $E_0$-pages.
\erem

\subsection{A formula for instanton knot homology}\label{subsec: An integral surgery formula for instanton knot homology}

The third exact sequence (\ref{eq: third sequence}) implies 
$$\sut{m-1}\cong H(\cone(\Psi_{-,m-1+2k}^{m-1+k}-\Psi_{+,m-1+2k}^{m-1+k}:\sut{m-1+k}\op\sut{m-1+k}\to \sut{m-1+2k}))$$
for any sufficiently large integer $k$. Since there are two copies $\sut{m-1+k}$, we can always regard the grading shifts of the maps $\Psi_{-,m-1+2k}^{m-1+k}$ as different ones by rescaling the grading of the first summand from $i$ to $2i-1$ and the second summand from $i$ to $2i$. Hence we do not need the assumption $mp-q\neq 0$ as in the previous mapping cone formula in Theorem \ref{thm: general mapping cone}. By Lemma \ref{lem: projectivity}, we can replace the minus sign with any coefficient. 

In this subsection, we restate this result in the language of bent complexes. The formula is inspired by Eftekhary's formula for knot Floer homology $\widehat{HFK}$ \cite[Proposition 1.5]{eft18bordered} (see also Hedden-Levine's work \cite{hedden21surgery}). Since $m$ can be any integer, we replace $m-1$ by $m$.
\bthm\label{thm: mapping cone for KHI}
Suppose $m,j\in\mathbb{Z}$. Define $$j^+= j-\frac{(m-1)p-q}{2}\aand j^-= j+\frac{(m-1)p-q}{2}.$$Then there exists an isomorphism $$\Xi_{m,j}^\p:H(B^+(j^+))\xra{\cong}H(B^-(j^-))$$such that$$(\sut{m},j)\cong H\bigg(\cone(I^-(j^-)+\Xi_{m,j}^\p\circ I^+(j^+):H(B^-(\le j^-))\oplus H(B^+(\ge j^+))\to H(B^-(j^-)))\bigg).$$ 

\ethm
\bpf
As mentioned before, we have$$\sut{m}\cong H(\cone(\Psi_{-,m+2k}^{m+k}-\Psi_{+,m+2k}^{m+k}))\cong H(\cone(\Psi_{-,m+2k}^{m+k}+\Psi_{+,m+2k}^{m+k}))$$for any sufficiently large integer $k$.

Since bypass maps are homogeneous, the above mapping cone splits into $\mathbb{Z}$-gradings (or $(\mathbb{Z}+\frac{1}{2})$-gradings). Hence we can use it to calculate $\sutg{m}{j}$. By Lemma \ref{lem: bypass n,n+1,mu}, the corresponding spaces are $$(\sut{m+k},j-\frac{kp}{2})\oplus(\sut{m+k},j+\frac{kp}{2})\aand  (\sut{m+2k},j).$$
From Proposition \ref{prop: subcomplex convergence} with $i=j\pm kp/2$, by fixing an inner product on $\sut{m+k}$, we know that $$(\sut{m+k},j-\frac{kp}{2})\cong H(B^-(\le j^-))\text{ for }j-\frac{kp}{2}<-\frac{p-\chi(S)}{2}+(m+k)p$$and$$(\sut{m+k},j+\frac{kp}{2})\cong H(B^+(\ge j^+))\text{ for }j+\frac{kp}{2}>\frac{p-\chi(S)}{2}-(m+k)p.$$Since $m$ is fixed, when $k$ is sufficiently large, we know that any $j$ with $\sutg{m}{j}$ nontrivial (i.e. $|j|\le (|mp-q|-\chi(S))/2$ by Lemma \ref{lem: vanishing grading}) satisfies the above inequalities. By Proposition \ref{prop: subcomplex convergence} again (fixing an inner product on $\sut{m+2k}$) and Remark \ref{rem: total complex}, for $k$ sufficiently large, we know that $$(\sut{m+2k},j)\cong H(B^-(j^-))\cong H(B^+(j^+))$$for such $j$. By unpackaging the construction of differentials $d_+$ and $d_-$ in \cite[Section 3.4]{LY2021large}, we know that the restrictions of maps $\Psi_{-,m+2k}^{m}$ and $\Psi_{+,m+2k}^m$ on the corresponding gradings coincide with the maps induced by the inclusions $I^-(j^-)$ and $I^+(j^+)$ under the canonical isomorphisms, respectively.


For $|j|\le (|mp-q|-\chi(S))/2$, let $$\Xi^\p_{m,j}:H(B^+( j^+))\xra{\cong}H(B^-(j^-))$$ be the isomorphism obtained from identifying both spaces to the corresponding grading summand of $\sut{m+2k}$. Note that it depends on inner products on $\sut{\mu},\sut{m+k}$ and $\sut{m+2k}$. For other $j$, we can take any isomorphism $\Xi_{m,j}^\p$ since the choice does not affect the computation of the mapping cone. Then we know that $$\begin{aligned}\sutg{m}{j}\cong & H(\cone(\Psi_{-,m+2k}^{m+k}+\Psi_{+,m+2k}^{m+k}|(\sut{m+k},j+\frac{kp}{2})\oplus(\sut{m+k},j-\frac{kp}{2})))\\\cong & H(\cone(I^-(j^-)+\Xi_{m,j}^\p\circ I^+(j^+))).\end{aligned}$$
\epf


\section{Dehn surgery and bypass maps}\label{sec: (-1)-surgery and bypass}

In this section, we prove a generalization of Lemma \ref{lem: -1 surgery and bypass} and Proposition \ref{prop: -1 surgery}.

Suppose $(M,\ga)$ is a balanced sutured manifold and $\al\subset\partial M$ is a connected simple closed curve that intersects the suture $\ga$ twice. There are two natural bypass arcs associated to $\al$, each of which intersects the suture at three points and induces a bypass triangle (\textit{c.f.} \cite[Section 4]{baldwin2018khovanov})
\begin{equation*}
	\xymatrix{
	\shi(-M,-\ga)\ar[rr]^{\psi_{\pm}}&&\shi(-M,-\ga_{2})\ar[dl]\\
	&\shi(-M,-\ga_{3})\ar[ul]&
	}
\end{equation*}
where $\ga_{2}$ and $\ga_{3}$ are the sutures coming from bypass attachments. Note that the two bypass exact triangles involve the same set of balanced sutured manifolds but have different maps between them.
Let $(M_0,\ga_0)$ be obtained from $(M,\ga)$ by attaching a contact $2$-handle along $\al$. From \cite[Section 3.3]{baldwin2016instanton}, it has been shown that a closure of $(-M_0,-\ga_0)$ coincides with a closure of the sutured manifold obtained from $(-M,-\ga)$ by $0$-surgery along $\al$ with respect to the surface framing. Hence there is also a surgery exact triangle (\textit{c.f.} \cite[Lemma 3.21]{LY2020})
\begin{equation*}
	\xymatrix{
	\shi(-M,-\ga)\ar[rr]^{H_{\al}}&&\shi(-M,-\ga_{2})\ar[dl]\\
	&\shi(-M_0,-\ga_0)\ar[ul]&
	}
\end{equation*}
The map $H_{\al}$ is related to the bypass maps $\psi_{\pm}$ as follows:
\bprop\label{prop: -1 surgery and bypass}
There exist $c_1,c_2\in\cstar$, such that
$$H_{\al}=c_1 \psi_+ + c_2 \psi_-.$$
\eprop

\brem
The proof of Proposition \ref{prop: -1 surgery and bypass} was developed through the discussions with John A. Baldwin and Steven Sivek.  
\erem

\bpf[Proof of Proposition \ref{prop: -1 surgery and bypass}]
Let $A\subset \partial M$ be a tubular neighborhood of $\al\subset\partial M$. Push the interior of $A$ into the interior of $M$ to make it a properly embedded surface. By a standard argument in \cite{honda2000classification}, we can assume that a collar of $\partial M$ is equipped with a product contact structure such that $\ga$ is (isotopic to) the dividing set, $\al$ is a Legendrian curve, $A$ is in the contact collar, and $A$ is a convex surface with Legendrian boundary that separates a standard contact neighborhood of $\al$ off $M$. The convex decomposition of $M$ along $A$ yields two pieces
$$M=M'\mathop{\cup}_AV,$$
where $M'$ is diffeomorphic to $M$ and $V$ is the contact neighborhood of $\al$. It is straightforward to check that, after rounding the corners, the contact structure near the boundary of $M'$ is still a product contact structure with $\partial M'$ being a convex boundary. Let $\ga'$ be the dividing set on $\partial M'$. Also, after rounding the corners, with the contact structure on $V\cong S^1\times D^2$, we can suppose $\partial V$ is a convex surface with dividing set being the union of two connected simple closed curves on $\partial V$ of slope $-1$. When viewing $V$ as the complement of an unknot in $S^3$, the dividing set coincides with the suture $\Gamma_1\subset V$, so from now on we call it $\Gamma_1$. By the construction of the gluing map in \cite{li2018gluing}, there exists a map
$$G_1:\shi(-M',-\ga')\otimes \shi(-V,-\Ga_1)\ra \shi(-M,-\ga).$$
As in \cite{li2018gluing}, the map $G_1$ comes from attaching contact handles to $(M',\ga')\sqcup (V,\Gamma_1)$ to recover the gluing along $A$. From \cite[Proposition 1.4]{li2019direct}, we know that
$$\shi(-V,-\Ga_1)\cong\mathbb{C}.$$
Note that $M'$ and $M$ are equipped with the product contact structure near their boundaries. From the functoriality of the contact gluing map in \cite{li2018gluing}, we know that $G_1$ is an isomorphism. Now both the $(-1)$-surgery along a push off of $\al$ and the bypass attachments can be thought of as happening in the piece $V$. Note that the result of both $(-1)$-surgery and the bypass attachments for $\Ga_1$ is $\Ga_2$. Hence we have the following commutative diagram.
\begin{equation}\label{eq: comm diag for gluing and surgery}
	\xymatrix{
	\shi(-M',-\ga')\otimes \shi(-V,-\Ga_1)\ar[rr]^{\quad\quad\quad\quad\quad G_1}\ar[dd]^{}_{\operatorname{Id}\otimes\widehat{H}_{\al}}&&\shi(-M,-\ga)\ar[dd]^{}_{H_{\al}}\\
	&&\\
	\shi(-M',-\ga')\otimes \shi(-V,-\Ga_2)\ar[rr]^{\quad\quad\quad\quad\quad G_2}&&\shi(-M,-\ga_2)
	}
\end{equation}
where $\widehat{H}_\al$ denotes the surgery map for the manifold $V$ and $G_2$ is the gluing map obtained by attaching the same set of contact handles as $G_1$. A similar commutative diagram holds when replacing $H_\al$ and $\widehat{H}_\al$ by $\psi_\pm$ and 
$$\hat{\psi}_{\pm}:\shi(-V,-\Ga_1)\ra \shi(-V,-\Ga_2)$$in (\ref{eq: comm diag for gluing and surgery}), respectively.

Since $G_1$ is an isomorphism, to obtain a relation between $H_{\al}$ and $\psi_{\pm}$, it suffices to understand the relation between $\widehat{H}_{\al}$ and $\hat{\psi}_{\pm}$. From \cite[Proposition 1.4]{li2019direct}, we know that
$$\shi(-V,-\Ga_2)\cong\mathbb{C}^2.$$
Moreover, the meridian disk of $V$ induces a $(\intg+\frac{1}{2})$ grading on $\shi(-V,-\Ga_2)$ and we have
$$\shi(-V,-\Ga_2)\cong \shi(-V,-\Ga_2,\frac{1}{2})\oplus \shi(-V,-\Ga_2,-\frac{1}{2}),$$
with
$$\shi(-V,-\Ga_2,\frac{1}{2})\cong \shi(-V,-\Ga_2,-\frac{1}{2})\cong\mathbb{C}.$$
Let $$\mathbf{1}\in \shi(-V,-\Ga_1)\cong\mathbb{C}$$be a generator. In \cite[Section 4.3]{li2019direct} it is shown that
$$\hat{\psi}_{-}(\mathbf{1})\in \shi(-V,-\Ga_2,\frac{1}{2}){\rm~and~}\hat{\psi}_{+}(\mathbf{1})\in \shi(-V,-\Ga_2,-\frac{1}{2})$$are non-zero.
Also, when viewing $V$ as the complement of the unknot $U$, there is an exact triangle
\begin{equation}\label{eq: unknot surgery}
\xymatrix{
\shi(-V,-\Ga_1)\ar[rr]^{\widehat{H}_{\al}}&&\shi(-V,-\Ga_2)\ar[dl]^{F_{2}}\\
&I^{\sharp}(-S^3)\ar[ul]^{G_1}&
}	
\end{equation}
as in Lemma \ref{lem: surgery triangles}. Comparing the dimensions of the spaces in (\ref{eq: unknot surgery}), we have $G_1=0$ and $\widehat{H}_\al$ is injective. From the fact that $\tau_I(U)=0$, we know from \cite[Corollary 3.5]{li2019tau} that
$$F_2\big|{\shi(-V,-\Ga_2,\frac{1}{2})}\neq0{\rm~and~}F_2\big|{\shi(-V,-\Ga_2,-\frac{1}{2})}\neq0,$$
By the exactness in (\ref{eq: unknot surgery}), we have $\ker( F_2)=\im (\widehat{H}_{\alpha})$ and then $\widehat{H}_\al(\mathbf{1})$ is not in $\shi(-V,-\Ga_2,\pm\frac{1}{2})$, \textit{i.e.}, it is a linear combination of generators of $\shi(-V,-\Ga_2,\pm\frac{1}{2})$. Hence we know that there are $c_1,c_2\in\cstar$ such that
$$\widehat{H}_{\al}(\mathbf{1})=c_1\hat{\psi}_+(\mathbf{1})+c_2\hat{\psi}_-(\mathbf{1}).$$
Then the proposition follows from the commutative diagram (\ref{eq: comm diag for gluing and surgery}).
\epf

In Remark \ref{rem: scalar ambiguity}, we discussed the ambiguity arising from scalars. It is worth mentioning that such ambiguity already exists in instanton theory. For example, if $M$ is the complement of a knot $K\subset S^3$ and $\ga$ consists of two meridians of the knot, which we denote by $\Gamma_{\mu}$, we can choose $\alpha$ to be a curve on $\partial (S^3\backslash N(K))$ of slope $-n$. Then we have a surgery triangle:
\begin{equation*}
	\xymatrix{
	\shi(-M,-\Ga_{\mu})\ar[rr]^{H_{n}}&&\shi(-M,-\Ga_{n-1})\ar[dl]\\
	&I^{\sharp}(-S^3_{-n}(K))\ar[ul]&
	}
\end{equation*}
Note that this triangle is not the one from Floer's original exact triangle, but rather one with a slight modification on the choice of 1-cycles inside the $3$-manifold that represents the second Stiefel-Whitney class of the relevant  $SO(3)$-bundle; see \cite[Section 2.2]{baldwin2020concordance} for more details. Floer's original exact triangle, on the other hand, yields a different triangle
\begin{equation*}
	\xymatrix{
	\shi(-S^3\backslash N(K),-\Ga_{\mu})\ar[rr]^{H'_{n}}&&\shi(-S^3\backslash N(K),-\Ga_{n-1})\ar[dl]\\
	&I^{\sharp}(-S^3_{-n}(K),\mu)\ar[ul]&
	}
\end{equation*}
where $\mu\subset -S^3_{-n}(K)$ denotes a meridian of the knot. Note the difference between $H_{\al}$ and $H'_{\al}$ is that they come from the same cobordism but the $SO(3)$-bundles over the cobordism are different. The local argument to prove Proposition \ref{prop: -1 surgery and bypass} works for both $H_{\al}$ and $H'_{\al}$. Hence there exists non-zero complex numbers $c_1,c_2,c_1',c_2'$ such that
$$H_{\al}=c_1\psp{\mu}{n}+c_2\psm{\mu}{n-1}~{\rm and~}H'_{\al}=c'_1\psp{\mu}{n}+c'_2\psm{\mu}{n-1}$$
where the maps
$$\psi^{\mu}_{\pm,n-1}:\shi(-S^3\backslash N(K),-\Ga_{\mu})\ra \shi(-S^3\backslash N(K),-\Ga_{n-1})$$
are the two related bypass maps. When $n\neq 0$, these two bypass maps have different grading shifting behavior, so by Lemma \ref{lem: projectivity}, different choice of non-zero coefficients does not change the dimensions of kernel and cokernel of the map. Hence we conclude that for $n\neq0$,
$$I^{\sharp}(-S^3_{-n}(K),\mu)\cong I^{\sharp}(-S^3_{-n}(K)).$$
However, when $n=0$, the two bypass maps $\psi^{\mu}_{\pm,n-1}$ both preserves gradings, making the coefficients significant, \textit{i.e.}, $I^{\sharp}(-S^3_{0}(K),\mu)$ and $I^{\sharp}(-S^3_{0}(K))$ might have different dimensions for different choices of coefficients. Indeed, it is observed by Baldwin-Sivek \cite{baldwin2020concordance} that for what they called as W-shaped knots (which is clearly a non-empty class, \textit{e.g.} the figure-8 knot \cite[Proposition 10.4]{baldwin2022concordance2}), these two framed instanton homologies have dimensions differing by $2$.

\section{Some exactness by diagram chasing}\label{sec: Surgery exact triangles and commutative diagrams}

\subsection{At the direct summand}\label{subsec: the direct summand space}

In this subsection, we prove Proposition \ref{prop: exactness at direct summand} by diagram chasing. We restate the result in Proposition \ref{prop2: exactness at direct summand}. We also adopt the conventions for scalars from Section \ref{subsec Fixing the scalars}, and this together with Lemma \ref{lem: com diag for n,n+1,n+2} implies that
\[\Psi_{+,n+2k_0}^{n+k_0}\circ\Psi_{-,n+k_0}^n=\Psi_{-,n+2k_0}^{n+k_0}\circ\Psi_{+,n+k_0}^n.
\]
for any $n$ and $k_0$.

\bprop[]\label{prop2: exactness at direct summand}
Given $n\in\mathbb{Z}$ and $k_0\in\mathbb{N}_+$, for any $c_1,c_2,c_3,c_4$ satisfying the equation\begin{equation*}\label{eq: setting}
    c_1c_3=-c_2c_4,
\end{equation*}
the following sequence is exact
$$\sut{n}\xra{(c_1\Psi_{+,n+k_0}^n,c_2\Psi_{-,n+k_0}^n)}\sut{n+k_0}\oplus\sut{n+k_0}\xra{c_3\Psi_{-,n+2k_0}^{n+k_0}+c_4\Psi_{+,n+2k_0}^{n+k_0}}\sut{n+2k_0}$$
\eprop
\bpf
For simplicity, we only prove the proposition for $n=0$. The proof for any general $n$ is similar (replacing all $\sut{m}$ below by $\sut{n+m}$ and modifying the notations for bypass maps). Also, we only prove the case when
$$c_1=c_2=c_3=1,c_4=-1.$$
The proof for general scalars can be obtained similarly. 

We prove the proposition by induction on $k_0$. We will use the exactness in Lemma \ref{lem: bypass n,n+1,mu} and the commutative diagrams in Lemma \ref{lem: comm diag for n,n+1,mu} and Lemma \ref{lem: com diag for n,n+1,n+2} for many times. For simplicity, we will use them without mentioning the lemmas.

First, we assume $k_0=1$. The proposition reduces to $$\ker(\psi_{-,2}^1-\psi_{+,2}^1)=\im((\psi_{+,1}^0,\psi_{-,1}^0)).$$The commutative diagram in Lemma \ref{lem: com diag for n,n+1,n+2} implies $$\ker(\psi_{-,2}^1-\psi_{+,2}^1)\supset\im((\psi_{+,1}^0,\psi_{-,1}^0)).$$
We then prove $$\ker(\psi_{-,2}^1-\psi_{+,2}^1)\subset\im((\psi_{+,1}^0,\psi_{-,1}^0)).$$Suppose $$(x_1,x_2)\in \ker(\psi_{-,2}^1-\psi_{+,2}^1),\text{\textit{i.e.}},\psi_{-,2}^1(x_1)-\psi_{+,2}^1(x_2)=0.$$Then we have $$\psi_{+,\mu}^1(x_1)=\psi_{+,\mu}^2\circ\psi_{-,2}^1(x_1)= \psi_{+,\mu}^2\circ\psi_{+,2}^1(x_2)=0.$$ By exactness, there exists $y\in \sut{0}$ such that $\psi_{+,1}^0(y)=x_1$. Then $$\psi_{+,2}^1\circ\psi_{-,1}^0(y)=\psi_{-,2}^1\circ\psi_{+,1}^0(y)=\psi_{-,2}^1(x_1)\aand \psi_{+,2}^1(x_2-\psi_{-,1}^0(y))=0.$$By exactness, there exists $z\in\sut{\mu}$ such that $$\psi_{+,1}^\mu(z)=x_2-\psi_{-,1}^0(y).$$
Let $y^\p=y+\psi_{+,0}^\mu(z).$ Then $$\psi_{+,1}^0(y^\p)=\psi_{+,1}^0(y)=x_1$$and$$\psi_{-,1}^0(y^\p)=\psi_{-,1}^0(y)+\psi_{-,1}^0\circ\psi_{+,0}^\mu(z)=\psi_{-,1}^0(y)+\psi_{+,1}^\mu(z)=x_2,$$which concludes the proof for $k_0=1$.

Suppose the proposition holds for $k_0=k$. We prove it also holds for $k_0=k+1$. The proof is similar to the case for $k_0=1$. Again by Lemma \ref{lem: com diag for n,n+1,n+2}, we have $$\ker(\Psi_{-,2k+2}^{k+1}- \Psi_{+,2k+2}^{k+1})\supset\im(( \Psi_{+,k+1}^0, \Psi_{-,k+1}^0)).$$
Then we prove $$\ker( \Psi_{-,2k+2}^{k+1}- \Psi_{+,2k+2}^{k+1})\subset\im(( \Psi_{+,k+1}^0,\Psi_{-,k+1}^0)).$$Suppose $$(x_1,x_2)\in \ker( \Psi_{-,2k+2}^{k+1}- \Psi_{+,2k+2}^{k+1}),\text{\textit{i.e.}}, \Psi_{-,2k+2}^{k+1}(x_1)- \Psi_{+,2k+2}^{k+1}(x_2)=0.$$Then we have $$\psi_{+,\mu}^{k+1}(x_1)=\psi_{+,\mu}^{2k+2}\circ\Psi_{-,2k+2}^{k+1}(x_1)= \psi_{+,\mu}^{2k+2}\circ\Psi_{+,2k+2}^{k+1}(x_2)=0.$$By exactness, there exists $y_1\in \sut{k}$ such that $\psi_{+,k+1}^k(y_1)=x_1$. By a similar reason, there exists $y_2\in\sut{k}$ such that $\psi_{-,k+1}^k(y_2)=x_2$. The goal is to prove $$\Psi_{-,2k}^k(y_1^\p)= \Psi_{+,2k}^k(y_2^\p)$$for some modifications $y^\p_1$ and $y^\p_2$ of $y_1$ and $y_2$ as for $y^\p$ in the case of $k_0=1$. Then the induction hypothesis will imply that there exists $w\in\sut{0}$ such that $$\Psi_{+,k}^0(w)= y_1^\p\aand \Psi_{-,k}^0(w)= y_2^\p.$$Hence we will have$$\Psi_{+,k+1}^0(w)= \psi_{+,k+1}^k(y_1^\p)= x_1\aand \Psi_{-,k+1}^0(w)=\psi_{-,k+1}^k(y_2^\p)= x_2.$$
This will conclude the proof for $k_0=k+1$.

Now we start to construct $y_1^\p$. We have $$\begin{aligned} \psi_{+,2k+2}^{2k+1}(\Psi_{+,2k+1}^{k+1}(x_2)-\Psi_{-,2k+1}^k(y_1))&= \psi_{+,2k+2}^{2k+1}\circ \Psi_{+,2k+1}^{k+1}(x_2)-\psi_{+,2k+2}^{2k+1}\circ \Psi_{-,2k+1}^k(y_1)\\&= \Psi_{+,2k+2}^{k+1}(x_2)-\psi_{+,2k+2}^{2k+1}\circ\Psi_{-,2k+1}^k(y_1)\\&= \Psi_{-,2k+2}^{k+1}(x_2)- \Psi_{+,2k+2}^{k+1}(x_1)\\&=0.\end{aligned}$$By exactness, there exists $z_1\in\sut{\mu}$ such that $$\psi_{+,2k+1}^\mu(z_1)=\Psi_{+,2k+1}^{k+1}(x_2)-\Psi_{-,2k+1}^k(y_1).$$
Let $y_1^\p=y_1+\psi_{+,k}^\mu(z_1).$ Then $$\psi_{+,k+1}^k(y_1^\p)= \psi_{+,k+1}^k(y_1)=x_1$$and
$$\begin{aligned}\Psi_{-,2k+1}^k(y_1^\p)&=\Psi_{-,2k+1}^k(y_1)+\Psi_{-,2k+1}^k\circ\psi_{+,k}^\mu(z_1)\\&=\Psi_{-,2k+1}^k(y_1)+\psi_{+,2k+1}^\mu(z_1)\\&=\Psi_{+,2k+1}^{k+1}(x_2),\end{aligned}$$

Then we start to construct $y_2^\p$. We have $$\begin{aligned}\psi_{-,2k+1}^{2k}(\Psi_{-,2k}^{k}(y_1^\p)-\Psi_{+,2k}^k(y_2))&=\Psi_{-,2k+1}^{k}(y_1^\p)-\psi_{-,2k+1}^{2k}\circ\Psi_{+,2k}^k(y_2)\\&=\Psi_{-,2k+1}^{k}(y_1^\p)-\Psi_{-,2k+1}^{k+1}(x_2)\\&=0.\end{aligned}$$By exactness, there exists $z_2\in\sut{\mu}$ such that $$\psi_{-,2k}^\mu(z_2)=\Psi_{-,2k}^{k}(y_1^\p)-\Psi_{+,2k}^k(y_2).$$
Let $y_2^\p=y_2+\psi_{-,k}^\mu(z_2).$ Then $$\psi_{-,k+1}^k(y_2^\p)= \psi_{-,k+1}^k(y_2)=x_2$$and
$$\begin{aligned}\Psi_{+,2k}^k(y_2^\p)&=\Psi_{+,2k}^k(y_2)+\Psi_{+,2k}^k\circ\psi_{-,k}^\mu(z_2)\\&=\Psi_{+,2k}^k(y_2)+\psi_{-,2k}^\mu(z_2)\\&=\Psi_{-,2k}^{k}(y_1^\p),\end{aligned}$$

Then we have the following commutative diagrams
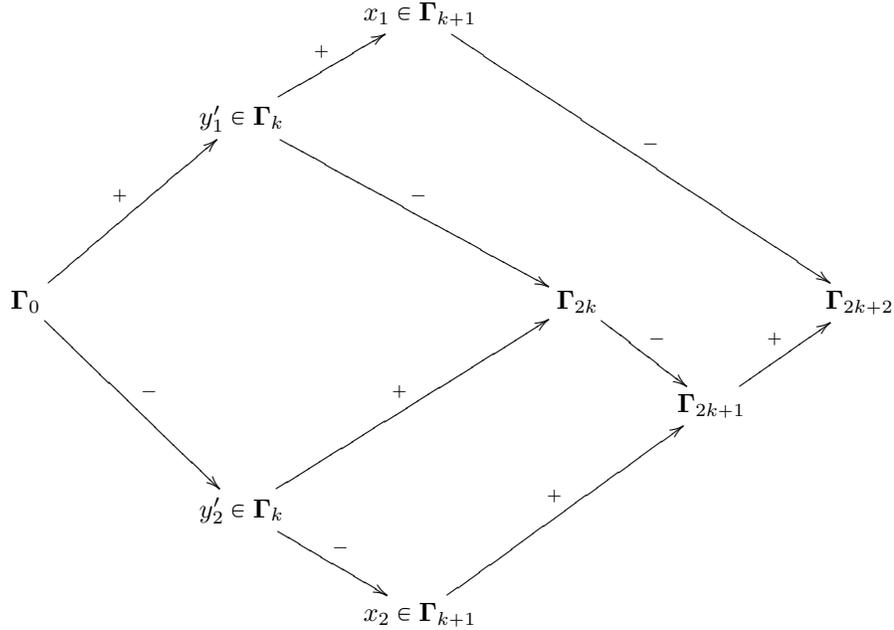
\begin{equation*}
\xymatrix{
&&&x_1\in\sut{k+1}\ar[dddrrr]^{-}&&&\\
&&y_1^\p\in\sut{k}\ar[ur]^{+}\ar[ddrr]^{-}&&&&\\
&&&&&&\\
\quad\sut{0}\quad\ar[rruu]^{+}\ar[ddrr]^{-}&&&&\sut{2k}\ar[dr]^{-}&&\sut{2k+2}\\
&&&&&\sut{2k+1}\ar[ur]^{+}&\\
&&y_2^\p\in\sut{k}\ar[dr]^{-}\ar[uurr]^{+}&&&&\\
&&&x_2\in\sut{k+1}\ar[uurr]^{+}&&&
}	
\end{equation*}By the induction hypothesis, there exists $w\in\sut{0}$ such that $$\Psi_{+,k}^0(w)= y_1^\p\aand \Psi_{-,k}^0(w)= y_2^\p,$$which concludes the proof for $k_0=k+1$.
\epf
\brem\label{rem: generalization}
By similar arguments, we can prove that the following sequence is exact for any $k_1,k_2\in\mathbb{N}_+$
$$\sut{n}\xra{(c_1\Psi_{+,n+k_1}^n,c_2\Psi_{-,n+k_2}^n)}\sut{n+k_1}\oplus\sut{n+k_2}\xra{c_3\Psi_{-,n+k_1+k_2}^{n+k_1}+c_4\Psi_{+,n+k_1+k_2}^{n+k_2}}\sut{n+k_1+k_2},$$where the scalars satisfies the equality $c_1c_3=-c_2c_4$.
\erem
\subsection{The second exact triangle}\label{subsec: second triangle}

In this subsection, we prove Proposition \ref{prop: exactness at other two, 2} by diagram chasing. For convenience, we restate it as follows, which is a little stronger than the previous version. Replacing the original knot in the proposition by the dual knot in the Dehn filling of slope $-(m+k)\mu+\lambda$ with framing $-\mu$ and setting $n=-1$ will recover Proposition \ref{prop: exactness at other two, 2}.

\bprop[]\label{prop2: exactness at other two, 2}
Suppose $$l^\p=\psi_{+,n-1}^{\mu}\circ\psi_{+,\mu}^{n+1}=\psi_{-,n-1}^{\mu}\circ\psi_{-,\mu}^{n+1}.$$Then for any $c_1,c_2,c_3,c_4\in\cstar$, the following sequence is exact
$$\sut{n}\oplus\sut{n}\xra{c_3\psi_{-,n+1}^n+c_4\psi_{+,n+1}^n}\sut{n+1}\xra{l^\p}\sut{n-1}\xra{(c_1\psi_{-,n}^{n-1},c_2\psi_{+,n}^{n-1})}\sut{n}\oplus\sut{n}.$$
\eprop
\bpf
We adopt the conventions from Section \ref{subsec Fixing the scalars}. We will use Lemma \ref{lem: bypass n,n+1,mu}, Lemma \ref{lem: com diag for n,n+1,n+2} and Lemma \ref{lem: comm diag for n,n+1,mu} without mentioning them. We prove the exactness at $\sut{n-1}$ first. We have$$\psi_{\pm,n}^{n-1}\circ l^\p= \psi_{\pm,n}^{n-1}\circ \psi_{\pm,n-1}^{\mu}\circ\psi_{\pm,\mu}^{n+1}=0.$$Hence $$\ker((c_1\psi_{-,n}^{n-1},c_2\psi_{+,n}^{n-1}))\supset \im(l^\p).$$Then we prove $$\ker((c_1\psi_{-,n}^{n-1},c_2\psi_{+,n}^{n-1}))\subset \im(l^\p).$$Suppose $$x\in \ker((c_1\psi_{-,n}^{n-1},c_2\psi_{+,n}^{n-1}))=\ker(\psi_{-,n}^{n-1})\cap \ker(\psi_{+,n}^{n-1}).$$By exactness, there exists $y\in\sut{\mu}$ such that $\psi_{+,n-1}^\mu(y)=x$. Then we have$$\psi_{+,n}^\mu(y)= \psi_{-,n}^{n-1}\circ\psi_{+,n-1}^\mu(y)=\psi_{-,n}^{n-1}(x)=0.$$By exactness, there exists $z\in \sut{n+1}$ such that $\psi_{+,\mu}^{n+1}(z)=y$. Thus, we have $l^\p(z)=x$, which concludes the proof for the exactness at $\sut{n-1}$.

Then we prove the exactness at $\sut{n+1}$. Similarly by exactness, we have $$\ker(l^\p)\supset \im(c_3\psi_{-,n+1}^n+c_4\psi_{+,n+1}^n)=\im(\psi_{-,n+1}^n)+\im(\psi_{+,n+1}^n).$$Suppose $x\in\ker(l^\p)$. If $\psi_{+,\mu}^{n+1}(x)=0$, then by the exactness, we know $x\in \im(\psi_{+,n+1}^n)$. If $\psi_{+,\mu}^{n+1}(x)\neq 0$, then by the exactness, there exists $y\in \sut{n}$ such that $$\psi_{+,\mu}^n(y)=\psi_{+,\mu}^{n+1}(x).$$Then we know$$x-\psi_{-,n+1}^n(y)\in\ker(\psi_{+,\mu}^{n+1})=\im(\psi_{+,n+1}^n).$$Thus, we have $$x\in \im(\psi_{-,n+1}^n)+\im(\psi_{+,n+1}^n),$$which concludes the proof for the exactness at $\sut{n+1}$.
\epf

\section{Some technical constructions}\label{sec: Some technical lemmas}
\subsection{Filtrations}
In this subsection, we study some filtrations on $\dehny{}$ and $\sut{\mu}$ that will be important in later sections. We continue to adopt conventions from Section \ref{subsec Fixing the scalars}. In particular, $K\subset Y$ is a rationally null-homologous knot and $S$ is a rational Seifert surface of $K$.

\blem
The maps $G_n$ in Lemma \ref{lem: surgery triangles} lead to a filtration on $\dehny{}$: for a sufficiently large integer $n_0$,
$$0=\ker G_{-n_0}\subset\dots\subset \ker G_{n}\subset \ker G_{n+1}\subset \dots\subset \ker G_{n_0}=\dehny{}.$$
\elem
\bpf
It follows from Lemma \ref{lem: F_n and G_n are iso when n large} that when $n_0$ is sufficiently large we have
$$0=\ker G_{-n_0}~{\rm and~}\ker G_{n_0}=\dehny{}.$$
It follows from Lemma \ref{lem: comm diag for n,n+1,dehn} that for any $n\in\intg$,
$$\ker G_{n}\subset \ker G_{n+1}.$$
\epf

\blem\label{lem: graded part of G_n}
For any $n\in\mathbb{Z}$, the map $G_n$ induces an isomorphism
$$G_{n}:\bigg(\ker G_{n+1}\slash \ker G_n\bigg)\xra{\cong} \ker \psp{n}{n+1}\cap\ker\psm{n}{n+1}.$$
\elem
\bpf
Suppose $x\in \ker G_{n+1}$. Then from Lemma \ref{lem: comm diag for n,n+1,dehn} we know that
$$\psi_{\pm,n+1}^n\circ G_{n}(x)=G_{n+1}(x)=0.$$
Hence we have 
$$G_n(\ker G_{n+1}))\subset \ker \psp{n}{n+1}\cap\ker\psm{n}{n+1}.$$
Clearly $G_n$ is injective on $\ker G_{n+1}\slash \ker G_n$ so it suffices to show that the image is $\ker \psp{n}{n+1}\cap\ker\psm{n}{n+1}$. To achieve this, for any element $x\in \ker \psp{n}{n+1}\cap\ker\psm{n}{n+1}$, Lemma \ref{lem: -1 surgery and bypass} implies that
$$x\in \ker H_{n}=\im G_{n}.$$
As a result, there exists $\alpha\in\dehny{}$ such that
$$x=G_{n}(\alpha).$$
Again from Lemma \ref{lem: comm diag for n,n+1,dehn} we know that
$$G_{n+1}(\al)=\psp{n+1}{n}\circ G_n(\al)=\psp{n+1}{n}(x)=0.$$
This implies that $\al\in\ker G_{n+1}.$
\epf

\blem
For any $n\in\mathbb{Z}$, the maps $\psi_{\pm,n}^{\mu}$ induce isomorphisms
$$\psp{\mu}{n}:\bigg(\im\psp{n+2}{\mu}\slash\im\psp{n+1}{\mu}\bigg)\xra{\cong} \ker \psp{n}{n+1}\cap\ker\psm{n}{n+1}$$
$$\psm{\mu}{n}:\bigg(\im\psm{n+2}{\mu}\slash\im\psm{n+1}{\mu}\bigg)\xra{\cong} \ker \psp{n}{n+1}\cap\ker\psm{n}{n+1}$$
\elem

\bpf
We only prove the lemma for positive bypasses. The proof for the negative bypasses is similar.
Let $u\in \im\psp{n+2}{\mu}$. By Lemma \ref{lem: bypass n,n+1,mu} and Lemma \ref{lem: comm diag for n,n+1,mu}, we have
$$\psp{n}{n+1}\circ\psp{\mu}{n}(u)=0\aand \psm{n}{n+1}\circ\psp{\mu}{n}(u)=\psp{\mu}{n+1}(u)=0.$$
Hence we know
$$\psp{\mu}{n}(\im\psp{n+2}{\mu})\subset \ker \psp{n}{n+1}\cap\ker\psm{n}{n+1}.$$
Since $\ker\psp{\mu}{n}=\im\psp{n+1}{\mu}$, the map $\psp{\mu}{n}$ is injective on $\im\psp{n+2}{\mu}\slash\im\psp{n+1}{\mu}$. To show it is surjective as well, pick $x\in \ker \psp{n}{n+1}\cap\ker\psm{n}{n+1}$. Note $x\in \ker\psp{n}{n+1}=\im\psp{\mu}{n}$ implies that there exists $u\in\sut{\mu}$ such that $\psp{\mu}{n}(u)=x$. Lemma \ref{lem: comm diag for n,n+1,mu} then implies that
$$\psp{\mu}{n+1}(u)=\psm{n}{n+1}\circ\psp{\mu}{n}(u)=\psm{n}{n+1}(x)=0.$$
As a result, $u\in \ker\psp{\mu}{n+1}=\im\psp{n+2}{\mu}$.
\epf

\bcor\label{cor: identifying filtrations on Y and Gamma_mu}
\begin{enumerate}
	\item For any $n\in\intg$, there is a canonical isomorphism
	$$\bigg(\ker G_{n+1}\slash\ker G_{n}\bigg)\cong\bigg(\im\psp{n+2}{\mu}\slash\im\psp{n+1}{\mu}\bigg)\cong\bigg(\im\psm{n+2}{\mu}\slash\im\psm{n+1}{\mu}\bigg).$$
	\item For sufficiently large $n_0$, there exists a (noncanonical) isomorphism
	$$\dehny{}\cong\bigg(\im\psp{n_0}{\mu}\slash\im\psp{-n_0}{\mu}\bigg)\cong\bigg(\im\psm{n_0}{\mu}\slash\im\psm{-n_0}{\mu}\bigg)$$
\end{enumerate}
\ecor


\bdefn\label{defn: F_n^i}
For any integer $n\in\intg$ and any grading $i$, define the map $F_n^i$ as the restriction
$$F_n^i=F_{n}|{\sutg{n}{i}}.$$
where $F_n$ is the map from Lemma \ref{lem: surgery triangles}.
\edefn



\blem\label{lem: ker of F_n^i}
Suppose $n_0\in\intg$ is small enough such that $F_{n_0}=0$ ({\it c.f.} Lemma \ref{lem: F_n and G_n are iso when n large}). Then for any integer $n\geq n_0$ and any grading $i$, we have
$$\psi^{n}_{\pm,\mu}(\ker F_n^i)=\im\bigg({\rm Proj}^{i\mp\frac{(n-1)p-q}{2}}_{\mu}\circ\psi^{n_0}_{\pm,\mu}\bigg),$$
where
$${\rm Proj}^{i\mp\frac{(n-1)p-q}{2}}_{\mu}:\sut{\mu}\to\sutg{\mu}{i\mp\frac{(n-1)p-q}{2}}$$
is the projection.
\elem

\bpf
We only prove the lemma for positive bypasses and the proof for negative bypasses is similar. First, suppose 
$$u\in \im\bigg({\rm Proj}^{i-\frac{(n-1)p-q}{2}}_{\mu}\circ\psi^{n_0}_{+,\mu}\bigg)=\im\psi^{n_0}_{+,\mu}\cap \sutg{\mu}{i-\frac{(n-1)p-q}{2}}.$$
Pick $x\in\sutg{n_0}{i-\frac{(n-n_0)p}{2}}$ such that
$$\psp{n_0}{\mu}(x)=u.$$
Taking $y=\Psm{n_0}{n}(x)$, we know from Lemma \ref{lem: bypass n,n+1,mu} that $y\in\sutg{n}{i},$ from Lemma \ref{lem: comm diag for n,n+1,dehn} that $F_n(y)=F_{n_0}(x)=0,$ and from Lemma \ref{lem: comm diag for n,n+1,mu} that $\psp{n}{\mu}(y)=u.$
As a result, we conclude $u\in \psp{n}{\mu}(\ker F_n^i)$.

Second, suppose $u\in\psp{n}{\mu}(\ker F_n^i)$ is non-zero. Pick $x_1\in\ker F_n^i$ such that
$$\psp{n}{\mu}(x_1)=u.$$
By Lemma \ref{lem: vanishing grading} and Lemma \ref{lem: bypass n,n+1,mu}, the fact that $\psp{n}{\mu}(x_1)=u\neq 0$ implies that
\begin{equation}\label{eq: lem 6.6, 1}
	-\frac{p-\chi(S)}{2}\le i-\frac{(n-1)p-q}{2}\le \frac{p-\chi(S)}{2}.
\end{equation}
Pick a sufficiently large integer $k$ and then take 
$$x_2=\Psm{n}{n+k}(x_1)\aand x_3=\Psp{n+k}{2n+2k-n_0}(x_2).$$
By Lemma \ref{lem: comm diag for n,n+1,dehn} we have $$F_{2n+2k-n_0}(x_3)=F_{n+k}(x_2)=F_{n}(x_1)=0.$$
Note that the grading $j$ of $x_3$ equals to 
\begin{equation}\label{eq: lem 6.6, 2}
	j=i+\frac{kp}{2}-\frac{(n+k-n_0)p}{2}=i-\frac{(n-n_0)p}{2}.
\end{equation}
Combining \ref{eq: lem 6.6, 1} and \ref{eq: lem 6.6, 2}, we obtain
\[\frac{(n_0-2)p-q-\chi(S)}{2}\leq j\leq \frac{n_0p-q-\chi(S)}{2}.\]
Note that we pick $k$ to be a sufficiently large integer. In particular, we can assume
\[
\frac{-(2n+2k-n_0)p+q+\chi(S)}{2}+1\leq \frac{(n_0-2)p-q-\chi(S)}{2}
\]
and
\[
\frac{n_0p-q-\chi(S)}{2}\leq \frac{(2n+2k-n_0)p-q+\chi(S)}{2}.
\]
Thus $j$ is in this range as well and then Lemma \ref{lem: F_n and G_n are iso when n large} implies that $F_{2n+2k-n_0}$ is injective on the grading $j$. Hence $x_3=0$. Then the following Lemma \ref{lem: having pre-image} applies to $(x,y)=(x_2,0)$ and there exists $x_4\in\sut{n_0}$ such that $$\Psm{n_0}{n+k}(x_4)=x_2.$$Thus by Lemma \ref{lem: comm diag for n,n+1,mu},
$$u=\psp{n}{\mu}(x_1)=\psp{n+k}{\mu}(x_2)=\psp{n_0}{\mu}(x_4)\in \im\bigg({\rm Proj}^{i-\frac{(n-1)p-q}{2}}_{\mu}\circ\psi^{n_0}_{+,\mu}\bigg).$$
\epf

\blem\label{lem: having pre-image}
Suppose $n\in\intg$ and $k_1,k_2\in\mathbb{N}_+$. Suppose $x\in\sut{n+k_1},y\in\sut{n+k_2}$ such that
$$\Psp{n+k_1}{n+k_1+k_2}(x)=\Psm{n+k_2}{n+k_1+k_2}(y)$$
Then there exists $z\in\sut{n}$ such that
$$\Psm{n}{n+k_1}(z)=x~{\rm and~}\Psp{n}{n+k_2}(z)=y.$$
\elem

\bpf
This is a restatement of Remark \ref{rem: generalization}. The proof is similar to that of Proposition \ref{prop2: exactness at direct summand}.
\epf
\subsection{Tau invariants in a general $3$-manifold}\label{subsec: tau inv}

\bdefn\label{defn: tau inv}
An element $\al\in \dehny{}$ is called {\bf homogeneous} if there exists an $n\in\intg$ and a grading $i$ such that $\al\in \im F_n^i$. Note that from \ref{lem: comm diag for n,n+1,dehn} and Corollary \ref{cor: psi^n_+,n+1 is an iso}, we know that
\[\al \in \im F_n^i\Rightarrow \al\in \im F_{n+1}^{i\pm\frac{p}{2}}.\]
For a homogeneous element $\al\in \dehny{}$, we pick a sufficiently large $n_0$ and define
$$\tau^+(\al)\deq\max_{i}\{i~|~\exists x\in\sutg{n_0}{i},~F_{n_0}(x)=\al\}-\frac{(n_0-1)p-q}{2}$$
$$\tau^-(\al)\deq\min_{i}\{i~|~\exists x\in\sutg{n_0}{i},~F_{n_0}(x)=\al\}+\frac{(n_0-1)p-q}{2}$$
$$\tau(\al)\deq 1+\frac{\tau^-(\al)-\tau^+(\al)+q}{p}=\frac{\min-\max}{p}+n_0.$$
\edefn
We will prove the independence of these $\tau$ invariants about $n_0$ later in Lemma \ref{lem: basic properties of tau}.
\brem
Here we fix the knot $K\subset Y$ and define the tau invariants for a homogeneous element $\al\in I^{\sharp}(Y)$. The reason why we go in this order is because (1) currently the definition of homogeneous elements depends on the choice of the knot and (2) in this paper we only focus on the Dehn surgeries of a fixed knot.
\erem
\brem
The normalization $\mp\frac{(n_0-1)p-q}{2}$ comes from the grading shifts of $\psi_{\pm,\mu}^{n_0}$ in Lemma \ref{lem: bypass n,n+1,mu}. When $K$ is a knot inside $Y=S^3$, we have that $\tau^\pm(\al)$ is equal to the tau invariant $\tau_I(K)$ defined in \cite{li2019tau}, where $\al$ is the unique generator of $I^\sharp(-S^3)\cong\mathbb{C}$ up to a scalar. Then $\tau(\al)=1-2\tau_I(K)$.
\erem

\blem\label{lem: homogeneous elements}
We have the following properties.
\begin{enumerate}
	\item Suppose $n_1,n_2$ are two integers and $i_1,i_2$ are two gradings such that there exist $x_1\in\sutg{n_1}{i_1}$ and $x_2\in\sutg{n_2}{i_2}$ with
	$$F_{n_1}(x_1)=F_{n_2}(x_2)\neq0.$$
	Then there exists an integer $N$ such that 
	$$i_2=i_1-\frac{(n_2-n_1)p}{2}+Np$$\textit{i.e.} when we send $x_1$ and $x_2$ into the same $\sut{n_3}$ with $n_3>n_1,n_2$ by bypass maps, then the difference of the expected gradings of the images is divisible by $p$ (the grading shifts of the bypass maps $\psi_{\pm,n+1}^n $ are $\mp p/2$).
	\item Suppose we have an integer $n_1$, a grading $i_1$, and an element $x_1\in\sutg{n_1}{i_1}$. Then for any integer $n_2\geq n_1$ and grading $i_2$ such that there exists an integer $N\in[0,n_2-n_1]$ with
	$$i_2=i_1-\frac{(n_2-n_1)p}{2}+Np,$$
	there exists an element $x_2\in\sutg{n_2}{i_2}$ such that
	$$F_{n_1}(x_1)=F_{n_2}(x_2).$$
	\item Suppose $n\in\intg$ and for $1\leq j\leq l$ we have a grading $i_j$ and an element $x_j\in\sutg{n}{i_j}$ such that $F_n(x_1)$,\dots, $F_n(x_{l})$ are linearly independent. Then the element
	$$\al=\sum_{j=1}^lF_n(x_j)$$
	is homogeneous if and only if for any $1\leq j\leq l$, we have
	$$i_j\equiv i_1~({\rm mod~}p)$$
\end{enumerate}
\elem
\bpf
(1). Take $n_0$ a sufficiently large integer. For $j=1,2$, take $i'_j\in(-\frac{p}{2},\frac{p}{2}]$ to be the unique grading such that there exists an integer $N_j$ with
$$i_j'=i_j-\frac{(n_0-n_j)p}{2}+N_jp.$$
Take
$$x_j'=\Psp{n_j+N_j}{n_0}\circ\Psm{n_j}{n_j+N_j}(x_j).$$
From Lemma \ref{lem: comm diag for n,n+1,dehn} we know that
$$x'_j\in \sutg{n_0}{i_j'}~{\rm and~}F_{n_0}(x_1')=F_{n_1}(x_1)=F_{n_2}(x_2)=F_{n_0}(x_2').$$
By Lemma \ref{lem: F_n and G_n are iso when n large}, we know that $x_1'=x_2'$ and in particular, $i_1'=i_2'$. As a result, we can take $N=N_1-N_2$ then it is straightforward to verify that
$$i_2=i_1-\frac{(n_2-n_1)p}{2}+Np.$$

(2). We can take
$$x_2=\Psm{n_1+N}{n_2}\circ\Psp{n_1}{n_1+N}(x_1)$$
Then it follows from Lemma \ref{lem: bypass n,n+1,mu} that $x_2\in\sutg{n_2}{i_2}$ and follows from Lemma \ref{lem: comm diag for n,n+1,dehn} that
$$F_{n_2}(x_2)=F_{n_1}(x_1).$$

(3). The proof is similar to that of (1).
\epf

\blem\label{lem: basic properties of tau}
For a homogeneous element $\al$, we have the following.
\begin{enumerate}
	\item $\tau^{\pm}(\al)$ and hence $\tau(\al)$ are well-defined. (\textit{i.e.} they are independent of the choice of the large integer $n_0$.)
	\item We have $\tau(\al)\in\intg$.
	\item For any integer $n$ and grading $i$, the following two statements are equivalent.
	\begin{enumerate}
		\item There exists $x\in\sutg{n}{i}$ such that $F_{n}(x)=\al$.
		\item We have $n\geq \tau(\al)$ and there exists $N\in\intg$ such that $N\in[0,n-\tau(\al)]$ and
	$$i=\frac{\tau^+(\al)+\tau^-(\al)-(n-\tau(\al))p}{2}+Np.$$
	\end{enumerate}
	\item We have
	$$\tau^+(\al)\geq -\frac{p-\chi(S)}{2}~{\rm and}~\tau^-(\al)\leq \frac{p-\chi(S)}{2}.$$
\end{enumerate}
\elem
\bpf
(1). Suppose $\al$ is a homogeneous element. Then by definition there exists $x\in\sutg{n}{i}$ for some integer $n$ and grading $i$ such that
$$F_n(x)=\al.$$
Then for sufficiently large $n_0$, we can take
$$y=\psp{n}{n_0}(x)$$
and from Lemma \ref{lem: comm diag for n,n+1,dehn} implies that
$$F_{n_0}(y)=\al$$
and hence $\tau^{\pm}(\al)$ exists.

 To show the value of $\tau^{\pm}(\al)$ is independent of $n_0$ as long as it is sufficiently large, a combination of Lemma \ref{lem: vanishing grading} and Lemma \ref{lem: bypass n,n+1,mu} implies that the map
$$\psm{n_0}{n_0+1}:\sutg{n_0}{i}\to\sutg{n_0+1}{i+\frac{p}{2}}$$
is an isomorphism for any $i>g-\frac{n_0p-q-1}{2}$. Then Lemma \ref{lem: comm diag for n,n+1,dehn} implies that $\tau^+$ is well-defined. The argument for $\tau^-$ is similar.

(2). It follows directly from Lemma \ref{lem: homogeneous elements} part (1).

(3). We first establish the following claim.

{\bf Claim}. There exists an element
\[z\in\sutg{\tau(\al)}{\frac{\tau^+(\al)+\tau^-(\al)}{2}}\]
such that
\[F_{\tau(\al)}(z)=\al.\]
\bpf[Proof of the claim]
Suppose $n_0\in\intg$ is sufficiently large and 
$$x_{\pm}\in\sutg{n_0}{\tau^{\pm}(\al)\pm\frac{(n_0-1)p-q}{2}}$$
such that $F_{n_0}(x_{\pm})=\al$. Note that the existence of $x_{\pm}$ follows from the definition of $\tau^{\pm}(\al)$. Let
$$x_{\pm}'=\Psi^{n}_{\pm,2n_0-\tau(\al)}(x_{\pm}).$$
It follows from Lemma \ref{lem: bypass n,n+1,mu} that
$$x_{\pm}'\in\sutg{2n_0-\tau(\al)}{\frac{\tau^+(\al)+\tau^-(\al)}{2}}.$$
From Lemma \ref{lem: comm diag for n,n+1,dehn} we know that
$$F_{2n_0-\tau(\al)}(x_+')=\al=F_{2n_0-\tau(\al)}(x_-').$$
By Lemma \ref{lem: F_n and G_n are iso when n large} this implies that
$$x_+'=x_-'.$$
Hence Lemma \ref{lem: having pre-image} applies and there exists $z\in\sut{\tau(\al)}$ such that
$$\Psi^{\tau(\al)}_{\mp,n}(z)=x_{\pm}.$$
Again Lemma \ref{lem: bypass n,n+1,mu} implies that $z$ is in the grading
$$z\in\sutg{\tau(\al)}{\frac{\tau^+(\al)+\tau^-(\al)}{2}}$$
and Lemma \ref{lem: comm diag for n,n+1,dehn} implies
$$F_{\tau(\al)}(z)=\al.$$
\epf

Now if an integer $n$ and a grading $i$ satisfy statement (b), then (a) is a direct consequence of the above claim and Lemma \ref{lem: homogeneous elements} part (2).

It remains to show that (a) implies (b). Suppose there exists $x\in\sutg{n}{i}$ such that $F_n(x)=\alpha$. From the above claim, we already know that there exists
$$z\in\sutg{\tau(\al)}{\frac{\tau^+(\al)+\tau^-(\al)}{2}}$$
such that
$$F_{\tau(\al)}(z)=\al$$
Hence Lemma \ref{lem: homogeneous elements} part (1) implies that there exists $N\in\intg$ such that
$$i=\frac{\tau^+(\al)+\tau^-(\al)-(n-\tau(\al))p}{2}+Np.$$
If $N>n-\tau(\al)$, we can take a sufficiently large $n_0$ and
$$x'=\Psm{n}{n_0}(x).$$
It follows from Lemma \ref{lem: bypass n,n+1,mu} that
$$x'\in\sutg{n_0}{i'}~{\rm with~}i'>\tau^+(\al)+\frac{(n_0-1)p-q}{2}.$$
Then Lemma \ref{lem: comm diag for n,n+1,dehn} implies that
$$F_{n_0}(x')=\al$$
which contradicts the definition of $\tau^+$ in Definition \ref{defn: tau inv}. Similarly if $N<0$ we can take
$$x'=\Psp{n}{n_0}(x)$$
which would be an element contradicting the definition of $\tau^-$. When $n<\tau(\al)$ we have $n-\tau(\al)<0$ so there is always a contradiction by the above argument. This concludes (b).

(4). It follows from the definition of $\tau^{\pm}$ and Lemma \ref{lem: F_n and G_n are iso when n large} that $F_{n_0}$ is an isomorphism when restricted to the direct sum of $p$ consecutive middle gradings of $\sut{n_0}$ when ${n_0}$ is large.
\epf

\blem\label{lem: im F_n is generated by homogeneous elements with tau<=n}
 For any $n\in\intg$ we have that 
	$$\im F_n={\rm Span}\{\al\in\dehny{}~|~\al~{\rm homogeneous~and}~\tau(\al)\leq n\}$$
\elem
\bpf
Suppose $\al\in\im F_n$. Let
$$\al=\sum_{i}\al_i~{\rm where~}\al_i\in \im F_n^i~{\rm is~homogeneous}.$$
From Lemma \ref{lem: basic properties of tau} we know that $\tau(\al_i)\leq n$ for all $i$.
On the other hand, suppose
$$\al=\sum_i\al_i~{\rm where}~\tau(\al_i)\leq n~{\rm for~all~}i.$$
By Lemma \ref{lem: basic properties of tau} part (3) we can pick $z_i\in\sut{\tau(\al_i)}$ such that
$$F_{\tau(\al_i)}(z_i)=\al_i.$$
Then from Lemma \ref{lem: comm diag for n,n+1,dehn} we know
$$\al=F_n(\sum_{i}\Psp{\tau(\al_i)}{n}(z_i)).$$

\epf

\subsection{A basis for framed instanton homology}\label{subsec: basis}
We pick a basis $\mathfrak{B}$ for $\dehny{}$ as follows. First
$$\mathfrak{B}=\mathop{\bigcup}_{n\in\intg}\mathfrak{B}_n.$$
To construct the set $\mathfrak{B}_n$, first, let $\mathfrak{B}_n=\emptyset$ if $F_n=0$. By Lemma \ref{lem: F_n and G_n are iso when n large} this means $\mathfrak{B}_n=\emptyset$ for all small enough $n$. Write
$$\mathfrak{B}_{\leq n}=\mathop{\bigcup}_{k\leq n}\mathfrak{B}_{k}.$$
We pick the set $\mathfrak{B}_n$ inductively. Note that we have taken $\mathfrak{B}_n=\emptyset$ for $n$ with $F_n=0$. Suppose we have already constructed the set $\mathfrak{B}_{\leq n-1}$ that consists of homogeneous elements and is a basis of $\im F_{n-1}$, we pick the set $\mathfrak{B}_{n}$ such that $\mathfrak{B}_{n}$ consists of homogeneous elements with $\tau=n$, and the set
$$\mathfrak{B}_{\leq n}=\mathfrak{B}_{\leq n-1}\cup \mathfrak{B}_{n}$$
forms a basis of $\im F_n$. Note that Lemma \ref{lem: im F_n is generated by homogeneous elements with tau<=n} implies that $\mathfrak{B}_{n}$ exists and
$$|\mathfrak{B}_{n}|=\dim_\mathbb{C}\bigg(\im F_n\slash \im F_{n-1}\bigg).$$

For any $n,k\in\intg$ such that $k\leq n-2$, define maps
$$\eta^n_{\pm,k}:\mathfrak{B}_n\to \sut{k}$$
as follows: for any $\al\in \mathfrak{B}_n\subset\im F_n$, since $\al$ is homogeneous and $\tau(\al)=n$, we can pick $$z\in\sutg{n}{\frac{\tau^+(\al)+\tau^-(\al)}{2}}$$by Lemma \ref{lem: basic properties of tau} part (3) such that $F_n(z)=\al$. Then define 
$$\eta^n_{\pm,k}(\al)=\psi^{\mu}_{\pm,k}\circ\psi^{n}_{\pm,\mu}(z)$$.

\blem\label{lem: the map eta^n_pm,k}
Suppose $n,k\in \intg$ such that $k\leq n-2$. 
\begin{enumerate}
	\item The maps $\eta^n_{\pm,k}$ are all well-defined.
	\item We have 
	$\eta^{n}_{+,n-2}=c_n\cdot \eta^n_{-,n-2}.$ for some scalar $c_n\in\cstar$.
	\item Elements in $\im\eta^n_{\pm,k}\subset\sut{k}$ are linearly independent.
	\item $\im\eta^n_{\pm,n-2}$ forms a basis for $\ker\psp{n-2}{n-1}\cap\ker\psm{n-2}{n-1}$.
	\item For any $\al\in\mathfrak{B}_n$ we have
	$$\eta^n_{\pm,k}(\al)\in\sutg{k}{\frac{\tau^+(\al)+\tau^-(\al)}{2}\mp\frac{(n-2-k)p}{2}}.$$
	\item We have
	$$\psi^{k-1}_{\mp,k}\circ\eta^n_{\pm,k-1}=\eta^{n}_{\pm,k},~{\rm and~}\psi^{k-1}_{\pm,k}\circ\eta^n_{\pm,k-1}=0.$$
\end{enumerate}
\elem

\bpf
(1). We only work with $\eta^n_{+,k}$ and the arguments for $\eta^n_{-,k}$ are similar. Suppose there are $z_1,z_2\in \sutg{n}{i}$ such that $F_n(z_1)=F_n(z_2)=\alpha$, where $i=\frac{\tau^+(\al)+\tau^-(\al)}{2}$. Then
$$z_1-z_2\in\ker F_n^i$$
and by Lemma \ref{lem: ker of F_n^i} we have
$$\psp{n}{\mu}(z_1-z_2)\in\psp{n}{\mu}(\ker F_n^i)\subset \im\psp{n_0}{\mu}\subset\im\psp{k+1}{\mu}.$$
Here $n_0\in\intg$ is a small enough integer. As a result,
$$\eta^n_{+,k}(\al)=\psp{\mu}{k}\circ\psp{n}{\mu}(z_1)=\psp{\mu}{k}\circ\psp{n}{\mu}(z_2)$$
is well-defined.

(2). This follows directly from Lemma \ref{lem: com diag for n,n+1,n+2}. Note that in Section \ref{subsec Fixing the scalars} we do not fix the scalars of the second commutative diagram of Lemma \ref{lem: com diag for n,n+1,n+2}, and hence a non-zero coefficient $c_n$ would possibly arise.

(3). We only work with $\eta^n_{+,k}$ and the arguments for $\eta^n_{-,k}$ are similar. Suppose 
$$\mathfrak{B}_n=\{\al_1,\dots,\al_{l}\},~{\rm where~}l=|\mathfrak{B}_{n}|=\dim_\mathbb{C}\bigg(\im F_n\slash \im F_{n-1}\bigg).$$
Suppose there exists $\lambda_1,\dots,\lambda_{l}$ such that
$$\sum_{j=1}^l\lambda_j\cdot\eta^n_{+,k}(\al_j)=0.$$
Pick $z_j\in\sutg{n}{\frac{\tau^+(\al_j)+\tau^-(\al_j)}{2}}$ such that $F_n(z_j)=\al_j$. Then we have
$$\psp{\mu}{k}\circ \psp{n}{\mu}\bigg(\sum_{j=1}^l\lambda_jz_j\bigg)=0.$$
As a result, there exists $x\in\sut{k+1}$ such that
$$\psp{k+1}{\mu}(x)=\psp{n}{\mu}\bigg(\sum_{j=1}^l\lambda_jz_j\bigg).$$
Note that, from Lemma \ref{lem: comm diag for n,n+1,mu}, we know
$$\psp{n}{\mu}\circ\Psm{k+1}{n}(x)=\psp{k+1}{\mu}(x)=\psp{n}{\mu}\bigg(\sum_{j=1}^l\lambda_jz_j\bigg)$$
so as a result there exists $y\in\sut{n-1}$ such that
$$\sum_{j=1}^l\lambda_jz_j=\Psm{k+1}{n}(x)+\psp{n-1}{n}(y).$$
Hence by Lemma \ref{lem: comm diag for n,n+1,dehn} we have
\beq
\sum_{j=1}^l\lambda_j\al_j&=F_n(\sum_{j=1}^l\lambda_jz_j)\\
&=F_n\circ\Psm{k+1}{n}(x)+F_n\circ\psp{n-1}{n}(y)\\
&=F_{k+1}(x)+F_{n-1}(y)\\
&\subset\im F_{n-1}.
\eeq
Since $\al_j$ form a basis of $\mathfrak{B}_n$, the sum cannot be in $\im F_{n-1}$ except  $\lambda_i=0$ for all $i$.

(4). For $\al\in\mathfrak{B}_n$, pick $z\in\sutg{n}{\frac{\tau^+(\al)+\tau^-(\al)}{2}}$ such that $F_n(z)=\al$. Then by definition
$$\eta^n_{+,n-2}(\al)=\psp{\mu}{n-2}\circ\psp{n}{\mu}(z).$$
Now we can compute
$$\psp{n-2}{n-1}\circ\eta^n_{+,n-2}(\al)=\psp{n-2}{n-1}\circ\psp{\mu}{n-2}\circ\psp{n}{\mu}(z)=0,$$
and by Lemma \ref{lem: comm diag for n,n+1,mu}
$$\psm{n-2}{n-1}\circ\eta^n_{+,n-2}(\al)=\psm{n-2}{n-1}\circ\psp{\mu}{n-2}\circ\psp{n}{\mu}(z)=\psp{\mu}{n-1}\circ\psp{n}{\mu}(z)=0.$$
Hence $$\eta^n_{+,n-2}(\al)\in \ker\psp{n-2}{n-1}\cap\ker\psm{n-2}{n-1}.$$Then (4) follows from (3), Lemma \ref{lem: graded part of G_n}, and $\im F_{n}=\ke G_{n-1}$.

(5). It follows directly from the construction of $\eta^{n}_{\pm,k}$ and Lemma \ref{lem: bypass n,n+1,mu}.

(6). It follows from the construction of $\eta^{n}_{\pm,k}$, the commutativity in Lemma \ref{lem: comm diag for n,n+1,mu} and the exactness in Lemma \ref{lem: bypass n,n+1,mu}.
\epf
\begin{conv}
We can define $$\tilde{\eta}_{+,k}^n=\eta_{+,k}^n\aand \tilde{\eta}_{-,k}^n=c_n\cdot \eta_{-,k}^n$$such that $$\tilde{\eta}_{+,n-2}^n=\tilde{\eta}_{-,n-2}^n$$ and the new maps satisfy all properties in Lemma \ref{lem: the map eta^n_pm,k} except (2). We will use $\eta_{+,k}^n$ to denote $\tilde{\eta}_{+,k}^n$ in latter sections.
\end{conv}

\section{The map in the third exact triangle}\label{sec: third triangle}

In this section, we construct the map $l$ in Proposition \ref{prop: exactness at other two} and Proposition \ref{prop: commute 2} and show it satisfies the exactness and the commutative diagram. We continue to adopt conventions from Section \ref{subsec Fixing the scalars}. We restate the propositions as follows and no longer use the notations $l,l^\p$ for maps.

\bprop\label{prop: exact triangle for n,n+k,n+2k}
Suppose $n\in\mathbb{Z}$ is fixed and $k\in\mathbb{Z}$ is sufficiently large. Then there is an exact triangle
\begin{equation*}
	\xymatrix{
	\sut{n}\ar[rr]^{\Phi^{n}_{n+k}}&&\sut{n+k}\oplus\sut{n+k}\ar[dl]^{\Phi^{n+k}_{n+2k}}\\
	&\sut{n+2k}\ar[ul]^{\Phi^{n+2k}_n}.&
	}
\end{equation*}
where two of the maps are already constructed
$$\Phi^n_{n+k}\deq (\Psp{n}{n+k},\Psm{n}{n+k}):\sut{n}\to\sut{n+k}\oplus\sut{n+k}$$
$$\Phi^{n+k}_{n+2k}\deq \Psm{n+k}{n+2k}-\Psp{n+k}{n+2k} :\sut{n+k}\oplus\sut{n+k}\to \sut{n+2k}.$$
\eprop
\bprop\label{prop: commutative diagram involving Phi^n+2k_n}
Suppose $n\in\mathbb{Z}$ is fixed and $k\in\mathbb{Z}$ is sufficiently large. Suppose $\Phi_{n-1}^{n+2k-1}$ is constructed in Proposition \ref{prop: exact triangle for n,n+k,n+2k}. Then, there are two commutative diagrams up to scalars.
\begin{center}
\begin{minipage}{0.4\textwidth}
	\begin{equation*}
\xymatrix{
\sut{\frac{2n+2k+1}{2}}\ar[dd]^{\Psp{n+k+1}{n+2k}\circ \psm{\frac{2n+2k+1}{2}}{n+k+1}}\ar[rrr]^{\psp{n+k+1}{\mu}\circ \psm{\frac{2n+2k+1}{2}}{n+k+1}}&&&\sut{\mu}\ar[dd]^{\psp{\mu}{n}}\\
&&&\\
\sut{n+2k}\ar[rrr]^{\Phi^{n+2k}_{n}}&&&\sut{n}
}	
\end{equation*}
\end{minipage}
\begin{minipage}{0.4\textwidth}
		\begin{equation*}
\xymatrix{
\sut{\frac{2n+2k+1}{2}}\ar[dd]^{\Psm{n+k+1}{n+2k+1}\circ\psp{\frac{2n+2k+1}{2}}{n+k+1}}\ar[rrr]^{\psm{n+k+1}{\mu}\circ\psp{\frac{2n+2k+1}{2}}{n+k+1}}&&&\sut{\mu}\ar[dd]^{\psm{\mu}{n}}\\
&&&\\
\sut{n+2k}\ar[rrr]^{\Phi^{n+2k}_{n}}&&&\sut{n}
}	
\end{equation*}\end{minipage}
\end{center}
\eprop

\subsection{Characterizations of the kernel and the image}
Before constructing $\Phi^{n+2k}_{n}$, we characterize the spaces $\ker \Phi ^n_{n+k}$ and $\im \Phi^{n+k}_{n+2k}$. These results will motivate the construction of $\Phi^{n+2k}_{n}$ to ensure that 
$$\im \Phi^{n+2k}_{n}=\ker \Phi ^n_{n+k}\aand \ker \Phi^{n+2k}_{n}=\im \Phi^{n+k}_{n+2k}.$$

Since $\Phi^{n}_{n+k}$ and $\Phi^{n+k}_{n+2k}$ are constructed using bypass maps, it suffices to consider their restrictions on each grading.
\blem\label{lem: ker of Phi^n_n+k}
Suppose $n\in\intg$ is fixed and $k\in\intg$ is sufficiently large. Let $${\rm Proj}^i_n:\sut{n}\ra\sutg{n}{i}$$
be the projection. Then we have
$$\ker\Phi^{n}_{n+k}\cap\sutg{n}{i}=\im\bigg({\rm Proj}^i_n\circ G_n\bigg).$$
\elem
\bpf
We need to apply Lemma \ref{lem: -1 surgery and bypass}. Following conventions in Section \ref{subsec Fixing the scalars}, we have
\begin{equation}
	H_n=\psp{n}{n+1}-\psm{n}{n+1}.
\end{equation}
Suppose $x\in\im\bigg({\rm Proj}^i_n\circ G_n\bigg)$. Pick $\al\in\dehny{}$ and $y\in\sut{n}$ such that
\[G_n(\al)=x+y\]
where ${\rm Proj}^i_n(y)=0$.
When $k$ is sufficiently large, we know from Lemma \ref{lem: F_n and G_n are iso when n large} that
$$G_{n+k}\equiv0.$$
In particular, from Lemma \ref{lem: comm diag for n,n+1,dehn}
$$\Psi^{n}_{\pm,n+k}(x)+\Psi^{n}_{\pm,n+k}(y)=G_{n+k}(\al)=0.$$
Since the maps $\Psi^{n}_{\pm,n+k}$ are homogeneous, we know that
$$\Psi^{n}_{\pm,n+k}(x)=0,$$
which implies that $x\in\ker\Phi^{n}_{n+k}\cap\sutg{n}{i}$.

Next, suppose $x\in \ker\Phi^{n}_{n+k}\cap\sutg{n}{i}$. We take $x^i_n=x$ and we will pick $x^j_n\in\sutg{n}{j}$ for all $j\neq i$ such that
$$\sum_{j}x^j_n\in \ker H_n=\im G_n.$$
We will use the notation $x_a^b$ to denote an element in $\sutg{a}{b}$. Recall that from Lemma \ref{lem: bypass n,n+1,mu}, the grading shifts of $\psi_{\pm,n+1}^n$ are $\mp \frac{p}{2}$. Take 
$$x^{i+\frac{(k-1)p}{2}}_{n+k-1}=\Psm{n}{n+k-1}(x)~{\rm and~}x^{i+\frac{(k+1)p}{2}}_{n+k-1}=0.$$
Since $x\in\ker\Phi^{n}_{n+k}\cap\sutg{n}{i}$ we know that
\begin{equation}\label{eq: lem 7.3, 1}
	\psm{n+k-1}{n+k}(x^{i+\frac{(k-1)p}{2}}_{n+k-1})=\Psm{n}{n+k}(x)=0=\psp{n+k-1}{n+k}(x^{i+\frac{(k+1)p}{2}}_{n+k-1}).
\end{equation}
Hence from Lemma \ref{lem: having pre-image}, there exists 
$$x_{n+k-2}^{i+\frac{kp}{2}}\in\sutg{n+k-2}{i+\frac{kp}{2}}$$
such that
$$\psm{n+k-2}{n+k-1}(x_{n+k-2}^{i+\frac{kp}{2}})=x^{i+\frac{(k+1)p}{2}}_{n+k-1}=0~{\rm and~}\psp{n+k-2}{n+k-1}(x_{n+k-2}^{i+\frac{kp}{2}})=x^{i+\frac{(k-1)p}{2}}_{n+k-1}.$$
Then we can take
$$x_{n+k-2}^{i+\frac{(k+2)p}{2}}=0~{\rm and}~x_{n+k-2}^{i+\frac{(k-2)p}{2}}=\Psm{n}{n+k-2}(x).$$
We can apply the same argument and use Lemma \ref{lem: having pre-image} to find
$$x_{n+k-3}^{i+\frac{(k+3)p}{2}},~x_{n+k-3}^{i+\frac{(k+1)p}{2}},~x_{n+k-3}^{i+\frac{(k-1)p}{2}},~x_{n+k-3}^{i+\frac{(k-3)p}{2}}\in\sut{n+k-3}$$
such that $\psi^{n+k-3}_{\pm,n+k-2}$ send them to corresponding elements in $\sut{n+k-2}$. Repeating this argument, we can obtain elements
\[x_n^{i+pj}\in\sutg{n}{i+pj}~{\rm for~}j\in[1,k]\cap\intg\]
such that $x_n^{i}=x$, $\psm{n}{n+1}(x_n^{i+pk})=0$, $\psp{n}{n+1}(x_n^{i-pk})=0$, and for any $j\in[1,k-1]\cap\intg$ we have
$$\psm{n}{n+1}(x_n^{i+pj})=\psp{n}{n+1}(x_n^{i+p(j+1)}).$$

Note that we obtain the above $x_n^{i+pj}$ for $j\in[1,k]\cap\intg$ essentially from the fact that $\Psi^n_{-,n+k}(x)=0$ as in Equation (\ref{eq: lem 7.3, 1}). However, $x\in {\rm ker}\Phi^n_{n+k}$ so we have $\Psi^n_{+,n+k}(x)=0$ as well. A similar argument as above then yields
\[x_n^{i+pj}\in\sutg{n}{i+pj}~{\rm for~}j\in[-k,-1]\cap\intg.\]
Together with $x_n^i=x$, we obtain $x_n^{i+pj}$ for all $j\in[-k,k]\cap\intg$.

It is then straightforward to check that
$$y=\sum_{j=-k}^k x_n^{i+pj}\in\ker(\psp{n}{n+1}-\psm{n}{n+1})=\ker H_n=\im G_{n}.$$
\epf
\blem\label{lem: image of Phi^n+k_n+2k}
Suppose $\al\in\dehny{}$ is a homogeneous element and
$$\al=\sum_{j=1}^l\lambda_j\cdot \al_j$$
where $\lambda_j\neq 0$ and $\al_j\in\mathfrak{B}$ for $1\leq j\leq l$. Let $n$ be an integer, $i$ be a grading and $k$ be a sufficiently large integer. For an element $x\in \sutg{n+2k}{i}$ such that $F_{n+2k}(x)=\alpha$,
the following is true.
\begin{enumerate}
	\item We have
	$$\tau^+(\al)=\min_{1\leq j\leq l}\{\tau^+(\al_j)\}~{\rm and~}\tau^-(\al)=\max_{1\leq j\leq l}\{\tau^-(\al_j)\}.$$
	\item We have $x\in\im\Phi^{n+k}_{n+2k}$ if and only if for any $1\leq j\leq l$, at least one of the following inequalities holds
	$$i\geq \tau^-(\al_j)-\frac{(n-1)p-q}{2}\aand i\leq \tau^+(\al_j)+\frac{(n-1)p-q}{2}.$$
	\item If $x\notin\im\Phi^{n+k}_{n+2k}$ then there exists $j,N\in\intg$ such that $1\leq j\leq l$, $0\leq N\leq \tau(\al_j)-n-2$, and 
	$$ i=\tau^+(\al_j)+\frac{(n-1)p-q}{2}+(N+1)p.$$
\end{enumerate}
\elem

\bpf
(1). We only demonstrate the proof of the result for $\tau^+$ and the proof for $\tau^-$ is similar. First, We make the following two claims. 

{\bf Claim 1}. For any homogeneous elements (not necessarily elements in $\mathfrak{B}$) $\al_1$ and $\al_2$ such that $\al_1+\al_2$ is also homogeneous, if
$\tau^+(\al_1)>\tau^+(\al_2)$ then $\tau^+(\al_1+\al_2)=\tau^+(\al_2).$

To prove Claim 1, let $n_0$ be sufficiently large. From Lemma \ref{lem: homogeneous elements} part (3) we know that
\begin{equation}\label{eq: lem 7.4, 1}
	\tau^+(\al_1)\equiv\tau^+(\al_2)\equiv\tau^+(\al_1+\al_2)~({\rm mod~}p).
\end{equation}
Assume $\tau^+(\al_1+\al_2)>\tau^+(\al_2)$. Let 
$$\tau^+=\min\{\tau^+(\al_1),\tau^+(\al_1+\al_2)\}>\tau^+(\al_2).$$

We claim that there exist
$$x_1,x_3\in\sutg{n_0}{\tau^++\frac{(n_0-1)p-q}{2}}$$
such that
$$F_{n_0}(x_1)=\al_1~{\rm and}~F_{n_0}(x_3)=\al_1+\al_2.$$
We prove only the existence of $x_1$, and the argument for the existence of $x_3$ is similar. By Definition \ref{defn: tau inv}, we know that 
\[\tau^+(\al_1)+\frac{(n_0-1)p-q}{2}=\frac{\tau^+(\al_1)+\tau^-(\al_1)-(n_0-\tau(\al_1))p}{2}+(n_0-\tau(\al_1))p.\]
Taking
\[N=(n_0-\tau(\al_1))-\frac{1}{p}(\tau^+(\al_1)-\tau^+),\]
we know that
\begin{equation}\label{eq: lem 7.4, 2}
	\tau^++\frac{(n_0-1)p-q}{2}=\frac{\tau^+(\al_1)+\tau^-(\al_1)-(n_0-\tau(\al_1))p}{2}+Np.
\end{equation}
Equation (\ref{eq: lem 7.4, 1}) implies that $N\in\intg$. The definition of $\tau^+$ makes sure that $N\leq n_0-\tau(\al_1)$. The fact that $n_0$ is sufficiently large and Lemma \ref{lem: basic properties of tau} part (4) implies that $N\geq 0$. Hence Lemma \ref{lem: basic properties of tau} part (3) implies the existence of $x_1$ such that
\[x_1\in \sutg{n_0}{\tau^++\frac{(n_0-1)p-q}{2}}\text{ and }F_{n_0}(x_1)=\al_1.\]

Now the existence of $x_1$ and $x_3$ implies that
$$F_n(x_3-x_1)=\al_2$$
which contradicts the definition of $\tau^+(\al_2)$.

{\bf Claim 2}. Suppose $\al_1,\dots,\al_u\in\mathfrak{B}$ are pairwise distinct elements in $\mathfrak{B}$ such that
$$\tau^+(\al_1)=\tau^+(\al_2)=\dots=\tau^+(\al_u)=\tau^+.$$
Suppose
$$\al^\p=\sum_{i=1}^u\lambda_i\cdot \al_i$$
and suppose it is homogeneous. Then $\tau^+(\al^\p)=\tau^+$.

To prove Claim 2, assume that $\tau^+(\al')>\tau^+$. Without loss of generality, assume that $\lambda_1\neq0$ and
$$\tau^-(\al_1)=\min_{1\leq j\leq u}\{\tau^-(\al_j)\}.$$
Then a similar argument as in the proof of Claim 1 implies that
$$\tau^-(\al')\leq\tau^-(\al_1).$$
Note that we have assumed $\tau^+(\al')>\tau^+=\tau^+(\al_1)$. Hence by Definition \ref{defn: tau inv}, $\tau(\al')<\tau(\al_1)$, which contradicts the construction of the set $\mathfrak{B}$.

Now we prove part (1). Suppose $\al_1,\dots,\al_l\in\mathfrak{B}$ are pairwise distinct elements in $\mathfrak{B}$. Let 
$$\al=\sum_{j=1}^l\lambda_j\cdot \al_j.$$
We want to show that
$$\tau^+(\al)=\min_{1\leq j\leq l}\{\tau^+(\al_j)\}.$$
To do this, relabel the elements $\al_j$ if necessary such that
$$\tau^+(\al_1)=\tau^+(\al_2)=\dots=\tau^+(\al_u)<\tau^+(\al_{u+1})\leq\tau^+(\al_{u+2})\leq\dots\leq\tau^+(\al_l).$$
Since $\al$ is homogeneous, from Lemma \ref{lem: homogeneous elements} part (3), we know that the sum $$\sum_{j=1}^v\lambda_j\cdot \al_j$$is also homogeneous for any $v=1,\dots,l$.
Applying Claim 2, we conclude that
$$\tau^+\bigg(\sum_{j=1}^u\lambda_j\cdot\al_j\bigg)=\tau^+(\al_1).$$
Hence we can apply Claim 1 repeatedly to conclude that
$$\tau^+\bigg(\sum_{j=1}^l\lambda_j\cdot\al_j\bigg)=\tau^+(\al_1)=\min_{1\leq j\leq l}\{\tau^+(\al_j)\}.$$

(2). If $x\in\im\Phi^{n+k}_{n+2k}$, then there exists $y\in\sutg{n+k}{i-\frac{kp}{2}}$ and $z\in\sutg{n+k}{i+\frac{kp}{2}}$ such that
$$x=\Psm{n+k}{n+2k}(y)+\Psp{n+k}{n+2k}(z).$$
By assumption $$F_{n+2k}(x)=\al=\sum_{j=1}^l\lambda_j\cdot \al_j$$with $\lambda_j\neq 0$ and $\al$ homogeneous. By Lemma \ref{lem: comm diag for n,n+1,dehn} we have
$$\al=F_{n+k}(y+z).$$
Since $\mathfrak{B}$ forms a basis for $\dehny{}$, we can write
$$F_{n+k}(y)=\sum_{j=1}^{l}\lambda_j^\p\cdot\al_j~{\rm and~}F_{n+k}(z)=\sum_{j=1}^l\lambda_j^\pp\cdot \al_j,$$where $l=|\mathfrak{B}|$. Then for any $1\le j\le l$, at least one of $\lambda_j^\p$ and $\lambda_j^\pp$ is non-zero. Since both $F_{n+k}(y)$ and $F_{n+k}(z)$ are homogeneous, from part (1) we know
$$i-\frac{kp}{2}\geq \tau^{-}(\al_j)-\frac{(n+k-1)p-q}{2}\text{ when }\lambda_j^\p\neq 0$$
$$\text{and }i+\frac{kp}{2}\leq \tau^{+}(\al_j)+\frac{(n+k-1)p-q}{2}\text{ when }\lambda_j^\pp\neq 0.$$

Conversely, suppose for any $1\leq j\leq l$ at least one of the following inequalities holds 
$$i\geq \tau^-(\al_j)-\frac{(n-1)p-q}{2}\aand i\leq \tau^+(\al_j)+\frac{(n-1)p-q}{2}.$$
We need to show $x\in \im\Phi^{n+k}_{n+2k}$. We deal with three cases.

{\bf Case 1}. The grading $i$ satisfies
\[
i \geq \frac{(n+2k-2)p-q+\chi(S)}{2}.
\]
We want to argue that
\[
\Psi^{n+k}_{-,n+2k}:\sutg{n+k}{i-\frac{kp}{2}}\to \sutg{n+2k}{i}
\]
is surjective and hence conclude that $x\in \im\Phi^{n+k}_{n+2k}$. To do this, note that $\Psi^{n+k}_{-,n+2k}|_{\sutg{n+k}{i-\frac{kp}{2}}}$ is the composition of maps $\psi^{n+k+j}_{-,n+k+j+1}|_{\sutg{n+k+j}{i-\frac{kp}{2}+\frac{jp}{2}}}$ for $j=0,1,...,k-1$. With the assumption of Case 1, we have
\[
i-\frac{kp}{2}+\frac{jp}{2}\geq \frac{(n+k+j-2)p-q+\chi(S)}{2}>-\frac{(n+k+j)p-q+\chi(S)}{2}.
\]
(Note that since $k$ is sufficiently large this is a very loose inequality.) Then Corollary \ref{cor: psi^n_+,n+1 is an iso} part (2) applies and we conclude that $\psi^{n+k+j}_{-,n+k+j+1}|_{\sutg{n+k+j}{i-\frac{kp}{2}+\frac{jp}{2}}}$ is an isomorphism for all $j=0,1,...,k-1$. Hence we conclude Case 1.

{\bf Case 2}. The grading $i$ satisfies
\[
i\leq -\frac{(n+2k-2)p-q+\chi(S)}{2}.
\]
The argument is similar to that for Case 1, except for using $\Psi^{n+k}_{+,n+2k}$ instead of $\Psi^{n+k}_{-,n+2k}$.

{\bf Case 3}. If the grading $i$ satisfies
\[
|i|< \frac{(n+2k-2)p-q+\chi(S)}{2}
\]
Under the assumption of Case 3, Lemma \ref{lem: F_n and G_n are iso when n large} part (1) implies that $F_{n+2k}$ is injective when restricted to $\sutg{n+2k}{i}$. 

Now, for each $j=1,2,...,l$, if we have 
\[i\geq \tau^-(\al_j)-\frac{(n-1)p-q}{2},\]
we claim that there exists $y_j\in\sutg{n+k}{i-\frac{kp}{2}}$ such that 
\[
F_{n+k}(y_j)=\lambda_j\cdot\al_j.
\]
If instead 
\[
i\leq \tau^+(\al_j)+\frac{(n-1)p-q}{2},
\]
we claim that there exists $z_j\in \sutg{n+k}{i+\frac{kp}{2}}$ such that
\[F_{n+k}(z_j)=\lambda_j\cdot\alpha_j\]
We will verify the existence of $y_j$ or $z_j$ in a moment, but for now let $y$ be the sum of all $y_j$'s and $z$ be the sum of all $z_j$'s. Then from Lemma \ref{lem: comm diag for n,n+1,dehn} it is straightforward to check that
$$F_{n+2k}(\Psm{n+k}{n+2k}(y)+\Psp{n+k}{n+2k}(z))=\alpha=F_{n+2k}(x).$$
Since in Case 3 the restriction of $F_{n+2k}$ on $\sutg{n+2k}{i}$ is injective, we conclude that
$$x=\Psm{n+k}{n+2k}(y)+\Psp{n+k}{n+2k}(z)\in\im\Phi^{n+k}_{n+2k}.$$

It remains to show that the desired $y_j$ or $z_j$ exists. We only prove the existence of $y_j$ and the argument for $z_j$ is similar. Now assume that
\[
i\geq \tau^-(\al_j)-\frac{(n-1)p-q}{2}.
\]
This implies that
\[
i-\frac{kp}{2}\geq \tau^-(\al_j)-\frac{(n+k-1)p-q}{2}.
\]
The hypothesis of the Lemma and the definition of $\tau^-(\al_j)$ in Definition \ref{defn: tau inv}, together with Lemma \ref{lem: homogeneous elements} part (3) imply that 
\[i\equiv \tau^-(\al_j)-\frac{(n-1)p-q}{2}~(\text{mod }p).\]
As a result, there exists an integer $N\geq 0$ such that
\[
i-\frac{kp}{2}=\tau^-(\al_j)-\frac{(n+k-1)p-q}{2}+Np.
\]
Note that, from Definition \ref{defn: tau inv}, we know
\[
\tau^-(\al_j)-\frac{(n+k-1)p-q}{2}=\frac{\tau^+(\al_j)+\tau^-(\al_j)-(n+k-\tau(\al_j))p}{2}.
\]
As a result, we have
\[
i-\frac{kp}{2}=\frac{\tau^+(\al_j)+\tau^-(\al_j)-(n+k-\tau(\al_j))p}{2}+Np.
\]
The assumption in Case 3 and Lemma \ref{lem: basic properties of tau} part (4) then implies that $N\leq (n+k)-\tau(\al_j)$. Hence Lemma \ref{lem: basic properties of tau} part (3) implies the existence of $y_j$.


(3). If $x\notin\im\Phi^{n+k}_{n+2k}$, then part (2) means that there exists some $j$ such that
$$\tau^+(\al_j)+\frac{(n-1)p-q}{2}<i<\tau^-(\al_j)-\frac{(n-1)p-q}{2}.$$
Note that, by Lemma \ref{lem: homogeneous elements} part (3), we must have
$$i\equiv \tau^-(\al_j)-\frac{(n-1)p-q}{2}\equiv \tau^+(\al_j)+\frac{(n-1)p-q}{2}~({\rm mod}~p).$$By direct calculation, we have $$(\tau^-(\al_j)-\frac{(n-1)p-q}{2})-(\tau^+(\al_j)+\frac{(n-1)p-q}{2})=(\tau(\al_j)-n)p$$Then we can choose $N$ with $0\le N\le \tau(\al_j)-n-2$ as desired.
\epf

\subsection{The construction of the map}\label{subsec: constructing Phi^n+2k_n}

Since $\Phi^{n}_{n+k}$ and $\Phi^{n+k}_{n+2k}$ are homogeneous, we can construct $\Phi^{n+2k}_{n}$ for each grading to achieve both the exactness and the commutativity. Given the grading shifts in Lemma \ref{lem: bypass n,n+1,mu} and Lemma \ref{lem: bypass n+1,n,2n+1/2}, the map $\Phi^{n+2k}_{n}$ preserves the gradings. From Lemma \ref{lem: vanishing grading}, for any grading $i$ with 
\[|i|>\frac{|np-q|-\chi(S)}{2},\]
we have $\sutg{n}{i}=0$. From Corollary \ref{cor: psi^n_+,n+1 is an iso}, we know either $\Psi_{+,n+2k}^{n+k}$ or $\Psi_{-,n+2k}^{n+k}$ is surjective onto $\sutg{n+2k}{i}$ for such grading $i$. Thus, on such grading $i$, the zero map satisfies the exactness for $\Phi^{n+2k}_{n}$ (though we still have to verify the commutativity in Proposition \ref{prop: commutative diagram involving Phi^n+2k_n}). 

On the other hand, from Lemma \ref{lem: F_n and G_n are iso when n large}, the restriction of  $F_{n+2k}$ on the consecutive $p$ middle gradings is an isomorphism. In particular, when $p=1$, it is an isomorphism when restricted to each middle grading. Also from Lemma \ref{lem: ker of Phi^n_n+k}, it seems that the definition of $\Phi^{n+2k}_{n+k}$ on $\sutg{n+2k}{i}$ should involve ${\rm Proj}^i_n\circ G_n$. However, if we simply take $${\rm Proj}^i_n\circ G_n\circ F_{n+2k}$$as the definition, the current techniques fall short of demonstrating exactness and commutativity. 



We resolve this issue by introducing an isomorphism
$$I:\dehny{}\xra{\cong}\dehny{}$$
and define \begin{equation}\label{defn: the map Phi^n+2k_n}
    \Phi_{n}^{n+2k}(x)={\rm Proj}^i_n\circ G_n\circ I\circ F_{n+2k}(x)\text{ for }x\in \sutg{n+2k}{i}.
\end{equation}
The construction of $I$ is noncanonical but it helps us to prove the exactness and commutativity.

\brem
In the first arXiv version of this paper, we deal with the special case $Y=S^3$. In this case $\dehny{}\cong\mathbb{C}$ so up to a scalar we have $I=\operatorname{Id}$. In this special case indeed we could prove the exactness and commutativity without explicitly writing down the isomorphism $I$ as follows.
\erem

We first define the map $I$ on the basis $$\mathfrak{B}=\mathop{\cup}_{n\in\intg}\mathfrak{B}_n$$of $\dehny{}$ chosen in Section \ref{subsec: basis} that consists of homogeneous elements and then extend the map on the whole space linearly. We will show it is an isomorphism. 

Fix $n_0\in\intg$ small enough such that Corollary \ref{cor: psi^n_+,n+1 is an iso} and Lemma \ref{lem: F_n and G_n are iso when n large} apply. For any $\al\in\mathfrak{B}_n$, there exists a grading $i(\al)\in(-\frac{p}{2},\frac{p}{2}]$ such that there exists $N(\al)\in\intg$ with
$$i(\al)=\frac{\tau^+(\al)+\tau^-(\al)}{2}-\frac{(\tau(\al)-2-n_0)p}{2}+N(\al)p.$$
Note that, from the above equality, we know that
\[N(\alpha)=\frac{\tau(\al)-2-n_0}{2}+\frac{2i(\al)-\tau^+(\al)-\tau^-(\alpha)}{2p}.\]
Note that, except for $n_0$, the rest of the terms are bounded, so $N(\alpha)\geq 0$ when $n_0$ is small enough. Similarly,
\[N(\alpha)+n_0=\frac{\tau(\al)-2+n_0}{2}+\frac{2i(\al)-\tau^+(\al)-\tau^-(\alpha)}{2p}\]
so we have $N(\alpha)+n_0\leq \tau(\alpha)-2$ when $n_0$ is small enough. As a result, by Lemma \ref{lem: the map eta^n_pm,k} part (2) and (6) (and the convention after the lemma), we know that
for any $\al\in\mathfrak{B}_n$ $$\Psm{n_0+N(\al)}{\tau(\al)-2}\bigg(\eta^{\tau(\al)}_{+,n_0+N(\al)}(\al)\bigg)=\eta^{\tau(\al)}_{+,\tau(\al)-2}(\al)=\eta^{\tau(\al)}_{-,\tau(\al)-2}(\al)=\Psp{\tau(\al)-2-N(\al)}{\tau(\al)-2}\bigg(\eta^{\tau(\al)}_{-,\tau(\al)-2-N(\al)}(\al)\bigg).$$
Then by Lemma \ref{lem: having pre-image}, there exists $w\in\sutg{n_0}{i(\al)}$ such that
\begin{equation}\label{eq: choice of beta}
	\Psp{n_0}{n_0+N(\al)}(w)=\eta^{\tau(\al)}_{+,n_0+N(\al)}(\al)~{\rm and}~\Psm{n_0}{\tau(\al)-2-N(\al)}(w)=\eta^{\tau(\al)}_{-,\tau(\al)-2-N(\al)}(\al).
\end{equation}
Let
$${\rm Proj}:\sut{n_0}\to\bigoplus_{i\in(-\frac{p}{2},\frac{p}{2}]}\sutg{n_0}{i}.$$
From Lemma \ref{lem: F_n and G_n are iso when n large}, we know
$${\rm Proj}\circ G_{n_0}:\dehny{}\to \bigoplus_{i\in(-\frac{p}{2},\frac{p}{2}]}\sutg{n_0}{i}$$
is an isomorphism. Hence we define
$$I(\al)=({\rm Proj}\circ G_{n_0})^{-1}(w).$$

The following diagram might be helpful for understanding the construction of $I$. (We write $n=\tau(\al)$, $n_1=\tau(\al)-2-N(\al)$, and $n_2=n_0+N(\al)$.)
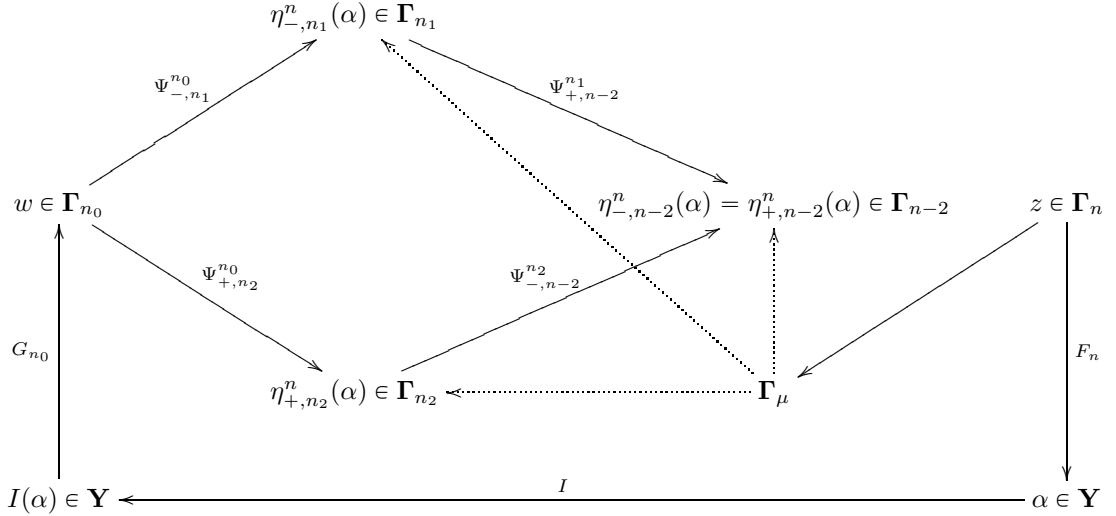
\begin{equation*}
	\xymatrix{
	&&\eta^{n}_{-,n_1}(\al)\in\sut{n_1}\ar[rrdd]^{\Psp{n_1}{n-2}}&&&\\
	&&&&&\\
	w\in\sut{n_0}\ar[uurr]^{\Psm{n_0}{n_1}}\ar[ddrr]^{\Psp{n_0}{n_2}}&&&&\eta^n_{-,n-2}(\al)=\eta^n_{+,n-2}(\al)\in\sut{n-2}&z\in\sut{n}\ar[ddd]^{F_n}\ar[ldd]\\
	&&&&&\\
	&&\eta^{n}_{+,n_2}(\al)\in\sut{n_2}\ar[uurr]^{\Psm{n_2}{n-2}}&&\sut{\mu}\ar@{.>}[luuluu]\ar@{.>}[uu]\ar@{.>}[ll]&\\
	I(\al)\in\dehny{}\ar[uuu]^{G_{n_0}}&&&&&\al\in\dehny{}\ar[lllll]_I\\
	}
\end{equation*}

\brem
For a general $3$-manifold $Y$, our construction of $I$ is noncanonical since there are many choices such as the basis $\mathfrak{B}$ and the element $w$ for each $\al\in\mathfrak{B}$. However, one could still ask whether we could simply pick $I=\operatorname{Id}$ or not. If we take $I=\operatorname{Id}$, then Proposition \ref{prop: commutative diagram involving Phi^n+2k_n} can finally be reduced to Conjecture \ref{conj: commutative diagram} which we state below. We believe that the following conjecture is true, though currently, we do not find a proof for it. Hence in order to fulfill the main purpose of the paper, we introduce the isomorphism $I$ to bypass this conjecture.
\erem

\begin{conj}\label{conj: commutative diagram}
	For any $\al\in\mathfrak{B}$, and any integer $n\leq \tau(\al)-2$, we have
$$\eta^{\tau(\al)}_{\pm,n}(\al)={\rm Proj}^{j_{\pm}}_n\circ G_{n}(\al),$$
where 
$$j_{\pm}=\frac{\tau^+(\al)+\tau^-(\al)}{2}\mp\frac{(\tau(\al)-2-n)p}{2}$$
and
$${\rm Proj}^{j_{\pm}}_n:\sut{n}\to\sutg{n}{j_{\pm}}$$
is the projection.
\end{conj}

\blem\label{lem: the isomorphism I}
We have the following.
\begin{enumerate}
	\item Suppose $\al\in\mathfrak{B}$ and $n_0$, $w$, $N(\al)$ are chosen as above. Suppose $n,k$ are two integers such that $n_0\leq k\leq n$. Then (a) $\Psm{k}{n}\circ\Psp{n_0}{k}(w)\neq 0$ if and only if (b) $k\leq n_0+N(\al)$ and $n-k\leq \tau(\al)-2-n_0-N(\al)$ (in particular, we have $n\leq \tau(\al)-2$).
	\item The map $I:\dehny{}\to\dehny{}$ is an isomorphism.
	\item For an element $\al\in \mathfrak{B}$, an integer $n$ and a grading $i$, the following two statements are equivalent.
	\begin{enumerate}
		\item We have ${\rm Proj}^i_n\circ G_n\circ I(\al)\neq 0.$
		\item We have $n\leq \tau(\al)-2$ and there exists $N\in\intg$ such that $N\in[0,\tau(\al)-2-n]$ and
	$$i=\frac{\tau^+(\al)+\tau^-(\al)-(\tau(\al)-2-n)p}{2}+Np.$$
	\end{enumerate}
	\item Suppose for an integer $n$ and a grading $i$ we have $\al_1,\dots,\al_L\in\mathfrak{B}$ such that ${\rm Proj}^i_n\circ G_n\circ I(\al_j)\neq 0$ for all $1\leq j\leq L$, then ${\rm Proj}^i_n\circ G_n\circ I(\al_1)$,\dots, ${\rm Proj}^i_n\circ G_n\circ I(\al_L)$ are linearly independent.
	\item Suppose $\al\in\mathfrak{B}$. For any $n\in\intg$ such that $n\leq\tau(\al)-2$, we have
	$${\rm Proj}^{i_{\pm}}_n\circ G_n\circ I(\al)=\eta^{\tau(\al)}_{\pm,n}(\al)~{\rm where~}i_{\pm}=\frac{\tau^{+}(\al)+\tau^{-}(\al)\mp(\tau(\al)-2-n)}{2}$$
\end{enumerate}
\elem
\bpf
(1). First, when $k>n_0+N(\al)$, from the construction of $w$, we know that
$$\Psp{n_0}{k}(w)=\Psp{n_0+N(\al)}{k}\circ\Psp{n_0}{n_0+N(\al)}(w)=\Psp{n_0+N(\al)}{k}\circ\eta^{\tau(\al)}_{+,n_0+N(\al)}(\alpha)=0.$$
The last equality is from Lemma \ref{lem: the map eta^n_pm,k} part (6). Similarly, if $n-k>\tau(\al)-2-n_0-N(\al)$, we know from Lemma \ref{lem: com diag for n,n+1,n+2} that
\beq
\Psm{k}{n}\circ\Psp{n_0}{k}(w)&=\Psp{n+n_0-k}{n}\circ\Psm{n_0}{n+n_0-k}(w)\\
&=\Psp{n+n_0-k}{n}\circ\Psm{\tau(\al)-2-N(\al)}{n+n_0-k}\circ\Psm{n_0}{\tau(\al)-2-N(\al)}(w)\\
({\rm Definition~of~}w)&=\Psp{n+n_0-k}{n}\circ\Psm{\tau(\al)-2-N(\al)}{n+n_0-k}\circ\eta^{\tau(\al)}_{-,\tau(\al)-2-N(\al)}(\al)\\
({\rm Lemma}~\ref{lem: the map eta^n_pm,k}~{\rm part~(6)})&=0.
\eeq

Next, we need to show that $\Psm{k}{n}\circ\Psp{n_0}{k}(w)\neq 0$ when $k\leq n_0+N(\al)$ and $n-k\leq \tau(\al)-2-n_0-N(\al)$. 
Again, from Lemma \ref{lem: com diag for n,n+1,n+2} we have
\beq
\Psp{n}{n+N(\al)+n_0-k}\circ\Psm{k}{n}\circ\Psp{n_0}{k}(w)&=\Psm{n_0+N(\al)}{n_0+N(\al)+n-k}\circ\Psp{n_0}{n_0+N(\al)}(w)\\
({\rm Definition~of~}w)&=\Psm{n_0+N(\al)}{n_0+N(\al)+n-k}\circ\eta^{\tau(\al)}_{+,n_0+N(\al)}(\al)\\
({\rm Lemma}~\ref{lem: the map eta^n_pm,k}~{\rm part~(6)})&=\eta^{\tau(\al)}_{+,n_0+N(\al)+n-k}(\al)\\
({\rm Lemma}~\ref{lem: the map eta^n_pm,k}~{\rm part~(3)})&\neq0.
\eeq

(2). Suppose $\mathfrak{B}=\{\al_1,\dots,\al_L\}$ where $L=\dim_\mathbb{C}\dehny{}$. We order the elements $\al_i$ such that
$$\tau(\al_i)\geq\tau(\al_{i+1}).$$
Let $w_i$, $N_i=N(\al_i)$ be the data associated to $\al_i$ as above. Since
$$w_i\in\bigoplus_{j\in(-\frac{p}{2},\frac{p}{2}]}\sutg{n_0}{j}$$
for any $i$, and by Lemma \ref{lem: F_n and G_n are iso when n large}, the map
$${\rm Proj}\circ G_{n_0}:\dehny{}\to \bigoplus_{j\in(-\frac{p}{2},\frac{p}{2}]}\sutg{n_0}{j}$$
is an isomorphism, in order to show that $I$ is an isomorphism, it suffices to show that $w_1$,\dots, $w_L$ are linearly independent.

Now suppose there are complex numbers $\lambda_1$,\dots, $\lambda_L$ such that
$$\sum_{i=1}^{L}\lambda_iw_i=0.$$
Our goal is to show that all $\lambda_i$ are zero. The idea is to apply various maps $\Psp{n_0+N_1}{\tau(\al_1)-2}\circ\Psm{n_0}{n_0+N_1}$ to filter out different indices by part (1) of the lemma.

Applying the map $\Psp{n_0+N_1}{\tau(\al_1)-2}\circ\Psm{n_0}{n_0+N_1}$, from the construction of $w_i$, the order of $\al_i$ and part (1) of the lemma, we know
\beq
0&=\Psp{n_0+N_1}{\tau(\al_1)-2}\circ\Psm{n_0}{n_0+N_1}(\sum_{i=1}^{L}\lambda_iw_i)\\
&=\sum_{\al_i}\lambda_i\eta^{\tau(\al_1)}_{\pm \tau(\al_1)-2}(\al_i)\\
\eeq
where the summation in the second line is over all $\al_i$ with
$$\tau(\al_i)=\tau(\al_1),~{\rm and}~N_i=N_1.$$
Note that, from Lemma \ref{lem: the map eta^n_pm,k} part (2) and the convention after the lemma, we know $\eta^{\tau(\al_1)}_{+ \tau(\al_1)-2}=\eta^{\tau(\al_1)}_{- \tau(\al_1)-2}$. From Lemma \ref{lem: the map eta^n_pm,k} again, we know that $\eta^{\tau(\al_1)}_{\pm \tau(\al_1)-2}(\al_i)$ are linearly independent, and as a result all relevant $\lambda_i$ must be zero.
Suppose $i_0$ is the smallest index in the rest. By our choice of $\al_i$, the element $\al_{i_0}$ has the largest $\tau$ among the rest of the $\al_i$. Hence we can apply the map $\Psp{n_0+N_{i_0}}{\tau(\al_{i_0})-2}\circ\Psm{n_0}{n_0+N_{i_0}}$ to filter out $\al_i$ with smaller $\tau$. Repeating this argument, we could prove that all $\lambda_i$ must be zero.

(3). Let $n_0$, $w$, $i(\al)$, and $N(\al)$ be constructed as above. We first prove that (b) $\Rightarrow$ (a). Note that, when constructing the isomorphism $I$, from Corollary \ref{cor: psi^n_+,n+1 is an iso} and Lemma \ref{lem: comm diag for n,n+1,dehn} we can take $n_0'=n_0-2$ and $w'=(\psp{n_0-2}{n_0-1})^{-1}\circ (\psp{n_0-1}{n_0})^{-1}(w)$ that will lead to the same $I$ as $n_0$, $w$. (Note that, by construction, passing from $n_0$ to $n_0-2$ will increase $N(\al)$ by $1$.) 

\brem
Note that the main goal of the current paper is to derive an integral surgery formula. For $n$ that is sufficiently large, we already know a large surgery as in \cite{LY2021large}. When $n$ is small enough, we can pass to $-n$ for the mirror of the knot. As a result, instead of changing $n_0$ for particular $n$, we could assume a universal bound for all the integers $n$ that we care about and make $n_0$ universally small.
\erem

As a result, we can always assume that $n_0$ is small enough compared with any given $n$. Now recall by construction
$$i(\al)=\frac{\tau^+(\al)+\tau^-(\al)}{2}-\frac{(\tau(\al)-2-n_0)p}{2}+N(\al)p$$
and by the assumption in (b) we have
$$i=\frac{\tau^+(\al)+\tau^-(\al)-(\tau(\al)-2-n)p}{2}+Np.$$
We can assume that $n_0$ is small enough such that $N(\al)>N$. Take $k=N(\al)-N+n_0$. Note that, by construction we have $w\in\sutg{n_0}{i(\al)}$, so from Lemma \ref{lem: bypass n,n+1,mu} we know that
\[\Psm{k}{n}\circ\Psp{n_0}{k}(w)\in \sutg{n}{i(\alpha)-\frac{(k-n_0)p}{2}+\frac{(n-k)p}{2}}=\sutg{n}{i}.\]
As a result, we conclude from the definition of $w$ and Lemma \ref{lem: comm diag for n,n+1,dehn} that
\[
\begin{aligned}
	\Psm{k}{n}\circ\Psp{n_0}{k}(w)&=\Psm{k}{n}\circ\Psp{n_0}{k}\circ {\rm Proj}_{n_o}^{i(\alpha)}\circ G_{n_0}\circ I(\alpha)\\
	&={\rm Proj}^i_n\circ \Psm{k}{n}\circ\Psp{n_0}{k}\circ G_{n_0}\circ I(\alpha)\\
	&={\rm Proj}^i_n \circ G_n\circ I(\al)
\end{aligned}
\]

Then it is straightforward to verify that $k-n_0\leq N(\alpha)$ and $n-k\leq \tau(\al)-2-n_0-N(\al)$. As a result, we conclude from part (1) that
$${\rm Proj}^i_n\circ G_n\circ I(\al)\neq 0.$$

Next we show that (a) $\Rightarrow$ (b). Again assume that $n_0$ is small enough compared with the given $n$. Then there exists $i'\in(-\frac{p}{2},\frac{p}{2}]$ such that there exists $N'\in\intg$ with
$$i'=i-\frac{(n-n_0)p}{2}+N'p.$$
By Lemma \ref{lem: comm diag for n,n+1,dehn}, we know that
$${\rm Proj}^{i}_n\circ G_n\circ I(\al)=\Psm{n_0+N'}{n}\circ\Psp{n_0}{n_0+N'}\circ {\rm Proj}^{i'}_{n_0}\circ G_{n_0}\circ I(\al).$$

From the construction of $I(\al)$ and Lemma \ref{lem: F_n and G_n are iso when n large} we know ${\rm Proj}^{i}_n\circ G_n\circ I(\al)\neq 0$ only if $i'=i(\al)$, in which case
$${\rm Proj}^{i}_n\circ G_n\circ I(\al)=\Psm{n_0+N'}{n}\circ\Psp{n_0}{n_0+N'}\circ {\rm Proj}^{i'}_{n_0}\circ G_{n_0}\circ I(\al)=\Psm{n_0+N'}{n}\circ\Psp{n_0}{n_0+N'}(w).$$
Hence ${\rm Proj}^{i}_n\circ G_n\circ I(\al)\neq 0$ implies that
$$N'\leq N(\al)~{\rm and~}n-N'\leq \tau(\al)-2-N(\alpha)$$
by part (1). Taking $N=N(\al)-N'$, it is then straightforward to check that
$$N\in[0,\tau(\al)-2-n]~{\rm and~}i=\frac{\tau^+(\al)+\tau^-(\al)-(\tau(\al)-2-n)p}{2}+Np.$$

(4). The proof is similar to that of (2).

(5). It follows from the proofs of parts (1) and (3).
\epf


\subsection{The exact triangle}
In this subsection, we prove the exact triangle. Note that we choose the basis  $\mathfrak{B}$ of $\dehny{}$ as in Section \ref{subsec: basis}.
\bpf[Proof of Proposition \ref{prop: exact triangle for n,n+k,n+2k}] We will verify the exactness at each space of the triangle.

{\bf The exactness at $\sut{n+k}\oplus\sut{n+k}$}. This follows from Proposition \ref{prop2: exactness at direct summand}.

{\bf The exactness at $\sut{n}$}. From Lemma \ref{lem: ker of Phi^n_n+k} and the construction of $\Phi^{n+2k}_n$ in (\ref{defn: the map Phi^n+2k_n}), we know that $\im\Phi^{n+2k}_n\subset \ker \Phi^{n}_{n+k}$. Now pick an arbitrary 
$$x\in\sutg{n}{i}\cap\ker \Phi^{n}_{n+k}=\im\bigg({\rm Proj}^i_n\circ G_n\bigg)$$
Since $I$ is an isomorphism, we can assume that
$$x=\sum_{j=1}^l{\rm Proj}^i_n\circ G_n(\lambda_j\cdot I(\al_j))$$
where $\al_j\in\mathfrak{B}$ and ${\rm Proj}^i_n\circ G_n\circ I(\al_j)\neq0$. From Lemma \ref{lem: the isomorphism I} part (3), we know that this implies that for any $j\in[1,l]\cap\intg$, we have $n\leq \tau(\al_j)-2$ and there exists $N_j\in\intg$ such that $N_j\in[0,\tau(\al_j)-2-n]$
$$i=\frac{\tau^+(\al_j)+\tau^-(\al_j)-(\tau(\al_j)-2-n)p}{2}+N_jp.$$
Now, for $k$ sufficiently large, we have $n+2k>\tau(\al_j)$. Taking
$$N_j'=n+k+1-\tau(\al_j)+N_j\in\intg,$$ 
it is straightforward to verify that when $k$ is sufficiently large, we have
$$N_j'\in[0,n+2k-\tau(\al_j)]~{\rm and~}~i=\frac{\tau^+(\al_j)+\tau^-(\al_j)-(n+2k-\tau(\al_j))p}{2}+N_j'p.$$
Hence by Lemma \ref{lem: basic properties of tau} part (3), there exists $y_j\in\sutg{n+2k}{i}$ such that $F_{n+2k}(y_j)=\al_j$. As a result, it is straightforward to check that
$$x=\Phi^{n+2k}_n\bigg(\sum_{j=1}^l\lambda_j\cdot y_j\bigg)\in\im\Phi^{n+2k}_n.$$

{\bf The exactness at $\sut{n+2k}$}. Suppose $x\in\sutg{n+2k}{i}$ and
$$F_{n+2k}(x)=\sum_{j=1}^l\lambda_j\cdot \al_j$$
with $\lambda_j\neq 0$ and $\al_j\in \mathfrak{B}$.

First, if $x\in \im\Phi^{n+k}_{n+2k}$, then from Lemma \ref{lem: image of Phi^n+k_n+2k} part (2), we know that for any $1\leq j\leq l$, we have
$${\rm either}~i\geq \tau^-(\al_j)-\frac{(n-1)p-q}{2}~{\rm or~}i\leq \tau^+(\al_j)+\frac{(n-1)p-q}{2}.$$
If we write
$$i=\frac{\tau^+(\al_j)+\tau^-(\al_j)-(\tau(\al_j)-2-n)p}{2}+N_jp,$$
for some $N_j$ then the inequality
$$i\geq \tau^-(\al_j)-\frac{(n-1)p-q}{2}$$
implies that
\beq
N_j&\geq \frac{1}{p}\bigg(\tau^-(\al_j)-\frac{(n-1)p-q}{2}-\frac{\tau^+(\al_j)+\tau^-(\al_j)-(\tau(\al_j)-2-n)p}{2}\bigg)\\
&=\frac{1}{2}\bigg(1+\frac{\tau^-(\al_j)-\tau^+(\al_j)+q}{p}\bigg)+\frac{\tau(\al_j)}{2}-1-n\\
&=\tau(\al_j)-1-n.
\eeq
Note that the last equality uses the definition of $\tau(\al)$ in Definition \ref{defn: tau inv}.
Similarly, we can compute that
$$i\leq \tau^+(\al_j)+\frac{(n-1)p-q}{2}$$
implies that
\beq
N_j&\leq \frac{1}{p}\bigg(\tau^+(\al_j)+\frac{(n-1)p-q}{2}-\frac{\tau^+(\al_j)+\tau^-(\al_j)-(\tau(\al_j)-2-n)p}{2}\bigg)\\
&=\frac{1}{2}\bigg(-1+\frac{\tau^+(\al_j)-\tau^-(\al_j)-q}{p}\bigg)+\frac{\tau(\al_j)}{2}-1\\
&=-1.
\eeq
In summary, $x\in\im\Phi^{n+k}_{n+2k}$ implies that for all $1\leq j\leq l$, either $N_j\geq \tau(\al_j)-1-n$ or $N_j\leq -1$. Hence from Lemma \ref{lem: the isomorphism I} part (3), we know that
$${\rm Proj}^i_n\circ G_n\circ I(\al_j)=0$$
for all $1\leq j\leq l$ and as a result, $\Phi^{n+2k}_n(x)=0$.

Second, suppose $x\notin\im\Phi^{n+k}_{n+2k}$. For any $1\leq j\leq l$, we can write
$$i=\frac{\tau^+(\al_j)+\tau^-(\al_j)-(\tau(\al_j)-2-n)p}{2}+N_jp$$
for some $N_j$. Then from Lemma \ref{lem: image of Phi^n+k_n+2k} part (3) we know that there exists $j$ such that $1\leq j\leq l$, and
$$N_j\in[0,\tau(\al_j)-2-n]\cap\intg.$$ 
Hence by Lemma \ref{lem: the isomorphism I} part (3) and (4) we know that
$${\rm Proj}^i_n\circ G_n\circ I (\al_j)\neq0\Rightarrow\Phi^{n+2k}_{n}(x)\neq 0.$$
Hence we conclude that
$$\im\Phi^{n+k}_{n+2k}=\ker\Phi^{n+2k}_n.$$
\epf

\subsection{The commutative diagram}
\label{subsec: commutative diagram}In this subsection, we will prove the commutative diagram presented at the beginning of the section. Note that we choose the basis  $\mathfrak{B}$ of $\dehny{}$ as in Section \ref{subsec: basis}.

\blem\label{lem: decomposition along basis}
Suppose $n\in\intg$ and $i$ is a grading. Suppose $x\in\sutg{n}{i}$ such that
$$F_n(x)=\sum_j^{l}\lambda_j\al_j,$$
with $\lambda_j\neq0$ and $\al_j\in\mathfrak{B}$ for all $1\leq j\leq l$. Then for any $1\leq j\leq l$, there exists $N_j\in[0,n+1-\tau(\al_j)]$ such that
$$i=\frac{\tau^+(\al_j)+\tau^-(\al_j)-(n+1-\tau(\al_j))p}{2}+N_jp.$$
\elem

\bpf
This is a combination of Lemma \ref{lem: homogeneous elements} part (3), Lemma \ref{lem: basic properties of tau} part (3), and Lemma \ref{lem: image of Phi^n+k_n+2k} part (1). The proof is similar to that of Lemma \ref{lem: image of Phi^n+k_n+2k} part (2).
\epf

\bpf[Proof of Proposition \ref{prop: commutative diagram involving Phi^n+2k_n}]
We only prove the first commutative diagram
\begin{equation*}
\xymatrix{
\sut{\frac{2n+2k+1}{2}}\ar[dd]^{\Psp{n+k+1}{n+2k}\circ \psm{\frac{2n+2k+1}{2}}{n+k+1}}\ar[rrr]^{\psp{n+k+1}{\mu}\circ \psm{\frac{2n+2k+1}{2}}{n+k+1}}&&&\sut{\mu}\ar[dd]^{\psp{\mu}{n}}\\
&&&\\
\sut{n+2k}\ar[rrr]^{\Phi^{n+2k}_{n}}&&&\sut{n}
}	
\end{equation*}
The other is similar. Note that at the end of Section \ref{subsec: basis}, we introduce new notations of $\eta_{\pm,n-2}^n$ to remove the scalars. Then the second commutative diagram only holds up to a scalar.

First, note that the maps from $\sut{\frac{2n+2k+1}{2}}$ to $\sut{\mu}$ and $\sut{n+2k}$ both factor through $\sut{n+k+1}$. As a result, we only need to prove the following commutative diagram for sufficiently large $k$.
\begin{equation*}
\xymatrix{
\sut{n+k+1}\ar[dd]^{\Psp{n+k+1}{n+2k}}\ar[rrr]^{\psp{n+k+1}{\mu}}&&&\sut{\mu}\ar[dd]^{\psp{\mu}{n}}\\
&&&\\
\sut{n+2k}\ar[rrr]^{\Phi^{n+2k}_{n}}&&&\sut{n}
}	
\end{equation*}
Now suppose $x\in\sutg{n+k+1}{i}$. Write
$$F_{n+k+1}(x)=\sum_{j=1}^l\lambda_i\cdot\al_j$$
with $\lambda_j\neq0$ and $\al_j\in\mathfrak{B}$ for $1\leq j\leq l$. We want to first establish an identity
\begin{equation}\label{eq: lem 7.9, 1}
	\psp{\mu}{n}\circ\psp{n+k+1}{\mu}(x)=\sum_{\substack{1\leq j\leq l\\n\leq \tau(\al_j)-2\\ N_j=n+k+1-\tau(\al_j)}}\lambda_j\cdot\eta^{\tau(\al_j)}_{+,n}(\al_j)
\end{equation}
and then show that the other composition has exactly the same expression.

From Lemma \ref{lem: decomposition along basis}, we know for any $1\leq j\leq l$, there exists $N_j\in[0,n+k+1-\tau(\al_j)]$ such that
\begin{equation}\label{eq: commutativity, N_j}
	i=\frac{\tau^+(\al_j)+\tau^-(\al_j)-(n+k+1-\tau(\al_j))p}{2}+N_jp.
\end{equation}
Taking $n'_j=\tau(\al_j)$ and $N_j'=0$, we can apply Lemma \ref{lem: basic properties of tau} part (3) to find an element
$$x_j\in \sutg{\tau(\al_j)}{\frac{\tau^{+}(\al_j)+\tau^{-}(\al_j)}{2}}$$
such that
$$F_{\tau(\al_j)}(x_j)=\al_j.$$
It is then straightforward to check that
\begin{equation}\label{eq: commutativity, hat x_j}
	y_j=\Psp{\tau(\al_j)+N_j}{n+k+1}\circ\Psm{\tau(\al_j)}{\tau(\al_j)+N_j}(x_j)\in\sutg{n+k+1}{i}.
\end{equation}
Write
$$y=x-\sum_{j=1}^l\lambda_j\cdot y_j\in\sutg{n+k+1}{i}.$$
From Lemma \ref{lem: comm diag for n,n+1,dehn} we know that
$$F_{n+k+1}(y)=0.$$

As a result, by Lemma \ref{lem: ker of F_n^i},
\beq
\psp{\mu}{n}\circ\psp{n+k+1}{\mu}(x)&=\sum_{j=1}^l\lambda_j\cdot\psp{\mu}{n}\circ\psp{n+k+1}{\mu}(y_j)\\
\eeq
Note that, unless $N_j=n+k+1-\tau(\al_j)$, we have
$$\psp{n+k+1}{\mu}\circ\Psp{\tau(\al_j)+N_j}{n+k+1}=0$$
by the exactness. As a result,
\beq
\psp{\mu}{n}\circ\psp{n+k+1}{\mu}(x)&=\sum_{\substack{1\leq j\leq l\\ N_j=n+k+1-\tau(\al_j)}}\lambda_j\cdot\psp{\mu}{n}\circ\psp{n+k+1}{\mu}\circ\Psm{\tau(\al_j)}{\tau(\al_j)+N_j}(x_j)\\
({\rm Lemma~\ref{lem: comm diag for n,n+1,mu}})&=\sum_{\substack{1\leq j\leq l\\n\leq \tau(\al_j)-2\\ N_j=n+k+1-\tau(\al_j)}}\lambda_j\cdot\psp{\mu}{n}\circ\psp{\tau(\al_j)}{\mu}(x_j)+\sum_{\substack{1\leq j\leq l\\n\geq \tau(\al_j)-1\\ N_j=n+k+1-\tau(\al_j)}}\lambda_j\cdot\psp{\mu}{n}\circ\psp{\tau(\al_j)}{\mu}(x_j)\\
\text{Equation } (\ref{eq: prop 7.2, 3})&=\sum_{\substack{1\leq j\leq l\\n\leq \tau(\al_j)-2\\ N_j=n+k+1-\tau(\al_j)}}\lambda_j\cdot\psp{\mu}{n}\circ\psp{\tau(\al_j)}{\mu}(x_j)\\
({\rm Definition~of~}\eta^{\tau(\al_j)}_{+,n})&=\sum_{\substack{1\leq j\leq l\\n\leq \tau(\al_j)-2\\ N_j=n+k+1-\tau(\al_j)}}\lambda_j\cdot\eta^{\tau(\al_j)}_{+,n}(\al_j)
\eeq
This verifies Equation (\ref{eq: lem 7.9, 1}) if we show that
\[
\sum_{\substack{1\leq j\leq l\\n\geq \tau(\al_j)-1\\ N_j=n+k+1-\tau(\al_j)}}\lambda_j\cdot\psp{\mu}{n}\circ\psp{\tau(\al_j)}{\mu}(x_j)=0.
\]
To verify this last equality, assume that $n\geq \tau(\al_j)-1$. Then from Lemma \ref{lem: comm diag for n,n+1,mu} and the exactness of the bypass maps we have
\begin{equation}\label{eq: prop 7.2, 3}
	\psp{\mu}{n}\circ\psp{\tau(\al_j)}{\mu}(x_j)=\psp{\mu}{n}\circ\psp{n+1}{\mu}\circ\Psm{\tau(\al_j)}{n+1}(x)=0.
\end{equation}

Now we deal with $\Phi^{n+2k}_{n}\circ\Psp{n+k+1}{n+2k}(x)$. Since $F_{n+k+1}(y)=0$, Lemma \ref{lem: F_n and G_n are iso when n large} implies that
$$\Psp{n+k+1}{n+2k}(y)=0.$$
Hence
$$\Psp{n+k+1}{n+2k}(x)=\sum_{j=1}^l\lambda_j\cdot\Psp{n+k+1}{n+2k}(y_j)$$
where $y_j$ is defined as in (\ref{eq: commutativity, hat x_j}). Note that, by definition we have $y_j\in \sutg{n+1+k}{i}$, so from Lemma \ref{lem: bypass n,n+1,mu}, we know
$$\Psp{n+k+1}{n+2k}(y_j)\in\sutg{n+2k}{i-\frac{(k-1)p}{2}}.$$
Note that, by (\ref{eq: commutativity, hat x_j}) and Lemma \ref{lem: comm diag for n,n+1,dehn}, we know that
$$F_{n+2k}\circ \Psp{n+k+1}{n+2k}(y_j)=F_{\tau(\al_j)}(x_j)=\al_j.$$
Hence
$$\Phi^{n+2k}_{n}\circ\Psp{n+k+1}{n+2k}(x)=\sum_{j=1}^l\lambda_j\cdot{\rm Proj}^{i-\frac{(k-1)p}{2}}_n\circ G_n\circ I(\al_j).$$
We write
$$i-\frac{(k-1)p}{2}=\frac{\tau^+(\al_j)+\tau^-(\al_j)-(\tau(\al_j)-2-n)p}{2}+N_j'p$$
Comparing the above formula with (\ref{eq: commutativity, N_j}), we know
$$N_j'=N_j+\tau(\al_j)-n-k-1.$$
Note that, by construction, $N_j\leq n+k+1-\tau(\al_j)$, which means $N_j'\leq 0$. Hence from Lemma \ref{lem: the isomorphism I} we know
$${\rm Proj}^{i-\frac{(k-1)p}{2}}_n\circ G_n\circ I(\al_j)\neq0$$
if and only if $N_j'=0$, \textit{i.e.}, $N_j=n+k+1-\tau(\al)$.
 Also when $N_j'=0$ from Lemma \ref{lem: the isomorphism I} part (5) we know
$${\rm Proj}^{i-\frac{(k-1)p}{2}}_n\circ G_n\circ I(\al_j)=\eta^{\tau(\al_j)}_{+,n}(\al_j).$$

Note that we could focus on indices $j$ such that ${\rm Proj}^{i-\frac{(k-1)p}{2}}_n\circ G_n\circ I(\al_j)\neq0.$ This is because if an index $j$ makes ${\rm Proj}^{i-\frac{(k-1)p}{2}}_n\circ G_n\circ I(\al_j)=0$, then on one hand it does not contribute to $\Phi^{n+2k}_{n}\circ\Psp{n+k+1}{n+2k}(x)$ since the corresponding summand is $0$, on the other hand, we have $N_j\neq n+k+1-\tau(\al)$ hence per Equation (\ref{eq: lem 7.9, 1}) it does not contribute to $\psp{\mu}{n}\circ\psp{n+k+1}{\mu}(x)$, either. Also, we know from Lemma \ref{lem: the isomorphism I} part (3) that when ${\rm Proj}^{i-\frac{(k-1)p}{2}}_n\circ G_n\circ I(\al_j)\neq0$ we must have $n\leq \tau(\al_j)-2$. As a result, we know
\beq
\Phi^{n+2k}_{n}\circ\Psp{n+k+1}{n+2k}(x)&=\sum_{j=1}^l\lambda_j\cdot{\rm Proj}^{j-\frac{(k-1)p}{2}}_n\circ G_n\circ I(\al_j)\\
&=\sum_{\substack{1\leq j\leq l\\n\leq\tau(\al_j)-2\\ N_j=n+k+1-\tau(\al_j)}}\lambda_j\cdot\eta^{\tau(\al_j)}_{+,n}(\al_j)\\
\text{Equation }(\ref{eq: lem 7.9, 1})&=\psp{\mu}{n}\circ\psp{n+k+1}{\mu}(x)
\eeq

\epf


\bibliography{ref.bib}

\end{document}